\documentclass[12pt]{amsart}

\usepackage{graphicx}

\newcommand{\R}{{\mathbb R}}

\newcommand{\Z}{{\mathbb Z}}
\newcommand{\SL}{{\mathcal S} \ell}
\newcommand{\DL}{{\mathcal D} \ell}

\newcommand{\cg}{\chi_{\partial \gamma}}

\newtheorem{theo}{{\sc Theorem}}[section]
\newtheorem{cor}[theo]{{\sc Corollary}}

\newtheorem{lem}[theo]{{\sc Lemma}}

\newtheorem{prop}[theo]{{\sc Proposition}}
\newtheorem{defn}[theo]{{\sc Definition}}
\newtheorem{rem}[theo]{{\sc Remark}}

\title[Inverse spectral problem I ]{Inverse spectral problem for analytic domains I: \newline
Balian-Bloch trace formula }

\author{Steve Zelditch}
\address{Department of Mathematics, Johns Hopkins University, Baltimore, MD
21218, USA}
\email{
zelditch@@math.jhu.edu }

\thanks{Research partially supported by  NSF grant \#DMS-0071358.}

\date{\today}

\begin{document}

\begin{abstract} This is the first in a series of papers \cite{Z3, Z4} on  inverse spectral/resonance problems
for analytic plane domains $\Omega$.  In this paper, we present a
rigorous version of the Balian-Bloch trace formula  \cite{BB1,
BB2}. It is an asymptotic   formula for the trace $Tr 1_{\Omega}
R_{\rho}(k + i \tau \log k)$ of the regularized resolvent of the
Dirichlet Laplacian of $\Omega$ as $k \to \infty$ with $\tau > 0$.
When the support of $\hat{\rho}$ contains  the length $L_{\gamma}$
of precisely one  periodic reflecting ray $\gamma$, then the
asymptotic expansion  of $Tr 1_{\Omega} R_{\rho}(k + i \tau \log
k)$ is essentially the same as the wave trace expansion at
$\gamma$.  The raison d'etre for this approach to wave invariants
is that they are explicitly computable.  Applications of the trace
formula will be given in the subsequent articles in this series.
For instance, in \cite{Z3, Z4} we will prove that analytic domains
with one symmetry are determined by their Dirichlet (or Neumann)
spectra. Although we only present details in dimension 2 for the
sake of simplicity, the methods and results extend with few
modifications to all dimensions.
\end{abstract}

\maketitle
\section{Introduction}  This paper is the first in a series of articles devoted to
  inverse spectral and resonance problems for
analytic domains $\Omega \subset \R^2$ with Dirichlet (or Neumann)
boundary conditions \cite{Z3, Z4}.  As in the earlier articles
\cite{Z1, Z2, ISZ}, the aim is to recover the domain $\Omega$ from
spectral invariants associated to special closed orbits of the
billiard flow, specifically  the wave trace invariants
$B_{\gamma^r, j}$ at iterates $\gamma^r$ of  a single bouncing
ball orbit $\gamma$. Such wave invariants are polynomials in the
Taylor coefficients $f_{\pm}^{(j)}(0)$ of the defining functions
$f_{\pm}$ of $\Omega$ at the endpoints of $\gamma$ and our goal is
to recover the coefficients $f_{\pm}^{(j)}(0)$ from the invariants
$B_{\gamma^r, j}$ as $r$ varies, and hence to recover the analytic
domain.  In this series, we will  achieve this goal  when $\Omega$
satisfies one symmetry condition.  Except at the very last
inverting  step of \cite{Z3, Z4}, we do not use any symmetry or
analyticity assumptions on $\Omega$ and our computations of wave
invariants are valid on all smooth plane domains. We restrict to
plane domains to simplify the exposition;  the  methods extend in
a straightforward way to domains in $\R^n$.

The path we take towards the
 wave invariants $B_{\gamma^r, j}$ is the  one  initiated by Balian-Bloch in the classic papers
 \cite{BB1, BB2}. As recalled below, these papers are concerned with the asymptotics of a (regularized) resolvent trace  rather than the trace of
the wave group, and are based on the
 Neumann series representation for the resolvent kernel in terms of the free resolvent.    Better known to physicists than to mathematicians, the Balian-Bloch papers were one of  the origins
of the Poisson relation for manifolds with boundary (see below and
\cite{CV2}).  The purpose of this first  paper in the series is to
give a rigorous version of  the Balian-Bloch approach. In
\cite{Z3}, we  use it to calculate wave invariants in terms of
Feynman diagrams and amplitudes and to  reduce the inverse
spectral problem to concrete combinatorial problems.

In  subsequent articles, we will  demonstrate the usefulness of
this reduction with the following applications:
\begin{itemize}

\item In  \cite{Z3} we  prove that simply connected analytic
domains with one symmetry (that reverses a bouncing ball orbit)
are determined among such domains by their Dirichlet spectra;

 \item In \cite{Z4} we  prove that exterior domains with one
mirror symmetry are determined by their resonance poles. This is a
resonance analogue of our inverse spectral results;

\item These results should admit extensions  to  domains formed by
flipping the graph of an analytic function around the $x$-axis,
allowing for corners at the $x$-intercept. It is easy to construct
examples of this kind.

\end{itemize}

The crucial advantage  of our present approach is the very
explicit nature of the trace formula, which allows us to obtain
relatively simple formula for lower order terms in wave trace
expansions. These lower order terms allow us to remove one
symmetry from the inverse results of ourself and of
Iantchenko-Sj\"ostrand-Zworski \cite{Z1, Z2, ISZ}. The rather
concrete formulae for wave invariants (see \ref{WTF}) are valid
for all plane domains and the symmetry assumption is only used to
reduce the amount of data required to determine the domain.

 We should emphasize that we view the Balian-Bloch asymptotic expansions as primarily a computational device. Although
one could give a new and self-contained proof of the Poisson
relation for bounded Euclidean domains by this approach, we do not
do so here. In particular, we use the Melrose-Sjostrand results on
propagation of singularities to microlocalize traces to orbits.
Also, we do not explain how to cancel so-called `ghost orbits' of
non-convex domains, i.e. closed billiard trajectories which do not
stay within the domain. To explain how the Balian-Bloch approach
improves on normal forms or microlocal parametrix constructions,
and to explain its connection to the usual Poisson relation, we
now give an informal
 exposition of the Balian-Bloch trace formula.

\subsection{Introduction to wave trace and  Balian-Bloch}  Let $\Omega$ denote a bounded $C^{\infty}$ plane domain, and  let
$$E_{\Omega} (t) = \cos t \sqrt{\Delta_{\Omega}}$$
denote the even part of the wave group of the Dirichlet  Laplacian
$\Delta_{\Omega}$.  As recalled in \S 3,  the singular points $t$
of the distribution trace $Tr E_{\Omega}(t)$ are contained in the
length spectrum $Lsp(\Omega)$ of $\Omega,$ i.e. the set of
 lengths $t = L_{\gamma}$ of closed orbits  of the billiard flow
$\Phi^t$ of $\Omega$ (i.e. straight-line motion in $\Omega$ with the Snell law of reflection at
the boundary; cf. \S 2). In particular, at lengths $L_{\gamma}$ of  periodic reflecting rays $\gamma$ the trace has a complete
singularity expansion; its  coefficients are known as the wave trace invariants of $\gamma$.

The wave trace invariants at a periodic reflecting ray may be obtained from a dual semiclassical asymptotic
expansions as $k \to \infty$ for the regularized trace of  the Dirichlet resolvent
$$R_{\Omega}(k + i \tau) : = - (\Delta_{\Omega} + (k + i \tau)^2)^{-1},\;\;\;\; (k \in \R, \tau \in \R^+). $$
We emphasize that  $\Delta_{\Omega}$ (the usual Dirichlet
Laplacian) is a negative operator, so the signs of the two terms
are opposite. In the classical work of Seeley and others on
resolvent traces, asymptotics are taken along the vertical line
(or a ray of non-zero slope) in the upper half plane $k + i \tau$
and traces are polyhomogeneous functions of $k$. Here, we are
taking asymptotics along horizontal lines (or  logarithmic curves)
in this half-plane and obtain oscillatory asymptotics reflecting
the behaviour of closed geodesics.

 We fix
 a non-degenerate periodic reflecting ray  $\gamma$ (cf. \S 2).
 For technical convenience we will assume that the reflection
 points of $\gamma$ are points of non-zero curvature of $\partial
 \Omega$. This is a minor assumption, but it simplifies one
 argument (Proposition \ref{COMP}) and we use it again in
 \cite{Z3}. We briefly explain there how to remove this assumption. We also
  let  $\hat{\rho} \in C_0^{\infty}(L_{\gamma} - \epsilon, L_{\gamma} + \epsilon)$ be a cutoff,  equal to one on an interval  $(L_{\gamma} - \epsilon/2, L_{\gamma} + \epsilon/2)$ which contains
no other lengths in Lsp$(\Omega)$ occur in its support.
We then define the  smoothed (and localized) resolvent by
\begin{equation}\label{RRHO} R_{\rho}(k + i \tau):= \int_{\R} \rho (k - \mu) (\mu + i \tau)  R_{\Omega}(\mu + i \tau) d \mu.  \end{equation}
When $\gamma, \gamma^{-1}$ are the unique closed orbits of length
$ L_{\gamma}$,  it follows from the Poisson relation for manifolds
with boundary (see \S 3 and \cite{GM, PS}; see also  Proposition
(\ref{BBL}))) that the trace $Tr 1_{\Omega} R_{\rho}((k + i \tau))$ of the regularized
resolvent on $L^2(\Omega)$  admits a
complete asymptotic expansion  of the form:
\begin{equation} \label{PR} Tr 1_{\Omega} R_{\rho}(k + i \tau) \sim  e^{ (i k - \tau) L_{\gamma}}  \sum_{j = 1}^{\infty} (B_{\gamma, j}
+  B_{\gamma^{-1}, j}) k^{-j},\;\;\; k \to \infty \end{equation}
with coefficients $B_{\gamma, j},  B_{\gamma^{-1}, j}$ determined by the jet of $\Omega$ at the reflection points of $\gamma$, in the sense that
$$Tr 1_{\Omega} R_{\rho}(k + i \tau) - e^{ (ik - \tau) L_{\gamma}}  \sum_{j = 1}^{R} (B_{\gamma, j} +  B_{\gamma^{-1}, j}) k^{-j} = {\mathcal O} (|k|^{-R}).$$
The coefficients $B_{\gamma, j},  B_{\gamma^{-1}, j}$ are thus essentially the same as the wave trace
 coeffficients at the singularity  $t = L_{\gamma}$. Our main goal in this paper is to give a useful algorithm for
calculating them explicitly in terms of the defining function of $\partial \Omega.$

In fact,  it is technically more convenient to consider
asymptotics of traces along logarithmic curves $k + i \tau \log k$
in the upper half plane. We therefore modify the regularization
(\ref{RRHO}) to
\begin{equation}\label{RRHOLOG} R_{\rho}(k + i \tau \log k ):= \int_{\R} \rho (k - \mu) (\mu + i \tau \log k)  R_{\Omega}(\mu + i \tau \log \mu) d \mu.  \end{equation}
In place of
(\ref{PR}) we will get
\begin{equation} \label{PRLOG} Tr 1_{\Omega} R_{\rho}(k + i \tau \log k ) \sim  e^{ (i k ) L_{\gamma}} k^{- \tau L_{\gamma}}  \sum_{j = 1}^{\infty} (B_{\gamma, j}
+  B_{\gamma^{-1}, j}) k^{-j},\;\;\; k \to \infty. \end{equation}
The additional power law decay $k^{- \tau L_{\gamma}}$ will not cause  problems in our study
of wave invariants at a fixed closed geodesic, because the
the  errors have the accuracy of $k^{- {\infty}}$.

In principle, one could obtain sufficiently explicit formulae for
$B_{\gamma, j}$ by
 applying the method of stationary phase to a microlocal parametrix at $\gamma$ (cf. \cite{GM} \cite{PS}),  or
by constructing a Birkhoff normal form for $\Delta$ at $\gamma$
\cite{Z1, Z2, SZ, ISZ}.  However, in practice we have found the
Balian-Bloch approach more effective. Its starting point is  to
write the Dirichlet Green's kernel $G_{\Omega}(k + i \tau, x, y)$
of $R_{\Omega}(k + i \tau)$ as a Neumann series (called the
multiple reflection expansion in \cite{BB1})  in terms of the free
Green's function $G_0(k + i \tau, x, y)$, i.e. the kernel of the
free resolvent $R_{\Omega}(k + i \tau) = - (\Delta_0 + (k + i
\tau)^2)^{-1}$ on $\R^2$:
\begin{equation}\label{MREone} \begin{array}{lll} G_{\Omega}(k + i \tau,x, y) &=&
G_0(k + i \tau, x, y) +  \sum_{M = 1}^{\infty} (-2)^M G_M(k + i \tau, x, y), \;\mbox{where}\\ & & \\
G_M(k + i \tau ,x,y )  & =  &  \int_{(\partial \Omega)^M} \partial_{\nu_y} G_0(k + i \tau, x, q_1)G_0(k + i \tau, q_M, y)  \\ & & \\
& \times &\Pi_{j = 1}^{M - 1}
\partial_{\nu_y} G_0(k + i \tau, q_{j + 1}, q_j) ds(q_1) \cdots ds(q_M),\end{array} \end{equation}
where $ds(q)$ denotes arclength on $\partial \Omega$ and
$\partial_{\nu_y}$ is the interior unit normal operating in the
second variable. (The only change in the case of the Neumann
Green's kernel is that the $M$th term has an additional factor of
$(-1)^M$; by making this change, all our methods and results
extend immediately to the Neumann case.)  The terms are
regularized as in (\ref{RRHO})  by setting
\begin{equation} \label{RGM} G_{M, \rho}(k + i \tau \log k)) = \int_{\R} \rho (k - \mu)  (\mu + i \tau \log \mu ) G_{M}(\mu + i \tau \log \mu) d\mu \end{equation}
There exists an  explicit formula for $G_0(k + i \tau)$ in terms of Hankel functions (\S 4), from which it appears that  the traces $Tr 1_{\Omega} G_{M, \rho}(k + i \tau \log k) = \int_{\Omega} G_{M, \rho}(k + i \tau \log k, x, x) dx$ are (formally)  oscillatory integrals with phases
$$L(x, q_1, \dots, q_M, x) = |x - q_1| + |q_1 - q_2| + \cdots + |q_M - x|,$$
equal to the length of the polygon with vertices at points $(q_1, \dots, q_M) \in (\partial \Omega)^M$ and $x \in \Omega.$ The smooth critical points correspond to the $M$-link periodic reflecting rays of the billiard
flow of $\Omega$ of length $L_{\gamma}$, satisfying Snell's law at each vertex (for short, we call such
polygons Snell polygons).   Since the amplitudes and phases only involve the free Green's function, they
are known explicitly and  the coefficients of the  stationary phase expansion of Tr $1_{\Omega} R_{\rho}((k + i \tau \log k))$ can be calculated
explicitly.

Thus, our plan for determining $\Omega$ from wave invariants at a bouncing ball orbit  $\gamma$ is as follows: We chose $\rho$ to localize at the length of the $r$th interate of $\gamma$.   We  then localize the integrals $G_{\rho, M}(k + i \tau \log k )$ to small intervals around
the endpoints $\{(0, -L/2), (0, L/2)\}$ and parametrize the two components, as above, by $y = f_{-}(x),$ resp. $y =  f_+(x).$
(see figure \ref{1}) The integral over $(\partial \Omega)^M$ is then reduced to an  integral over $(-\epsilon, \epsilon)^M$ with phase
and amplitude given by canonical functions of $f_-(x), f_+(x).$  Applying the stationary phase method, we obtain
coefficients which are polynomials  in the data
$f_+^{(j)}(0), f_-^{(j)}(0)$.   By examining the behaviour of the coefficients under iterates $\gamma^r$ of
the bouncing ball orbit, we try to  determine all Taylor coefficients $f_{\pm}^{(j)}(0)$ from these stationary phase coefficients.  In \cite{Z3, Z4, Z5}, we will see that this method suceeds at least  if the domain has one symmetry (so
that $f_+ = - f_-$).


A number   of technical problems must be overcome  in turning this
idea into a proof. The problem is that the free Green's kernel
(and its normal derivative) only possesses a WKB formula away from
the diagonals $q_i = q_{i + 1}$, and worse, it is  singular along
these diagonals. Hence, in the Balian-Bloch approach (unlike the
parametrix approach), the following problems arise:
\begin{itemize}

\item The individual  terms $G_{M, \rho}$ must be regularized, i.e. converted into sums
of standard oscillatory integrals. Until then, it is not even clear  that the trace of  each term
has asymptotic an asymptotic expansion, or that it  localizes at critical points;

\item There are non-smooth critical points (corresponding to Snell polygons in which at least one edge collapsed),
and  their contribution  to the stationary phase expansion must be determined;

\item Since the multiple reflection expansion (\ref{MREone}) is an infinite series, it must be
explained how the full series produces an asymptotic expansion (\ref{PR})  of $Tr 1_{\Omega} R_{\rho}(k + i \tau \log k )$, and why each term
in (\ref{PR})   depends on only the trace of  a finite number of terms of (\ref{MREone}).
In particular, the `tail' trace must be estimated.

\end{itemize}

Let us briefly outline how we deal with these difficulties. At the
same time, we will justify the effort involved by explaining how
the new results are obtained by this method.

\subsection{Regularizing the terms of the  Neumann series}

We first explain how we  regularize the individual  terms in the
Neumann series. For expository reasons, we  suppress the role of
the interior variable $x$ in the introduction.

As will be explained in \S \ref{LP}, the multiple reflection
expansion is derived from an exact formula  of potential theory
(\ref{NS}) for the Dirichlet resolvent:
\begin{equation} \label{POT}  R_{\Omega}(k + i \tau) =
 R_0(k + i \tau) - {\mathcal D} \ell(k + i \tau) (I +  N(k + i \tau))^{-1} r_{\Omega}
 {\mathcal S} \ell^{tr}(k + i \tau), \end{equation}
by expanding the operator $(I +  N(k + i \tau))^{-1} $ in a
geometric series. Here,
 ${\mathcal D} \ell (k + i \tau)$ (resp. ${\mathcal S} \ell (k + i \tau)$) is the double
 (resp. single) layer potential (see (\ref{layers})), $ {\mathcal S}^{tr}(k + i \tau)$ is the
 transpose,   and
$N(k + i \tau)$ is the boundary integral operator  on
$L^2(\partial \Omega)$ induced by  ${\mathcal D} \ell (k + i
\tau)$ (see (\ref{blayers}). Also, $R_0(k + i \tau)$ is the free
resolvent on $\R^2,$ and $r_{\Omega}$ is the restriction to the
boundary. The existence of the inverse $(I + N(k + i \tau))^{-1} $
is guaranteed by Fredholm theory.

As will be discussed in detail in \S \ref{LP} (see also
\cite{HZ}), the operator $N(k + i \tau)$ is a hybrid Fourier
integral operator. It has the singularity of a homogeneous
pseudodifferential operator of order $-1$ on the diagonal (in
fact, it is of order $-2$ in dimension $2$, see Proposition
\ref{N}). This is the way it is normally described in potential
theory \cite{T}. However, away from the diagonal, it has a WKB
approximation which exhibits it as a semi-classical Fourier
integral operator with phase $d_{\partial \Omega}(q,q') = |q -
q'|$ on $\partial \Omega \times \partial \Omega$, the boundary
distance function of $\Omega$. To make this precise,  we introduce
a cutoff $\chi(k^{1 - \delta} |q - q'| )$ to the diagonal, where
$\delta > 1/2$ and where $\chi \in C_0^{\infty}(\R)$ is a cutoff
to a neighborhood of $0$. We then put
\begin{equation} \label{N01}  N(k + i \tau) = N_0(k + i \tau) + N_1(k + i
\tau), \;\; \mbox{with}  \end{equation}
\begin{equation} \label{N01DEF} \left\{\begin{array}{l}  N_0(k + i
\tau, q, q') = \chi(k^{1 - \delta} |q - q'| ) \;  N(k + i \tau, q,
q'), \\ \\ N_1(k + i \tau, q, q') = (1 - \chi(k^{1 - \delta} |q -
q'| ))\; N(k + i \tau, q, q'). \end{array} \right.
\end{equation}
As will be shown in Proposition \ref{NSYM},  $ N_1((k + i \tau),
q, q') ) $ is a semiclassical Fourier integral operator with phase
equal to $d_{\partial \Omega}(q,q')$.

Roughly speaking, the boundary distance function generates the
billiard map of $\partial \Omega$ (see \S \ref{BILL} for the
definition). This is literally
 correct only on convex domains, since $d_{\partial \Omega}(q,q')$
generates both the interior and exterior billiard map, and hence
on non-convex domains its canonical relation contains `ghost
orbits' (orbits which may exit and re-enter the domain). Because
we are microlocalizing to one periodic reflecting ray $\gamma$,
ghost orbits play no essential role in this paper and we think of
$N_1(k + i \tau)$ as `quantizing' the billiard map.

Now consider the powers $N(k + i \tau)^M$ which arise when
expanding $(I - N(k + i \tau))^{-1}$ in a geometric series. We
write
\begin{equation} \label{BINOMIAL} (N_0 + N_1)^M = \sum_{\sigma: \{ 1, \dots, M\} \to \{0, 1\}}
N_{\sigma(1)} \circ N_{\sigma(2)} \circ \cdots \circ
N_{\sigma(M)}. \end{equation} To regularize $N^M$ is essentially
to remove all of the factors of $N_0$ from each of these terms.
This is obviously not possible for the term $N_0^M$ but  it is
possible for the other terms. In Proposition \ref{NPOWERFINAL}, we
show is that $N_0 N_1$ and $N_1 N_0$ are semiclassical Fourier
integral operators of the same type as $N_1$ (and with the same
phase), but with an amplitude of one lower degree in $k$. Thus,
the term $N_1^M$ is of the highest order in the sum. We will see
that it is the only important term for the inverse spectral
problem. In Proposition \ref{COMP}, we similarly break up the
layer potentials in (\ref{POT}) and analyse compositions with
these. Further we compose with a special kind of semiclassical
cutoff operator $\chi (x, k^{-1} D_x)$ to a neighborhood of the
orbit on both sides of (\ref{POT}).

\subsection{The M-aspect}

Since the geometric series expansion of $(I +  N(k + i
\tau))^{-1}$ is very slowly converging, we write it as a partial
geometric series and a remainder:
\begin{equation}\label{GS} \begin{array}{l} (I +  N(k + i \tau))^{-1} = \sum_{M = 0}^{M_0} (-1)^M N(k + i \tau)^M \;
+\; {\mathcal R}_{M_0},\\ \\
\mbox{where}\;\; {\mathcal R}_{M_0} = N(k + i \tau)^{M_0 + 1} \; (I +  N(k + i \tau))^{-1}. \end{array} \end{equation}

The partial geometric series is regularized by the methods
described in  the previous section. We now explain how to estimate
the tail when evaluating wave invariants at a closed orbit
$\gamma$. In our applications, $\gamma$ is a bouncing ball orbit,
so we will assume it is one in the remainder estimate.

We  use elementary inequalities on traces to reduce the estimate
of the tail trace
\begin{equation} \label{TAILTRACE} Tr 1_{\Omega} {\mathcal R}_{M_0}(k + i \tau \log k) \chi (x, k^{-1} D_x) \end{equation}
to an estimate of $Tr N^{M_0}  \cg (k + i \tau \log k) \cg^*
N^{M_0 *}$, where $\cg$ is a semiclassical  cutoff operator on
$\partial \Omega$ to the periodic orbit of the billiard map
corresponding to $\gamma$.  As will be verified in \S \ref{M},
after expanding as in (\ref{BINOMIAL}), the operator $N^{M_0}  \cg
(k + i \tau \log k) \cg^* N^{M_0 *} $ can be regularized as above
as a sum over $\sigma$ of semiclassical Fourier integral operators
whose phases involve lengths ${\mathcal L}$ of $M_0 -
|\sigma|$-link billiard trajectories, where $|\sigma|$ represents
the number of $N_0$ factors in a term. The cutoff operator will
force the direction of the first link to point nearly along
$\gamma$ and hence will force each link of  the critical $M_0 -
|\sigma|$-link paths to have length $\sim  \; L_{\gamma}$. We then
use the fact that the spectral parameter varies over a logarithmic
curve: due to the imaginary part of $\tau \log k$ of the phase,
each link in the critical path gets weighted by the factor $e^{-
\tau  L_{\gamma} \log k} = k^{- \tau  }$. On the other hand, each
removal of a factor $M_0$ lowers the order in $k$ by one. Combinng
these two effects,   we show in Proposition \ref{CRUCIAL} and
Lemma \ref{MASPECT} that for $M_0$ sufficiently large, the
remainder term will be of lower order than any prescribed power
$k^{-R}$.

\subsection{Main results}

This regularization procedure provides a method of explicitly
calculating wave invariants associated to an $M$-link periodic
reflecting ray. After regularizing the initial terms of the
multiple reflection expansion and discarding the remainder, we may
apply stationary phase to rather simple and canonical oscillatory
integrals.

  Our main result may be summarized in the
following  theorem. The term `canonical' means `independent of the
domain $\Omega.$'

\begin{theo} \label{SUM} Let $\gamma$ be a primitive non-degenerate  $m$-link
 reflecting ray, whose reflection points are points of non-zero curvature of $\partial \Omega$,
 and let $\hat{\rho} \in C_0^{\infty}(\R)$ be a  cut off
 satisfying supp $\hat{\rho} \cap Lsp(\Omega) = \{ r L_{\gamma}
 \}.$
Then there exists an effective algorithm for obtaining the wave invariants at $\gamma$ consisting of the following steps:\\

\noindent (i)  Each term $G_{M, \rho}(k + i \tau
\log k)$ of (\ref{MREone}) defines a kernel  of trace class;\\

\noindent (ii) The traces $Tr 1_{\Omega} G_{M, \rho}(k + i \tau
\log k)$  can  be regularized in a canonical way as  oscillatory
integrals with  canonical amplitudes and phases;\\

\noindent (iii) The stationary phase method applies to these oscillatory integrals.
The coefficients $B_{\gamma^r, j}$ of a given order $j \leq R$ are obtained by summing the expansions
for $ Tr 1_{\Omega} G_{M, \rho}(k + i \tau \log k)$ for $M \leq M_R$. There exists
$M_0$ such that the tail trace (\ref{TAILTRACE})  is $O(k^{- R})$. \\

\noindent(iv)
The coefficients of the term $B_{\gamma^r, j}$  are universal  polynomials in the $2j + 2$-jet of the defining
function of $\partial \Omega$ at the reflection points. \\

\noindent  (v) If  $\gamma$  is  a bouncing ball orbit, then
modulo an error term  $R_{2r} (j^{2j - 2} f(0))$ depending only on
the $(2j - 2)$-jet of $f$ at the endpoints, we obtain the formula
$$ B_{\gamma^r, j - 1} + B_{\gamma^{-r}, j - 1} = a_{j, r, +}
f_+^{(2j)}(0) + a_{j, r, -} f_-^{(2j)}(0) + b_{j, r, +}
f_+^{(2j-1)}(0) + b_{j, r, -} f_-^{(2j-1)}(0),$$ where the
coefficients are polynomials in the matrix elements $h^{pq}$ of
the inverse of the Hessian of the length function $L$ at
$\gamma^r$ with universal coefficients.
\end{theo}

For instance, if $\gamma$ is elliptic and  invariant under an
isometric involution $\sigma$ of $\Omega$, then \begin{equation}
\label{WTF}  B_{\gamma^r, j - 1} + B_{\gamma^{-r}, j - 1} =  r \{
2 (h^{11})^j f^{(2j)}(0) + \{2 (h^{11})^j \frac{1}{2 - 2 \cos
\alpha/2} + (h^{11})^{j - 2} \sum_{q = 1}^{2r} (h^{1 q})^3\}
f^{(3)}(0) f^{(2j - 1)}(0)\}\}. \end{equation}  Here, $e^{\pm 2
\pi i \alpha}$ are the eigenvalues of the Poincare map
$P_{\gamma}$ of $\gamma$. To our knowledge, no such explicit formula for
the wave invariants has been available before.

The most important aspect of Theorem (\ref{SUM}) is the algorithm
for  computing the coefficients  of $B_{\gamma^r, j}$ in terms of
the defining functions of $\partial \Omega$, as illustrated in (v)
and by the explicit formula (\ref{WTF})  . That is truly what
distinguishes the result from previous approaches to wave
invariants and the Poisson relation.  For instance,  recovering a
$\Z_2$-invariant $\partial \Omega$ is tantamount to sifting the
Taylor coefficients of $f$ at $0$ from the expression (\ref{WTF}),
and we see that this amounts to a fine study of sums of powers  in
the matrix coefficients $h^{ij}$ of the inverse  Hessian of
${\mathcal L}$. We needed the complete symbol of the wave trace to
obtain this result, and no such details have previously been
obtained by any other method. The gain in effectiveness is borne
out in the subsequent articles \cite{Z3, Z4}. Although we will not
do so here, the proof extends with no  essential changes
 bounded smooth domains in $\R^n.$

The author would like to thank A. Hassell for collaboration on a
related project \cite{HZ} which clarified many aspects of the
boundary integral operators. He would also like to thank J. Wunsch
for discussion of singularities of boundary problems.  Finally,
the author would like to thank Y. Colin de Verdiere for many
comments, corrections and constructive criticisms on earlier
versions of this  article, both in its written form and in various
verbal presentations. His encouragement and insight during the
course of this project was much appreciated.

\tableofcontents

\section{\label{BACKGROUND} Billiards and the length functional}

We collect here some notation and background results on plane billiards which we will need below, mainly following
the reference Kozlov-Trechev \cite{KT}.

Let $\Omega$ denote a simply connected analytic plane domain with
smooth boundary $\partial \Omega$ of length $2\pi$. We denote by
${\bf T} = \R \backslash 2 \pi \Z$ the unit circle and parametrize
the boundary counter-clockwise by arc-length starting at some
point $q_0 \in
\partial \Omega$:
\begin{equation} q: {\bf T} \to \partial \Omega \subset \R^2,\;\;\; q(\phi) = (x(\phi), y(\phi)),\;\;\; |\dot{q}(\phi)| = 1, \;\;q(0) = q_0.
\end{equation}
We similarly identify the $m$-fold Cartesian product $(\partial \Omega)^m$ of $\partial \Omega$ by ${\bf T}^m,$
and denote a point of the latter by $(\phi_1, \dots, \phi_m).$

By an $m$-link  {\it periodic reflecting ray}  of $\Omega$ we mean a billiard trajectory $\gamma$ which intersects $\partial \Omega$
transversally at $m$ points $q(\phi_1), \dots, q(\phi_m)$ of intersection, and  reflects
off $\partial \Omega$ at each point according to Snell's law
\begin{equation}\label{SN} \frac{q(\phi_{j + 1}) - q(\phi_j)}{|q(\phi_{j + 1}) - q(\phi_j)|} \cdot \nu_{q(\phi_j)} =
 \frac{q(\phi_{j }) - q(\phi_{j - 1})}{|q(\phi_{j }) - q(\phi_{j - 1})|} \cdot \nu_{q(\phi_j)}. \end{equation}
Here,  $\nu_{q(\phi)}$ is the inward  unit normal to $\partial \Omega$ at $q(\phi)$. We refer to the segments
 $q(\phi_{j }) - q(\phi_{j - 1})$ as the {\it links} of the trajectory.
An $m$-link periodic reflecting ray is thus the same as
 an  $m$-link polygon in which the Snell law holds at each vertex. Since they will come up often, we make:
\begin{defn} By  $P_{(\phi_1, \dots, \phi_m)}$ we denote  the polygon with
consecutive vertices at the points $( q(\phi_1), \dots, q(\phi_m)) \in (\partial \Omega)^n$. The polygon is called:
\begin{itemize}

\item Non-singular if   $\phi_j \not= \phi_{j + 1}$ for
all $j$

\item Snell if   $P_{(\phi_1, \dots, \phi_m)}$  is non-singular and if (\ref{SN}) holds for each pair of
consecutive links   $\{q(\phi_{j  }) - q(\phi_{j - 1}), q(\phi_{j + 1}) - q(\phi_j)\}$;

\item Singular Snell if $P_{(\phi_1, \dots, \phi_m)}$ has fewer than $n$ distinct vertices, but  each non-singular
pair of consecutive links satisfies Snell's law.

\end{itemize}
\end{defn}

We will denote the acute angle between the link $q(\phi_{j + 1 }) - q(\phi_{j })$ and
the inward unit normal $\nu_{q(\phi_{j + 1})}$ by $\angle ( q(\phi_{j + 1}) - q(\phi_{j }), \nu_{q(\phi_j)})$
 and that between  $q(\phi_{j + 1}) - q(\phi_{j })$ and
the  inward unit  normal  at $q(\phi_{j })$ by $\angle (q(\phi_{j + 1 }) - q(\phi_{j }), \nu_{q(\phi_{j })}$, i.e. we put
\begin{equation}  \frac{q(\phi_{j + 1}) - q(\phi_j)}{|q(\phi_{j + 1}) - q(\phi_j)|} \cdot \nu_{q(\phi_j)} =
\cos \angle (q(\phi_{j + 1 }) - q(\phi_{j }), \nu_{q(\phi_{j })}. \end{equation}
We also use
the notation $\angle (q(\phi_{j + r} - q(\phi_{j }), \nu_{q(\phi_j)} )$ for the  angle between  the link $q(\phi_{j + r }) - q(\phi_{j })$ and the unit inward normal at $q(\phi_j).$

 The function  $\angle ( q(\phi_{j + 1}) - q(\phi_{j }), \nu_{q(\phi_j)})$ is well-defined
 on ${\bf T}^m$ minus the diagonals $\Delta_{j, j + 1} = \{\phi_j = \phi_{j + 1}\}$.  It has a continuous extension
across the diagonals according to the following
\begin{prop} $ \cos \angle ( q(\phi) - q(\phi')), \nu_{q(\phi)}) = -\frac{1}{2} \kappa(\phi) |\phi' - \phi|
+ O(|\phi' - \phi|^2).$ \end{prop}

\begin{proof} (cf. \cite{AG, EP} and \cite{Z6} ) We have:
$$(q(\phi')) - q(\phi)) \cdot \nu_{q(\phi)} = -\frac{1}{2} (\phi - \phi')^2 \kappa(\phi) + O((\phi -\phi')^3). $$
Now divide by $|q(\phi) - q(\phi')|$.

\end{proof}

\subsection{Length functional}

We first define a length functional on ${\bf T}^M$ by:

\begin{equation}\label{LENGTH}    L(\phi_1, \dots, \phi_M) =  |q(\phi_1) - q(\phi_2)| + \dots + |q(\phi_{M - 1} ) -q(\phi_M)| .\end{equation}

It is clear  that $L$ is a smooth function away
 from the  `large diagonals' $\Delta_{j, j + 1}:= \{\phi_j = \phi_{j + 1}\}$, where it   has  $|x|$   singularities .
 We have:
\begin{equation}\label{Oneder}  \begin{array}{l} \left\{ \begin{array}{l} \frac{\partial}{\partial \phi_j} |q(\phi_j) - q(\phi_{j - 1})| =
 - \sin \angle ( q(\phi_{j }) - q(\phi_{j -1 }), \nu_{q(\phi_j)}),\\ \\
\frac{\partial}{\partial \phi_{j }} |q(\phi_j) - q(\phi_{j + 1})| =
 \sin \angle ( q(\phi_{j + 1}) - q(\phi_{j }), \nu_{q(\phi_{j +
 1})})\end{array} \right.
\\ \\ \implies
 \frac{\partial}{\partial \phi_j} L  = \sin \angle ( q(\phi_{j + 1}) - q(\phi_{j }), \nu_{q(\phi_{j + 1})}) - \sin \angle ( q(\phi_{j }) - q(\phi_{j -1 }), \nu_{q(\phi_j)}) .\end{array}\end{equation}
  The condition that $\frac{\partial}{\partial \phi_j} L
= 0$ is thus that  the $2$-link defined by the triplet $(q(\phi_{
j- 1}, q(\phi_j), q_{i + 1})$ is  Snell  at $\phi_j$.  A smooth
critical point of $L$ on ${\bf T}^M$ is thus the same as an
$M$-link Snell polygon.

We will also be  concerned with the length functional $L : \Omega \times {\bf T}^M\to \R^+ $ defined by:

\begin{equation}\label{PHASE} \begin{array}{l}
  L(x, \phi_1, \dots, \phi_M) = |x - q(\phi_1)| + |q(\phi_1) - q(\phi_2)| + \dots + |q(\phi_{M - 1} - q(\phi_M)| + |q(\phi_M)- x|,
\end{array}\end{equation}
which is smooth away from the diagonals $x = q(\phi_1), x = q(\phi_M)$ together with $\Delta_{j,j + 1}.$
Its gradient in the  $x$-variable is given by
\begin{equation}\label{LXZERO}  \nabla_x L = \frac{x - q(\phi_1)}{|x - q(\phi_1)|} + \frac{x - q(\phi_M)}{|x - q(\phi_M)|}, \end{equation}
 so that
 a smooth critical point $x$ of $L(x, \phi_1, \dots, \phi_M)$
corresponds to triple $(\phi_M, x, \phi_1)$ whose $2$-link   is
{\it straight} at $x$.
We sum up in the  following well-known proposition,  due to Poincare. For background, see \cite{KT}.
\begin{prop} A   smooth critical point $(x, \phi_1, \dots, \phi_M)$ of $L$ on $\Omega \times {\bf T}^M$ corresponds to an $M$-link  Snell polygon with vertices $(x, \phi_1, \dots, \phi_M)$. \end{prop}
We will also need the formula for the interior normal derivative $\frac{\partial}{\partial \nu_y} = \nu_{q(\phi_{j + 1})} \cdot \nabla_y$
along the boundary of the link-lengths :

\begin{equation}\label{NORMDER} \begin{array}{l}  \frac{\partial}{\partial \nu_y} |q(\phi_j) - y|_{y = q(\phi_{j + 1})} =   \frac{q(\phi_j) - q(\phi_{j + 1})}{|q(\phi_j) - q(\phi_{j + 1})|}  \cdot \nu_{q(\phi_{j + 1})} =
\cos \angle ( q(\phi_{j + 1}) - q(\phi_{j }), \nu_{q(\phi_j)})
\end{array} \end{equation}

\subsection{\label{BILL} Billiard flow and length spectrum }

The (geometer's) billiard flow $\Phi^t$ of $\Omega$ is the flow
on $T^*\Omega$
 which is partially defined by Euclidean motion in the interior
and  Snell's law of reflection at the boundary. We refer to the
billiard orbits as trajectories or rays, and when they have only
transversal intersections with  the boundary we refer to them as
transversally reflecting rays. The straight line segments between
intersection points are called links.

 At (co-) vectors tangential
to $\partial \Omega$ this law does not uniquely define the flow
unless $\Omega$ is convex, in which case tangentially intersecting
 rays can only travel along the boundary.
 In the non-convex case, there exist rays
which intersect $\partial \Omega$ tangentially and at such points
the geometer's billiard flow is not uniquely defined. The
propagation of singularities theorem for domains with boundary
(cf. \cite{AM, MS, PS}) largely resolves this ambiguity, and
completely resolves it for analytic  domains. It defines the
billiard flow as the broken bicharacteristic flow of the wave
operator, i.e. the trajectories along which singularities of
solutions of the wave equation move. Roughly speaking, the
trajectories are transversally reflecting rays and limits of such
rays with many small links. The limit rays intersect the boundary
tangentially and then glide for some time along the boundary and
then re-enter the domain.  Since small links can only occur in the
convex part of the boundary, the entrance and exit points to the
boundary of a plane domain  occur at its inflection points. There
exists unique continuation of geodesics unless $\partial \Omega$
has infinite order contact with a tangent line, and of course this
cannot occur for analytic domains. For background and further
discussion we refer to \cite{PS, GM, M}.
 Pictures of gliding rays may be found
\cite{GM, M}.

In the Poisson relation for the wave equation, it is the analysts'
billiard flow (propagation of singularities) which is relevant,
and henceforth we assume the billiard flow defined as the broken
bicharacteristic flow for the wave equation.  We then define the
{\it length spectrum} $Lsp(\Omega)$ to be the set of lengths of
periodic orbits of the  billiard flow (\cite{PS}, Definition
(1.2.9)).

By the billiard map $\beta$ of $\Omega$ we mean the  map induced
by $\Phi^t$ on $B^*
\partial \Omega$: if $(q, \eta) \in B^*(\partial \Omega)$, we may
add a multiple of the unit normal to obtain an inward pointing
unit vector $v$ at $q$. We then follow the billiard trajectory of
$v$ until it hits the boundary, and then define $\beta(q, \eta)$
to be its tangential projection.

In the case of strictly convex domains, periodic orbits are either
periodic reflecting rays $\gamma$ or closed geodesics on $\partial
\Omega$ (\cite{PS}, Ch. 7). By periodic $n$-link reflecting ray,
we mean a  periodic orbit of the billiard flow $\Phi^t$ on
$T^*\Omega$ whose projection to $\Omega$ has only transversal
intersections with $\partial \Omega$. That is, $\gamma$ is a Snell
polygon with $n$ sides.  (Here, and henceforth, we often do not
distinguish notationally between an orbit of $\Phi^t$ and its
projection to $\Omega.$)

\subsubsection{Poincare map and Hessian of the length functional}

The linear Poincare map $P_{\gamma}$ of $\gamma$ is the derivative at $\gamma(0)$ of the first return map to a
transversal to $\Phi^t$ at $\gamma(0).$
By a non-degenerate periodic reflecting ray $\gamma$ we mean one whose linear Poincare map $P_{\gamma}$ has no eigenvalue
equal to one. For the definitions and background, we refer to \cite{PS}\cite{KT}.

There is an important relation between the spectrum of the Poincare map $P_{\gamma}$ of a periodic $n$-link
reflecting ray and the Hessian $H_n$ of the length functional at the corresponding critical point of $L: {\bf T}^n \to \R.$
For the following, see \cite{KT} (Theorem 3).

\begin{prop} We have: $$\det (I - P_{\gamma}) = - \det (H_n) \cdot (b_1 \cdots b_n)^{-1},$$
where $b_j = \frac{\partial^2 | q(\phi_{j + 1} - q(\phi_j)|}{\partial \phi_j \partial \phi_{j + 1}}.$\end{prop}

\section{Wave trace and resolvent trace asymptotics}

The spectral invariants we will use in determining $\Omega$ are
essentially the wave trace invariants associated to bouncing ball
orbits.  As will be recalled below, these invariants are
coefficients of the asymptotic expansion of the trace $Tr
1_{\Omega} E_{\Omega}(t)$ of the Dirichlet wave group around in
singularities at lengths of periodic billiard trajectories. Dual
to the wave trace singularity expansion, at least formally, is the
asymptotics as $ k \to \infty$ of the trace $Tr  1_{\Omega}
R_{\Omega}(k + i \tau)$ of the Dirichlet resolvent.  Further,  we
will regularize the trace and relate the resolvent trace
coefficients at a periodic reflecting ray to the corresponding
wave trace invariants.

\subsection{Resolvent and Wave group}

By the Dirichlet Laplacian $\Delta_{\Omega}$ we mean the Laplacian $\Delta =  \frac{\partial^2}{\partial x^2} +
\frac{\partial^2}{\partial y^2}$ with domain $\{u \in H^1_0(\Omega): \Delta u \in L^2\}$; thus, in our notation,
$\Delta_{\Omega}$ is a negative operator.
We denote by  $E_{\Omega}(t,x,y) = \cos t \sqrt{ - \Delta_{\Omega}} (x, y)$  the fundamental solution of the mixed wave equation with Dirichlet
boundary conditions:
\begin{equation}\label{DWK}\left\{ \begin{array}{ll} \frac{\partial^2 E_{\Omega}}{\partial t^2} =
\Delta E & \mbox{on} \; \R \times \Omega \times \Omega \\ &
\\ E_{\Omega}(0, x,y) = \delta(x - y) & \frac{\partial E_{\Omega}}{\partial t}(0,x,y) = 0 \\
& \\
E_{\Omega}(t,x,y) = 0 & (t, x, y) \in \R \times \partial \Omega \times \Omega. \end{array} \right. \end{equation}

The resolvent of the Laplacian $\Delta_{\Omega}$ on $\Omega$ with
Dirichlet boundary conditions  is the operator on $L^2(\Omega)$
defined by $$R_{\Omega}(k + i \tau) = - (\Delta_{\Omega} + (k + i
\tau)^2)^{-1}, \;\;\;\;\; \tau > 0.$$
   The resolvent
kernel, which we refer to as the   {\it Dirichlet Green's function} $G_{\Omega}(k + i \tau, x, y)$ of $\Omega \subset \R^2$,  is by definition the solution of the boundary problem:
\begin{equation}\label{GREEN}  \left\{ \begin{array}{l} (\Delta_x + (k + i \tau)^2) G_{\Omega}(k + i \tau, x, y) = - \delta(x - y),\;\;\;
(x, y \in \Omega) \\ \\
G_{\Omega}(k + i \tau, x, y) = 0, \;\;\; x \in \partial \Omega. \end{array} \right. \end{equation}

 To clarify our sign conventions, let us specify them in the  case
of the free  Laplacian $\Delta_0$ on $\R^2$.  Our $\Delta$ is negative, so
the symbol of  the free resolvent is $(|\xi|^2 - (k + i \tau)^2)^{-1}$.  The free Green's function has the asymptotic behaviour
$$G_0(k + i \tau, x, y) \sim \frac{e^{i (k + i \tau) |x - y|}}{[(k + i \tau) |x - y|]^{1/2}},\;\;\;
([(k + i \tau) |x - y|] \to \infty) $$
hence is oscillatory in $k$ and  has the  exponentially decay $e^{- \tau |x - y|]}$
when  $\tau > 0.$  We will later place the spectral parameter on
the logarithmic curve $k + i \tau \log k$, where it has a power law decay in $k$ as
well as oscillatory behaviour. For further discussion of the signs, see  \cite{T} (p. 142).

The resolvent may be expresed  in terms of the (even) wave operator as
\begin{equation} \label{LAPLACE} R_{\Omega} (k + i \tau ) =  \frac{ 1}{k + i \tau }
\int_0^{\infty} e^{i ( k + i\tau) t} E_{\Omega} (t) dt,\;\;\;\;\; (\tau > 0 ) \end{equation}
which  holds because
$$\frac{1}{\lambda^2 - (k + i \tau)^2 } = \frac{ 1}{k + i \tau }
\int_0^{\infty} e^{i ( k + i \tau) t} \cos \lambda t dt,\;\;\;\;( \forall \lambda \in \R, \; \tau \in \R^+).$$

Given  $\hat{\rho}  \in C_0^{\infty}(\R^+)$, we have  defined the smoothed resolvent $R_{\rho}(k + i \tau)$ in
(\ref{RRHO}). By (\ref{LAPLACE}) we can rewrite it in terms of the wave kernel as:
\begin{equation}\label{RTAU} \begin{array}{lll}   R_{\rho} (k + i \tau ) & =  &
\int_0^{\infty}\int_{\R} \rho(k - \mu)  e^{i ( \mu +i \tau) t} E_{\Omega} (t) dt d\mu \\& & \\
&=& \int_0^{\infty} \hat{\rho} (t)   e^{i ( k +i \tau) t} E_{\Omega} (t) dt \\ & & \\
& = & \rho(k + i \tau + \sqrt{\Delta_{\Omega}}) + \rho(k + i \tau - \sqrt{\Delta_{\Omega}}) \end{array} \end{equation}
This is essentially the  smoothing  used in the study of wave invariants in  \cite{DG}.

\subsection{Wave trace and resolvent trace asymptotics}

We now recall the classical results about the asymptotics of $Tr 1_{\Omega}  R_{\rho}(k + i \tau)$  (cf.
\cite{GM} \cite{PS}). From the last formula in (\ref{RTAU}) we note that
\begin{equation} Tr 1_{\Omega} R_{\rho}(k + i \tau) = \sum_{j =1}^{\infty}[ \rho (k - \lambda_j - i \tau) +
\rho (k + \lambda_j - i \tau)].  \end{equation}
Note that $\rho$ is an entire function since $\hat{\rho} \in C_0^{\infty}(\R)$. We also
note that
\begin{equation} Tr 1_{\Omega}  R_{\rho}(k + i \tau) = \sum_{j =1}^{\infty} \rho (k - \lambda_j - i \tau)
+ O(k^{\infty}), \;\; (k \to \infty, \tau > 0)\end{equation}
since $\lambda_j > 0$ for all $j$. This is essentially the same expression studied in \cite{DG, GM}.

Dual to the sum is the   trace of the even part of the Dirichlet wave group, i.e.  the distribution in $t$ defined by
$$Tr 1_{\Omega}  E_{\Omega} (t) :=   \int_{\Omega} E_{\Omega}(t, x,x) dx = \sum_{j = 1}^{\infty}
\cos t \lambda_j$$
where
$$\Delta \phi_j = \lambda_j^2 \phi_j,\;\;\; \phi_j |_{\partial \Omega} =0,\;\; \langle \phi_i, \phi_j \rangle =
\delta_{ij},\;\;\;\; E_{\Omega} (t,x,y) = \sum_{j} \cos t \lambda_j \phi_j(x)
\phi_j(y).$$

The singular support of the wave trace is contained in the set $
Lsp(\Omega)$ of lengths of generalized broken geodesics (\cite{AM,
GM}, Theorem; and \cite{PS}): More precisely, for  any bounded
smooth domain, we have
$$singsupp Tr 1_{\Omega}  E_{\Omega}(t) \subset Lsp(\Omega).$$
If $\Omega$ belongs to a certain residual set $\mathcal R$, then
$$singsupp Tr 1_{\Omega} E_{\Omega} (t) = Lsp(\Omega).$$

When $L_{\gamma}$ is the length of a non-degenerate periodic
reflecting ray $\gamma$, and when $L_{\gamma}$ is not the length
of any other generalized periodic orbit, then $Tr 1_{\Omega}  E_{\Omega}(t)$
is a Lagrangean distribution in the interval  $(L_{\gamma} -
\epsilon, L_{\gamma} + \epsilon)$ for sufficiently small
$\epsilon$,  hence $Tr 1_{\Omega}  E_{\rho}(k + i \tau)$ has a complete
asymptotic expansion in powers of $k^{-1}.$  Let us recall the
precise statement  (see \cite{GM}, Theorem 1, and also page 228;
see also \cite{PS} Theorem 6.3.1).
\medskip

{\it Let $\gamma$ be a non-degenerate billiard trajectory whose length $L_{\gamma}$
is isolated and of multiplicity one in $Lsp (\Omega)$.  Then for
$t$ near $L_{\gamma}$, the trace of the wave group has the singularity expansion
$$ Tr 1_{\Omega}  E_{\Omega} (t) \sim  a_{\gamma} (t - L_{\gamma} + i0)^{-1} + a_{\gamma 0} \log (t - L_{\gamma} + i 0) +
\sum_{k = 1}^{\infty}a_{\gamma k}(t - L_{\gamma} + i 0)^k \log (t - L_{\gamma} + i 0)$$
where the coefficients $a_{\gamma k}$ (the wave trace invariants) are calculated by the stationary phase method
from the Lagrangean parametrix $\hat{E}$. }
\medskip

Recall that $\hat{E}(t)$ is a microlocal parametrix in that it approximates $E_{\Omega}(t)$
modulo regular kernels in a sufficiently small conic neigbhorhood $\Gamma_L$
of $\R^+ \gamma$.

We will need the following equivalent statement:

\begin{cor} \label{BBL} Assume that $\gamma$ is a non-degenerate periodic reflecting ray, and let $\hat{\rho } \in C_0^{\infty}(L_{\gamma} - \epsilon, L_{\gamma} + \epsilon)$, equal to one on $(L_{\gamma} - \epsilon/2, L_{\gamma} + \epsilon/2)$ and
with no other lengths in its support.  Then $Tr 1_{\Omega}
R_{\rho}(k + i \tau)$  admits a complete asymptotic expansion  of
the form (\ref{PR}). The coefficients $B_{\gamma; j}$ are
canonically related to the wave invariants $a_{\gamma; j}$.
\end{cor}

\begin{proof}  By (\ref{RTAU}), we have  $R_{\rho} = E_{\rho}$ where
\begin{equation} \label{ERHO}  E_{\rho}(k + i \tau) : = \int_{\R} e^{ (i k - \tau)  t} \hat{\rho}(t)  E_{\Omega}(t)  dt. \end{equation}
The corollary thus follows immediately from the Poisson relation.

\end{proof}

As mentioned in the introduction, we actually use a variant of this result:

\begin{cor} \label{BBLLOG} Under the same assumptions,  $Tr 1_{\Omega} R_{\rho}(k + i \tau \log k)$  admits a complete
asymptotic expansion  of the form (\ref{PRLOG}) with the same coefficients $B_{\gamma, j}$.
\end{cor}

\begin{proof}  We use (\ref{ERHO}) but with $\tau \log k$ in place of $\tau$. We
then subsitute the microlocal  parametrix $\tilde{E}(t)$
and calculate
$$ Tr 1_{\Omega} \int_{\R} e^{ (i k - \tau \log k)  t} \hat{\rho}(t)  \tilde{E}_{\Omega}(t)  dt$$
asymptotically by the stationary phase method. Since supp $\hat{\rho}$ is contained
in $\R_+$, the factor $e^{ - \tau \log k  t}$ may be absorbed into the amplitude and
decreases its order. The result follows exactly as in Corollary (\ref{BBL}) by the
stationary phase method.

\end{proof}

\subsection{Microlocal cutoff}

In this section, we explain  that the regularized wave trace
expansion $Tr 1_{\Omega} R_{\rho}(k)$ at a periodic reflecting ray
$\gamma$ can be microlocalized to $\gamma$, i.e. equals $Tr 1
R_{\rho}(k) \chi(k)$, where $\chi(k)$ is a (specially adapted)
cutoff to $\gamma$. In principle this should be obvious, but we
include some details since we were unable to find a suitable
reference. For simplicity we assume that $\gamma$ is a bouncing
ball orbit.

The cutoff  consist of there terms, $$\chi(k) = \chi_+ (k) +
\chi_0 (k) + \chi_- (k) $$  corresponding to the top, middle and
bottom of the orbit. Let $U = U_+ \cup U_0 \cup U_-$ be a small
strip around $\pi(\gamma)$ in $\R^2$, with $U_{\pm}  $ be a small
neighborhood of $U \cap
\partial \Omega_{\pm}$ (the top/bottom boundary component).
Further, let $(\theta, r) \to q(\theta) + r \nu_{q(\theta)}$
denote Fermi normal coordinates along $U_+$ and we denote the dual
symplectic coordinates by $p_r, p_{\theta}$.  We use the same
notation for $U_-$, anticipating that no confusion will arise. We
then define semiclassical pseuoddifferential cutoff operators of
the form:
\begin{itemize}

\item $\chi_{\pm} (r, \theta, k^{-1} D_{\theta})$ on $U_{\pm}$, with $\chi_+(r,
\phi, p_{\theta})$ supported in  $U_+ \times \{|p_{\theta}| <
\epsilon\}$ and with $\chi_{\pm} (r, \theta, p_{\theta}) \equiv 1$
for all $0 \leq r < \epsilon, \phi \in (-\epsilon, \epsilon),
|p_{\phi}| \leq \epsilon/2. $

\item $\chi_0(x, k^{-1} D_x) $ is properly supported in $U_0 \times U_0$
and $\chi_0(x, \xi)$ is supported in a small neighborhood of $U_0
\times \{(0, 10\}.$

\end{itemize}

One can apply $\chi(k) $ to functions supported in $\overline{U
\cap \Omega}.$ In particular, the cutoff resolvent is defined near
the boundary by
$$\chi_{\pm} (k)  R_{\Omega, \rho}(k + i \tau)(r, \theta, r', \theta') = \int_{\R} \int_{\R} e^{i k (\theta - \theta'')
p_{\theta} } R_{\Omega} (k + i \tau, r, \theta''; r', \theta')
\chi_{\pm}(r, \theta, p_{\theta}) d p_{\theta} d \theta''. $$

Naturally, the  introduction of the cutoff does not change the
trace modulo negligeable terms.  For the sake of completeness, we
sketch the proof that the full resolvent trace is unchanged.

\begin{lem}\label{EASYMICRO}  $Tr 1_{\Omega} R_{\Omega, \rho}(k + i \tau \log k) - Tr 1_{\Omega} \chi(k)  R_{\Omega, \rho}(k + i
\tau \log k) = O(k^{-\infty}). $ \end{lem}

\begin{proof} It is known \cite{AM} (see also \cite{GM})   that $$Tr 1_{\Omega}  R_{\Omega, \rho}(k + i \tau)
\sim Tr 1_{\Omega}  \tilde{R}_{\Omega, \rho}(k + i \tau)$$ where
$\tilde{R}_{\Omega, \rho}(k + i \tau)$ is a microlocal parametrix
for the semiclassical Dirichlet resolvent in neighborhood of
$\gamma$. To be precise, Andersson-Melrose \cite{AM} proved that
there exists a microlocal parametrix for $E_D(t)$, the even
Dirichlet wave group, near any transversal reflecting ray. The
Fourier-Laplace transform of this parametrix is a semiclassical
resolvent parametrix. The parametrix is defined throughout $U$.

 We  verify this using (\ref{RTAU}).
First, we may write
\begin{equation} \label{RESCHI}  R_{\rho D}(k )
\chi (k) = \int_0^{\infty} \hat{\rho}(t) E_{D}(t) \chi(r, y,
|D_t|^{-1} D_y)
 e^{i (k + i \tau ) t} dt. \end{equation}
Calculating the expansion in the statement of the Lemma is the
same as computing the  wave front set of the trace of
(\ref{RESCHI}) near $r L_{\gamma}$. Thus, the statement is
equivalent to saying that

\begin{equation} \label{WF} WF[ Tr E_{D}(t) (I - \chi(r, y,
|D_t|^{-1} D_y) ] \cap (rL_{\gamma} - \epsilon, r L_{\gamma} +
\epsilon) = \emptyset.
\end{equation}
To prove this, we recall that $WF( E_D(t, x, y))$ is the
space-time graph
\begin{equation} \Gamma =  \{(t, \tau, x, \xi, x', \xi') \in T^*(\R \times \Omega^c \times
\Omega^c):\;  \tau = - |\xi|, \;\; G^t(x, \xi) = (x', \xi')\},
\end{equation}
of the generalized billiard flow. Now in Fermi coordinates $(r,
\rho dr, y, \eta dy), \R^+ \gamma$ is ray in the direction
 of $dr$. Since the space-time
 graph of  the cotangent bundle along $\gamma$ may be described in normal coordinates as a
 neighborhood of
$$\{(t, \tau, t, \tau, 0, 0)  \},$$
a conic neighborhood  may be described in these coordinates  by
$|y| \leq \epsilon, |\eta/\tau| \leq \epsilon$. This is precisely
the set to which $\tilde{\chi}_{\gamma}(r, y, |D_t|^{-1} D_y)$
microlocalizes. Emptyness of the WF in (\ref{WF}) follows from the
calculus of wave front sets, which implies that only diagonal
points in $\Gamma$ contribute, i.e. periodic orbits of the
billiard flow, and from our assumption that $\gamma^r$ is the only
orbit with period in the given set.

\end{proof}

 The only potentially confusing
issue is in the choice of cutoff operator, which must be rather
special since it operators on a manifold with boundary. Let us
verify in another way that this kind of cutoff operator acts as a
microlocal cutoff to $\gamma$.

According to \cite{GM, AM}, we can calculate the trace using a
microlocal  parametrix
$$\tilde{E}(t, x, y) = \int_{\R^2} e^{i \phi(t, x, y, \xi)} a(t,
x, y, \xi) d \xi$$ for the Dirichlet wave kernel. We take its
Fourier-Laplace transform to get a semiclassical  microlocal
parametrix
$$\tilde{R}_{\Omega, \rho} (k + i \tau, x, y) = \int_{\R^2}\int_0^{\infty} \hat{\rho}(t)
 e^{i (k + i \tau) t}  e^{i \phi(t, x, y, \xi)} a(t,
x, y, \xi) d \xi.$$ We now change variables $\xi \to k \xi$ to
obtain
$$\tilde{R}_{\Omega, \rho} (k + i \tau, x, y) =k^2  \int_{\R^2}\int_0^{\infty} \hat{\rho}(t)
 e^{i k \phi(t, x, y, \xi) } e^{- \tau t} a(t,
x, y, k \xi) d \xi,$$ where the phase is $\Phi = t  +  \phi(t, x,
y, \xi)$.  We now apply $\chi(k)$ to get
$$\chi(k) \; \tilde{R}_{\Omega, \rho} (k + i \tau, x, y) = k^2  \int_{\R^2}\int_0^{\infty} \hat{\rho}(t)
 e^{i k [t  +  \phi(t, x, y, \xi)]} \tilde{\chi}(x, d_x \phi) e^{- \tau t} a(t,
x, y, k \xi) d \xi,$$ where $$\tilde{\chi}(x, d_x \Phi) = \chi(x,
d_x \phi) A(k, x, d_x \phi)$$ with $A$ a symbol of order $0$. The
assumption on the cutoff implies that $\tilde{\chi} \equiv 1$ near
$\gamma$. Since the trace is computed by applying stationary phase
to the trace of this oscillatory integral, it is unchanged modulo
rapidly decaying errors by the cutoff.

\section{\label{LP} Multiple reflection expansion of the resolvent}

The purpose of this section is to review the `multiple-reflection'
expansion of the Dirichlet Green's function of a bounded plane
domain. This is the term in \cite{BB1, BB2} for the Neumann series
expression for the Dirichlet Green's function in terms of double
layer potentials.   The same method also works for Neumann
boundary conditions, but for simplicity we only explicitly treat
the Dirichlet case.  We refer to  \cite{T}, Chapter 5 for
background in  potential theory.

 The method
of layer potentials (\cite{T} II, \S 7. 11) seeks to solve (\ref{GREEN})  in terms
of the `layer potentials' $G_0(k + i \tau, x, q), \partial_{\nu_y} G_0(k + i \tau, x, q) \in {\mathcal D}'(\Omega \times \partial \Omega)$,  where $\nu$ is the interior unit normal to $\Omega$,  where $\partial_{\nu} = \nu \cdot \nabla,$ and where $G_0(k + i \tau, x, y)$ is the `free'
Green's function of $\R^2$,
i.e. of the kernel of the free resolvent $- (\Delta_0 + (k + i \tau)^2)^{-1}$ of the Laplacian $\Delta_0$ on $\R^2$.
 The free Green's function in dimension two  is given by:
$$\begin{array}{l} G_0(k + i \tau, x, y) =
  H^{(1)}_0((k + i \tau) |x - y|)  = \int_{\R^2} e^{i \langle x - y, \xi \rangle} (|\xi|^2 - (k + i \tau)^2)^{-1} d \xi. \end{array}$$
Here, $H^{(1)}_0(z)$ is the Hankel function of index $0$.  In general, the Hankel function of index $\nu$ has the integral representations (\cite{T}, Chapter 3, \S 6)
\begin{equation} \label{HANKEL} \begin{array}{lll} H^{(1)}_{\nu} (z) & = & \frac{2 e^{-i \pi \nu} }{ i \sqrt{\pi} \Gamma(\nu + 1/2)} (\frac{z}{2})^{\nu}
\int_1^{\infty} e^{i z t}  (t^2 - 1)^{\nu -1/2} dt \\& &  \\
 & = &  (\frac{2}{\pi z})^{1/2} \frac{e^{i(z - \pi \nu/2 -  \pi/4)}}{ \Gamma(\nu + 1/2)}
\int_0^{\infty} e^{-s} s^{-1/2} (1 - \frac{s}{2 i z})^{\nu -1/2} ds.
\end{array} \end{equation}
From the first, resp. second, representation we derive the asymptotics:
\begin{equation}\label{LK} H^{(1)}_{\nu} ((k + i \tau) |x - y|) \sim \left\{ \begin{array}{l}
 - \frac{1}{2\pi} \ln (|k + i \tau| |x - y|) \;\rm{as}\; |k + i \tau| |x - y| \to 0, \; \mbox{if}\; \nu = 0\\ \\
 - \frac{i \Gamma(\nu)}{\pi} (\frac{2}{|k + i \tau| |x - y|} )^{\nu}  \;\rm{as}\; |k + i \tau| |x - y| \to 0,  \; \mbox{if}\; \nu > 0\\ \\
 e^{i ((k + i \tau) |x - y| - \nu \pi/2 - \pi/4} \frac{1}{(|k + i \tau| |x - y|)^{1/2}}\;\rm{as}\; |k + i \tau| | |x - y| \to \infty.  \end{array} \right. \end{equation}
Here it is assumed that $\tau > 0$.

The single layer, respectively double layer, potentials are  the operators
\begin{equation} \label{layers} \begin{array}{ll} {\mathcal S} \ell(k + i \tau)  f(x) = \int_{\partial \Omega} G_0(k + i \tau, x, q) f(q)
d s (q), &  \;\;  {\mathcal D} \ell(k + i \tau)f(x) =  \int_{\partial \Omega}\frac{\partial}{\partial \nu_y} G_0(k + i \tau, x, q) f(q) ds(q),
 \end{array} \end{equation}
where $ds(q)$ is the arc-length measure on $\partial \Omega$.
They induce boundary operators
\begin{equation} \label{blayers} \begin{array}{ll} (i) & S(k + i \tau)  f(q) = \int_{\partial \Omega} G_0(k + i \tau, q, q') f(q')
d s (q'),\\ & \\(ii)  &    N(k + i \tau)f(q) =  2 \int_{\partial \Omega}\frac{\partial}{\partial \nu_y} G_0(k + i \tau, q, q') f(q') ds(q') \end{array} \end{equation}
which map $H^s(\partial \Omega) \to H^{s + 1}(\partial \Omega).$ Furthemore one has (cf. \cite{T}II, Proposition 11.5)
\begin{equation}  \label{ISO} \begin{array}{ll} (i) &
S(k + i \tau), \; N(k + i \tau) \in \Psi^{-1}(\partial \Omega), \\ &  \\ (ii) &
 (I + N(k + i \tau)): H^s(\partial \Omega) \to H^s(\partial \Omega)\;\; \mbox{is\; an\;
isomorphism}. \end{array} \end{equation}

To understand the role of these operators, it helps to recall  that the Poisson integral operator
\begin{equation} \label{INTPI} PI_{\Omega}(k + i \tau): H^s(\partial \Omega) \to
H^{s + 1/2} (\Omega),\;\; PI_{\Omega}(k + i \tau) u(x) = \int_{\partial \Omega} \partial_{\nu} G_{\Omega}(k + i \tau, x, q) u(q) ds (q) \end{equation}
 may be expressed in the form
\begin{equation} \label{PIF} PI_{\Omega}(k + i \tau)   = 2 {\mathcal D}\ell(k + i \tau)(I + N(k + i \tau))^{-1}. \end{equation}
We refer to  \cite{T}I (see  Chapter 5, Proposition 1.7; see also
\cite{LM}  ) for background on layer potentials.

We include some other  facts about the kernel $N(k + i \tau)$
which will be used to estimate remainders in  the multiple
reflection expansion (see Proposition (\ref{FIRST}) and Lemma
(\ref{HSUSED}). Further and more precise estimates will be given
in \cite{Z5}.

\begin{prop} \label{N} Suppose that $\partial \Omega$ is $C^1$. Then:

\noindent{\bf (i)}  $N(k + i \tau, q(\phi_1), q(\phi_2)) \in C^{1 - \epsilon}({\bf T} \times {\bf T}).$

\noindent{\bf (ii)} $N(k + i \tau)$ is a Hilbert-Schmidt operator
on $L^2(\partial \Omega),$ with $ ||N(k + i \tau)||_{HS} \leq C
|k|^{1/2}. $

\noindent{\bf (iii)}  $N(k + i \tau) \in \Psi^{-2}(\partial
\Omega),$ hence it is a trace-class operator on $L^2(\partial
\Omega).$

 \end{prop}

\begin{proof}

\noindent{\bf (i)}
By definition,
$$\begin{array}{lll} N(k + i \tau, q(\phi_1), q(\phi_2)) &=& \partial_{\nu_y} G_0(k + i \tau, q(\phi_1), y)|_{y = q(\phi_2)} \\ &&\\
&=& - (k + i \tau) H^{(1)}_1 (|k + i \tau| |q(\phi_1) - q(\phi_2)|) \cos \angle(q(\phi_2) -q (\phi_1), \nu_{q(\phi_2)}). \end{array}$$
  Now at $r \sim 0,$ we have (see (\ref{LK}),  also \cite{AS}, 9.1)
and \cite{EP}),
\begin{equation} \label{HANKELSING} H^{(1)}_1(r) = \frac{2i}{\pi r} + {\mathcal O}(r ( 1 + \log r)). \end{equation}
We correspondingly define
$$\begin{array}{l} N(k + i \tau) = N_{sing} + N_{reg}(k + i \tau),\;\;\;\mbox{with}
\\ \\ N_{sing}(\phi, \phi') = \frac{1}{2\pi |q(\phi) - q(\phi')|} \nu_{q(\phi')} \cdot
 \frac{q(\phi) - q(\phi')}{|q(\phi) - q(\phi')|}
= -\frac{1}{4\pi} \kappa(\phi) + O((\phi - \phi')).\end{array}$$
In fact, the vanishing of the numerator to order $2$ implies that
$N_{sing}$ is a $C^{\infty}$ kernel. Hence the smoothness of $N(k
+ i \tau, q(\phi_1), q(\phi_2))$ equals that of $ N_{reg}(k + i
\tau, q(\phi_1), q(\phi_2))$. It follows that the kernel has the
regularity of $x \log x$, so  $N(k + i \tau, q(\phi_1),
q(\phi_2))$ is Lipschitz continuous of any exponent $\alpha < 1.$
\qed

\noindent{\bf (ii)}  Obviously the kernel is Hilbert-Schmidt. We
need to estimate \begin{equation} \label{SQUARE} \int_{\partial
\Omega \times
\partial \Omega} |N(k + i \tau, q_1, q_2)|^2 ds(q_1) ds(q_2).
\end{equation}

 For the sake of brevity, we only  sketch the norm estimate,
referring to \cite{Z5} for further  details.  We break up the
domain of integration into the three regions
$$|q - q'| \leq |k + i \tau|^{-1},  \;\;\; |k + i \tau|^{-1} \leq
|q - q'| \leq |k + i \tau|^{-2/3},  \;\;\; |q - q'| \geq |k + i
\tau|^{-2/3}.
$$

In the first (near diagonal) region, the kernel is bounded above
by a constant independent of $k + i \tau$. Indeed, both the
singularity and the factor of $(k + i \tau)$ cancel in the most
singular term, as can be seen from (\ref{HANKELSING}) or
(\ref{LK}). The smoother terms also  cancel the factor of $(k + i
\tau)$ in region (i).

In region (ii), the Hankel factor has a uniform upper bound. The
cosine factor puts in $|q - q'|$. Hence, in any  region $|k + i
\tau|^{-1} \leq |q - q'| \leq |k + i \tau|^{-r}  $ the integral is
dominated by
$$|k + i \tau|^2 \int_0^{|k + i \tau|^{-r}} x^2 dx = O(|k + i
\tau|^{2 - 3 r}). $$ So it is uniformly bounded if $r \geq 2/3$.

In region (iii), we can (and will) use the WKB expansion of the
Hankel function. The square of the Hankel function contributes the
factor $|(k + i \tau)|^{-1}$, cancelling one power in
(\ref{SQUARE}). Taking the square root of the result gives the
stated estimate. \qed

\noindent{\bf (iii)} It is standard that $N(k + i \tau) \in
\Psi^{-1}(\partial \Omega)$ \cite{T}.  The  smoothness of
$N_{sing}$ noted in (i) shows that the  symbol of $N$ of order
$-1$ vanishes; hence,  it is actually of order $-2$. Therefore, it
is of trace class.

An alternative proof is by the  Hille-Tamarkin theorem (cf.
\cite{GD}) theorem, which states that a Hilbert-Schmidt kernel
$K(x,y)$
 on an interval (or circle) is trace class if it is $  Lip_{\alpha}$ in one of its variables with $\alpha > 1/2$.
 \qed

This completes the proof of the Proposition.

\end{proof}

In \cite{Z5}, we prove that the operator norm of   $N(k + i \tau)$
 on $L^2(\Omega)$ is uniformly bounded in $k$ for $\tau
> 0$. On the other hand, it appears that  $N(k + i \tau)$ is not a
contraction for any fixed $\tau > 0$.

\subsection{Multiple-reflection expansion: Dirichlet boundary conditions}

To solve (\ref{GREEN}), one puts $\Gamma(k + i \tau, x, y) = G_{\Omega}(k + i \tau, x, y) - G_0(k + i \tau, x, y)$ and solves
the boundary value problem
\begin{equation} \label{DP} \left\{ \begin{array}{l} (\Delta_x - (k + i \tau)^2) \Gamma (k + i \tau, x, y) = 0,\;\;
x, y \in \Omega  \\ \\
r_{\Omega} \Gamma(k + i \tau, q, y) = - r_{\Omega} G_0 (k + i
\tau, q, y),\;\;\; q \in \partial \Omega.  \end{array} \right.
\end{equation} Here,
$$  r_{\Omega} u = u|_{\partial \Omega}$$
is the restriction operator acting in the first variable.

We try  to solve (\ref{DP}) with $\Gamma(k + i \tau, \cdot, y)$ of
the form: $\Gamma(k + i \tau, x, y) = {\mathcal D} \ell(k + i
\tau)  \mu.$ If we denote by
$$f_{+}(q) = \lim_{x \to q, x \in \Omega} f(x) \in C(\bar{\Omega}),$$
then one has (cf. \cite{T}II, Proposition 11.1)
$$[{\mathcal D}\ell(k + i \tau) f]_+(x) = 1/2 f(x) + 1/2 N(k + i \tau) f(x).$$
Hence the equation for $\mu(k + i \tau, q, y)$ is given by:
\begin{equation} \label{MU} (I +  N(k + i \tau)) \mu(k + i \tau, q, y) = - r_{\Omega} G_0(k + i \tau, q, y),\;\;\;\; (q, y) \in \partial \Omega \times
\Omega. \end{equation}
By (\ref{ISO}) there exists a unique solution,
\begin{equation} \label{MUmore}  \mu(k + i \tau, q, y) = - (I +  N(k + i \tau))^{-1} r_{\Omega} G_0(k + i \tau, q, y),\;\;\;\; (q, y) \in \partial \Omega \times \Omega, \end{equation}
where by (\ref{blayers}), $(I + N(k + i \tau))^{-1} \in \Psi^0(\partial \Omega).$  It follows that
\begin{equation} \label{NS} \begin{array}{lll} R_{\Omega}(k + i \tau) &=&
 R_0(k + i \tau) - {\mathcal D} \ell(k + i \tau) (I +  N(k + i \tau))^{-1} r_{\Omega}  R_0(k + i \tau)
\\ & & \\
& =& R_0(k + i \tau) + \sum_{M = 0}^{\infty}(-2)^M  {\mathcal D}
\ell(k + i \tau)  N(k + i \tau)^{M} r_{\Omega} R_0(k + i \tau).
\end{array}
\end{equation} More precisely, the kernels are related by the
Neumann series
\begin{equation} \label{MRE} \begin{array}{l} G_{\Omega}(k + i \tau, x, y) =
\sum_{M = 0}^{\infty} (-2)^M G_M(k + i \tau, x, y), \;\;\;
\mbox{where, for} \;\; M \geq 1,\\ \\
G_M(k + i \tau, x, y) = \int_{( \partial \Omega)^M}
\frac{\partial}{\partial \nu} G_0(k + i \tau, x, q_1) \Pi_{j =
1}^{M-1}  \frac{\partial}{\partial \nu} G_0(k + i \tau, q_j, q_{j
+ 1}) G_0(k + i \tau, q_M, y) \Pi_{j = 1}^M ds(q_j). \end{array}
\end{equation} Here, $\frac{\partial}{\partial \nu}$ is short for $\frac{\partial}{\partial \nu_y}$ This expansion of  $G_{\Omega}(k + i \tau, x, y)$
of $R_{\Omega}(k + i \tau)$ is referred to in \cite{BB1} as the
{\it multiple-reflection expansion}.
We regularize $G_M$ to $G_{M, \rho}$ as in (\ref{RGM}).
By substituting in the multiple reflection expansion (\ref{MRE}), we obtain a multiple reflection
expansion for $R_{\rho}$.

In the following proposition, we assume the boundary is smooth since we are dealing  with analytic
boundaries. Some additional work is needed if the boundary is only assumed piecewise smooth, as in \cite{Z4}.

\begin{prop}\label{FIRST} Suppose that $\partial \Omega$ is $C^{\infty}$. Then, for $\tau \geq 0$,  $G_{M, \rho}(k + i \tau, x, y)$ defines for each $M$
a trace class operator on $L^2(\Omega)$.

\end{prop}

\begin{proof}

First, the $M = 0$  term $ \int_{\R }\rho(k - \mu) 1_{\Omega} (\mu
+ i \tau) R_0(\mu + i \tau) 1_{\Omega} d\mu$ is easily seen
 to be trace class since it is the restriction of a smoothing operator to $\Omega.$

Next, we observe that   for $M \geq 1$,   $G_{M}(k + i \tau, x,
y)$ is the Schwartz kernel of ${\mathcal D} \ell (k + i \tau) N(k
+ i \tau)^{M - 1} r_{\Omega} R_0(k + i \tau) 1_{\Omega}$. Let us
consider each factor. Using (\ref{PIF}) it follows that
$${\mathcal D} \ell (k + i \tau): H^s(\partial \Omega) \to H^{s +
1/2} (\Omega)\;\; \mbox{continuously}. $$ Further (see \cite{T},
Chapter 4, Proposition 4.5), the restriction operator satisfies
 $r_{\Omega} : H^s(
\bar{\Omega}) \to H^{s - 1/2} (\partial \Omega),$ for $ s > 0. $
(See also \cite{LM}).  And $1_{\Omega} R_0(k + i \tau)1_{\Omega}:
H^s(\Omega) \to H^{s+2}(\Omega). $

It follows that  for $M \geq 2,$ $${\mathcal D} \ell (k + i \tau)
N(k + i \tau)^{M - 1} r_{\Omega} R_0(k + i \tau) 1_{\Omega} :
H^s(\Omega) \to H^{s + 2 (M - 1)  + 2}(\Omega).$$ Hence, this
operator is certainly   of trace class on $L^2(\partial \Omega)$
if $M \geq 2.$

The case $M = 1$ is different from the others, so it seems worth
considering it separately. The operator
 $ {\mathcal D} \ell
(\mu)  r_{\Omega} R_0(\mu) 1_{\Omega} $ has   Schwartz kernel
$$\int_{\R} \int_{\partial \Omega}\rho (k -\mu) (\mu + i \tau)   \partial_{\nu_y} G_0(\mu, x, q(\phi)) G_0(\mu, q(\phi), y) d \phi.$$
The double layer potential has kernel
\begin{equation} \label{NORDERGR} \begin{array}{lll} \partial_{\nu_y} G_0(\mu, x, q(\phi)) &=&
-( \mu + i \tau) H^{(1)}_1 ((\mu + i \tau)|q(\phi_1) - x|) \cos
\angle(q(\phi) - x,  \nu_{q(\phi)}).\end{array} \end{equation} and
as for $N(k + i \tau)$ above, we write
$$ \partial_{\nu_y} G_0(k + i \tau, x, q(\phi)) = K_{reg}(k + i \tau, x, q(\phi)) + K_{sing}(k + i \tau, x, q(\phi))$$
where
$$\begin{array}{lll} K_{sing}(k + i \tau, x, q(\phi)) =  \frac{1}{2\pi} \frac{1}{|x - q(\phi)|}
\frac{x - q(\phi)}{|x - q(\phi)|} \cdot \nu_{q(\phi)}.
\end{array}$$ It is clear that
$$\int_{\partial \Omega} K_{reg}(\mu + i \tau, x, q(\phi)) G_0(\mu + i \tau, q(\phi), y) d \phi$$
is a composition of a Hilbert-Schmidt kernel from $L^2(\Omega) \to
L^2(\partial \Omega)$ and one from  $L^2(\partial \Omega) \to
L^2(\Omega)$. Hence, this term is of trace class. The $K_{sing}$
term has kernel
$$ \int_{\partial \Omega}    \frac{1}{|x - q(\phi)|}
\frac{x - q(\phi)}{|x - q(\phi)|} \cdot \nu_{q(\phi)} \{ \int_{\R}
\rho (k -\mu) (\mu + i \tau)  G_0(\mu + i \tau, q(\phi), y)d \mu
\} d \phi $$ is of trace class.   However, the bracketed operator
is  a smoothing operator from $L^2(\Omega) \to L^2(
\partial \Omega).$ Hence,  the composition
 is of trace class.

\end{proof}

\subsection{The operators $N_0$ and $N_1$}

We now  describe $N_1(k + i \tau)$ of (\ref{N01}) - (\ref{N01DEF})
as a semiclassical Fourier integral operator. For related results,
see \cite{HZ}.
 We denote by $\chi$ a smooth bump function which is supported on $[-1,1]$
and equals $1$ on $[-1/2, 1/2]$. We also now put in $h = k^{-1 + \delta}$.

We denote by $S^p_{\delta}({\bf T}^m)$ the class of symbols $a(k, \phi_1, \dots, \phi_m)$ satisfying:
\begin{equation}\label{SYMBOL}  |(k^{-1} D_{\phi})^{\alpha} a(k, \phi)| \leq C_{\alpha} |k|^{p - \delta |\alpha|},\;\;\;
(|k| \geq 1).  \end{equation}

\begin{prop}\label{NSYM} $N_1(k + i \tau)$ is a semiclassical
Fourier integral operator of order $0$ associated to the billiard
map. More precisely,
 $$\begin{array}{l}  N_1(k + i \tau, q(\phi_1), q(\phi_2)) = (1 - \chi(k^{1 -
\delta} (\phi_1 - \phi_2) )) \\ \\ \cdot (k + i
 \tau)^{\frac{1}{2}} a_1(k + i \tau, q(\phi_1),
q(\phi_2)) e^{i (k + i \tau) |q(\phi_1) - q
(\phi_2)|}\end{array},$$ with $ a_1(k + i \tau, q(\phi_1),
q(\phi_2)))   \in S^{0}_{\delta}({\bf T}^2). $
\end{prop}

\begin{proof}

We begin by analyzing  the amplitudes of the Hankel functions.

\begin{lem} \label{HANKSYM} There exist amplitudes $a_0, a_1$ such that:\\

\noindent(i)  $H^{(1)}_0( (k + i \tau) z) = e^{i(k + i \tau) z}
a_0((k + i \tau)z),$
where  $(1 - \chi(k^{1 - \delta} z))  a_0((k + i \tau)z)   \in S^{-1/2}_{\delta} (\R) $;\\

\noindent(ii) $ (k + i \tau) H^{(1)}_1( (k + i \tau) z)  = (k + i
\tau)^{\frac{1}{2}} e^{i((k + i \tau) z} a_1((k + i \tau)z),$ with
$ (1 - \chi(k^{1 - \delta} z)) a_1((k + i \tau)z) \in
S^{0}_{\delta} (\R)$.   \end{lem}

\begin{proof} By the explicit formula (\ref{LK}),

\begin{equation} a_0((k + i \tau)z) = (\frac{2}{\pi (k + i \tau) z})^{1/2}  \int_0^{\infty} e^{-s} s^{-1/2}
(1 - \frac{s}{2 i (k + i \tau) z})^{-1/2} ds.
 \end{equation}
We claim that  the integral defines a polyhomogeneous symbol  of
order $0$ in $(k + i \tau)z$ as $|(k + i \tau)z| \to \infty$.
Indeed, applying  the binomial theorem to the factor $(1 -
\frac{s}{2 i (k + i \tau) z})^{-1/2}$ gives
\begin{equation}\begin{array}{lll}  a_0((k + i \tau)z) & = &  (\frac{2}{\pi (k + i \tau) z})^{1/2}
 \int_0^{\infty} e^{-s} s^{-1/2} \sum_{j = 0}^N C(-1/2, j)
 ( \frac{s}{2 i (k + i \tau) z})^j ds  +  R_N( (k + i \tau)z)\\ & &  \\
&=&  \sum_{j = 0}^N c_j ((k + i \tau) z)^{-j - 1/2} + R_N((k + i
\tau)z), \end{array}
 \end{equation}
where $C(-1/2, j)$ are binomial coefficients, where $c_j$ are the
resulting  constants, and where $R_N(s, (k + i \tau)z) =
\int_0^{\infty} e^{-s} s^{-1/2} R_N((k + i \tau) z, s) ds$ with
$R_N((k + i \tau) z, s)$ the $N$th order remainder in the Taylor
series expansion of $(1 - \frac{s}{2 i (k + i \tau) z})^{-1/2}$.
It is evident that $R_N(s, (k + i \tau)z) = O(|(k + i \tau)z|^{-(N
+ 1/2) })$.  Differentiating with $k^{-1} D_z$ similarly gives
\begin{equation}\begin{array}{lll}  (k^{-1} D)^{\alpha}_z a_0((k + i \tau)z) &=& \sum_{j = 0}^N c_j ((k + i \tau) z)^{-j - 1/2 - |\alpha|} +
(k^{-1} D)^{\alpha}_z R_N(kz)\\ & &\\
&=& O(|kz|^{-1/2 + |\alpha}). \end{array}
 \end{equation}

Now  $1 - \chi(k^{1 - \delta} z)$ clearly belongs to
$S^0_{\delta}(\R)$. Since $|k z| \geq k^{ \delta}$ on supp $1 -
\chi(k^{1 - \delta} z)$, we conclude that
 $$(k^{-1} D)^{\alpha}_z (1 - \chi(k^{1 - \delta} z)) a((k + i \tau)z) = O( k^{- (1/2 + |\alpha| \delta)}),$$
proving (i).

By definition, $H^{(1)}_1(  z) = - \frac{d}{dz}
 H^{(1)}_0( z)$, so (ii) follows immediately from (i).
\end{proof}

We now complete the proof of the Proposition. The amplitude of
$N_1$ is then \begin{equation} a_1(k + i \tau, q(\phi_1),
q(\phi_2)) :=  a_1 ((k + i \tau)|q(\phi_1) - q(\phi_2)|) \cos
\vartheta_{1,2}, \end{equation} where $\vartheta_{1,2} = \angle
q(\phi_2) - q(\phi_1), \nu_{q(\phi_2)} ). $

Since $q: {\bf T} \to
\partial \Omega$ is smooth, the metrics $|\phi_1 -  \phi_2|_{{\bf
T}}$ and $|q(\phi_1) - q(\phi_2)|_{\R^2}$ are equivalent,  and we
use the former. Thus, we need to check that
\begin{equation} \label{ORD}  |(k^{-1} D_{\phi})^{n}
  a_1 ((k + i \tau)|q(\phi_1) - q(\phi_2)|) \cos
\vartheta_{1,2}| \leq C k^{- n
  \delta}. \end{equation}

By repeatedly differentiating (\ref{Oneder}) away from the
diagonal, we  find
\begin{equation}\label{PHISYM}  | D_{\phi}^{\alpha} (q(\phi_1) -  q(\phi_2)) | \leq C_{\alpha} |\phi_1 - \phi_2|^{1 - |\alpha|}. \end{equation}
It follows from (\ref{Oneder}) and by the chain rule that
\begin{equation} \label{CHAIN} \begin{array}{lll}  (k^{-1} D_{\phi})^{\alpha }
  a_1 ((k + i \tau)|q(\phi_1) - q(\phi_2)|) & = &
\sum_{j: j \leq |\alpha|}
 (k^{-1} D_z))^{j }  a_1( (k + i \tau) z)) |_{z = |q(\phi_1) -  q(\phi_2)|} \\ & & \\
 & \times & \sum_{ \gamma_1, \dots, \gamma_j,
|\gamma| = |\alpha| - j} C_{\alpha, j}  \Pi_{\ell = 1}^j  (k^{-1}
D_{\phi})^{\gamma_{\ell}} \sin \vartheta_{1, 2}. \end{array}
\end{equation}

By Lemma \ref{HANKSYM}(ii),
\begin{equation}\label{HAN} | (1 - \chi(k^{1 - \delta} (\phi_1 - \phi_2) ) )(k^{-1} D_z))^{j }  a_1( (k + i \tau) z)) |_{z = |q(\phi_1) -  q(\phi_2)|}|
\leq C_{j} k^{ - j \delta}. \end{equation} Since $|\phi_1 -
\phi_2| \geq k^{-1 + \delta}$ on supp $(1 - \chi(k^{1 - \delta}
(\phi_1 - \phi_2) ))$ it follows from   (\ref{PHISYM}) that
\begin{equation} \label{COSDECAY} | (1 -  \chi(k^{1 - \delta} (\phi_1 - \phi_2) ))(k^{-1} D_{\phi})^{\gamma_{\ell}}
\sin \vartheta_{12}| \leq C_{|\gamma|} k^{  - |\gamma_{\ell}|
\delta}.
\end{equation} The same kind of estimate is correct for the factor $\cos \vartheta_{12}$. The proof folows from the combination of
(\ref{CHAIN})-(\ref{HAN})-(\ref{COSDECAY}) with the fact that  $(1
- \chi(k^{1 - \delta} (\phi_1 - \phi_2) )) \in S^0_{\delta}({\bf
T}^2)$.

\end{proof}

\subsection{Layer potentials as semi-classical Fourier integral
operators}

In this section, we give a description of the layer potentials
$\SL, \DL$ as Fourier integral operators from $L^2(\partial
\Omega) \to L^2(\Omega)$ which parallels that of Proposition
\ref{NSYM}. These operators also appear in the trace formula, and
we will be needing the results of this section in \S 7.  First, we
must introduce suitable cutoff operators away from the diagonal.

To define the cutoff operator, we need to introduce coordinates.
We separate $\Omega$ into two zones depending on the distance
$r(x,
\partial \Omega)$ to $\partial \Omega.$ We set
\begin{equation} \label{ZONES} \left\{\begin{array}{l} (\partial \Omega)_{\epsilon}
 = \{x \in \Omega: r(x, \partial \Omega) < \epsilon \} \\ \\
\Omega_{\epsilon} = \Omega \backslash (\partial
\Omega)_{\epsilon}.
\end{array}\right. \end{equation} Here,   $\epsilon$ is
 sufficiently small so that
the exponential map
$$ \exp : N^{\bot} (\partial \Omega) \to \Omega,\;\;\;  \exp_{q(\phi_0) }r \nu_{q(\phi_0)} = q(\phi_0) + r \nu_{q(\phi_0)}  $$
from the inward normal bundle along the boundary to the interior
is a diffeomorphism from vectors of length $< \epsilon$ to
$(\partial \Omega)_{\epsilon}$. Also, as above,  $\nu_{q(\phi_0)}$
is the interior unit normal at $q(\phi_0)$. The exponential map
induces Fermi normal coordinates $(\phi, r)$ on the     annulus
$(\partial \Omega)_{\epsilon}$, with   $r = r(x, \partial \Omega)$
the  distance from $x$ to $\partial \Omega.$  We denote the
Jacobian of the exponential map along the boundary by $J(\phi,
r)$. We  introduce the corresponding  cutoff
\begin{equation}
\chi_{\partial \Omega}^{k}(x) := \chi( k^{1 - \delta} r),
\end{equation}
where it is understood that the cutoff is supported in a tube
around the boundary where the Fermi normal coordinates are
defined. We further introduce an angular cutoff as for the
boundary integral operators.  Near the boundary we may use Fermi
normal
 coordinates $(r, \theta)$ and we introduce  $\chi(k^{1 -
 \delta}(\theta - \phi)$ where $q = q(\phi)$ in arclength
 coordinates.

We  then break up each layer potential into several pieces using
the cutoffs. For instance, we first  make a radial cutoff of the
single layer:
\begin{equation} \label{RADCUT}  \chi_{\partial \Omega}^{k^{-1 + \delta}}(x)) \SL +   \;\;
 (1 - \chi_{\partial \Omega}^{k^{-1 + \delta}}(x)) ) \SL. \end{equation}
We then further break up $ \chi_{\partial \Omega}^{k^{-1 +
\delta}}(x))
 \SL$ into
 \begin{equation}\label{TANCUT} \begin{array}{ll} \chi_{\partial \Omega}^{k^{-1 + \delta}}(x))
\SL(k + i \tau, (r, \theta), q(\phi)) & =
  \chi_{\partial \Omega}^{k^{-1 + \delta}}(x)) \chi(k^{1 -
 \delta}(\theta - \phi) \SL(k + i \tau, (r, \theta), q(\phi)) \\ & \\ & +
  \chi_{\partial \Omega}^{k^{-1 + \delta}}(x))
  (1 -  \chi(k^{1 -
 \delta}(\theta - \phi)) \SL(k + i \tau, (r, \theta), q(\phi))  \end{array} \end{equation}
We do likewise with the double layer potential.

We now show that, when suitably cutoff away from the diagonal
singularities, the  operators $\SL, \DL$ are semiclassical Fourier
integral operators  with phase function $|x - q|: \Omega \times
\partial \Omega \to \R.$ The phase parametrizes the canonical
relation \begin{equation} \label{GAMMA} \Gamma = \{(x, \frac{x -
q}{|x - q|}, q, - \frac{x - q}{|x - q|} \cdot T_q) \} \subset T^*
\Omega \times T^*
\partial \Omega \end{equation} which is the graph of the interior-to-boundary
billiard map $\beta$, which takes an interior vector $(x, \frac{x
- q}{|x - q|})$ to the tangential component of the tangent
vector(s) to the billiard ray it generates at its intersection
point(s) with the boundary. The natural projections $p: \Gamma \to
T^*\Omega$, resp. $ q: \Gamma \to T^*\partial \Omega$ have
singularities at points where the billiard ray intersects the
boundary tangentially. Away from these singular points, $p$ resp.
$q$ are submersions with discrete, resp. $1$ dimensional fibers.
Phases of this type arose in the work of Carleson-Sjolin (though
not in relation to boundary value problems) and are discussed in
detail in Sogge \cite{So}. Since we are cutting off the diagonal
we have eliminated grazing rays. In   the case of convex domains,
the cutoff  layer potentials are non-degenerate Fourier integral
operators,  while in the non-convex they have degeneracies at
points $(x, \frac{x - q}{|x - q|}, q, 0) \in \Gamma$. These too
will be eliminated when we microlocalize to the orbit $\gamma$.

 We have:

\begin{prop}\label{LSYM} The operators
\begin{itemize}

\item $ (1 - \chi_{\partial \Omega}^{k^{-1 + \delta}}(x)) ) \DL(k + i \tau)$ (resp. $\SL(k + i \tau)$);

\item $\chi_{\partial \Omega}^{k^{-1 + \delta}}(x))
  (1 -  \chi(k^{1 -
 \delta}(\theta - \phi)) \DL(k + i \tau, (r, \theta), q(\phi))$
(resp. $\SL$)

\end{itemize}
are  semiclassical Fourier integral operator of order $-1/4$
($\DL$), resp. $-5/4$ ($\SL$), associated to the canonical
relation $\Gamma$.  More precisely, there exist amplitudes such
that
 $$\begin{array}{ll}(i) &  (1 - \chi_{\partial \Omega}^{k^{-1 + \delta}}(x)) ) \DL(k + i
 \tau, x, q(\phi))   \\ \\ & =   (1 - \chi_{\partial \Omega}^{k^{-1 + \delta}}(x)) )
 (k + i
 \tau)^{\frac{1}{2}} A_1 ((k + i \tau), x,
q(\phi))
 e^{i (k + i \tau) |x - q (\phi)|},
 \\ &   \\ (ii) & \chi_{\partial \Omega}^{k^{-1 + \delta}}(r))
 (1 - \chi(k^{1 -
\delta} (\phi - \theta) )) \DL(k + i
 \tau, x, q(\phi))  \\ \\ &  =  (k + i
 \tau)^{\frac{1}{2}} A_2(k + i \tau, r,  \theta, \phi)  e^{i (k + i \tau) |q(\theta) + r\nu_{q(\theta)} - q
(\phi)|}\end{array},$$ with $ A_j(k + i \tau, r, \theta, \phi )
\in S^{0}_{\delta}({\bf T}^2). $ The analogous statements are true
for $\SL$.
\end{prop}

\begin{proof} The proof is similar to that of Proposition
\ref{NSYM}.  Due to the cutoffs, the kernel is in its
semiclassical regime where the WKB approximation is valid and the
orders of the amplitudes can be read off from Lemma \ref{HANKSYM}.
 After substituting $z = |x - q(\phi)|$ we get the
amplitudes:

\begin{itemize}

\item $a_0((k + i \tau)| x - q(\phi)|)$ for $\SL$;

\item $(k + i \tau) a_1 ((k + i \tau)| x - q(\phi)|) \cos \angle ( x
- q(\phi), \nu_{q(\phi)}) $ for $\DL$.

\end{itemize}

We note that the order convention on Fourier integral operators
operating between spaces of unequal dimensions $n_1, n_2$ is that
$k^{n_1/4 + n_2/4} $ times an amplitude of order $0$ on the
stationary phase set defines an operator of order $0$. Thus, the
cutoff  $\DL$ has an order of $1/4$ less than $N_1$. Since the
normal derivative introduces a factor of $(k + i \tau)$ in $\DL$
but not in $\SL$, the order of $\DL$ is one higher than $\SL$.

  To
complete the proof, we need to check, that differentiations in
$k^{-1} D_x, k^{-1} D_{\phi}$ in case (i), resp.   $k^{-1}
D_{\theta}, k^{-1} D_r, k^{-1} D_{\phi}$ in case (ii),  lower the
order of the amplitudes by $k^{- \delta}$. The only difference to
Proposition \ref{NSYM} is that $ \angle q(\phi_2) - q(\phi_1),
\nu_{q(\phi_2)} )$  there is replaced by $\angle ( x - q(\phi),
\nu_{q(\phi)}) = \angle ( q(\theta) + r \nu_{q(\theta)}  -
q(\phi), \nu_{q(\phi)}) $. As before, each derivative at most
increases the number of factors of  $|x - q(\phi)|^{-1}$ by one.
On the support of the cutoff, each such factor counts $k^{1 -
\delta}$. Due to the accompanying factor of $k^{-1}$ we see that
each derivative $k^{-1} D$ decreases the order of the amplitude by
$k^{-\delta}$.

\end{proof}

\subsection{Integral formulas}

In \S 6,  we are going to need an asymptotic formula and remainder
estimate for some integrals involving  Hankel functions. The
reasons for the assumptions on the parameters $a,b$ will be
clarified at that time.  The constraint $1 > \delta > 1/2$ on $h =
k^{-1 + \delta}$ originates in the following.

\begin{prop}\label{INTFORM} Let $a \in (-1 , 1  ),$
let $\Re b \geq 1, \Im b > 0$ and let $R \in {\bf N}$. Also, let
$\chi \in C_0^{\infty}(\R)$ be an even cutoff function, equal to
one near $0$. Then
$$ \begin{array}{l}  \int_{0}^{\infty}  \chi(k^{- \delta}  x) \cos ( a x) H_0^{(1)}(b x) dx  =
\frac{1}{\sqrt{b^2 - a^2}} + O(k^{- R \delta} |a - b|^R) + O(k^{-
R \delta}),
\end{array} $$
where the $O$-symbol is uniform.  The equation  may be
differentiated any number of times in $a$ with the same remainder
estimate.
  \end{prop}

\begin{proof}

The result is suggested by the standard formula
$$\int_{0}^{\infty}  \cos ( a x) H_0^{(1)}(b x) dx  =
\frac{1}{\sqrt{b^2 - a^2}} $$ for the Fourier transform of a
Hankel function \cite{O, AG}. To deal with the cutoff we proceed
as follows.

Since $1_{[1, \infty]}(t) ( t^2 - 1)^{- 1/2}  $ is a tempered
distribution on $\R$, and since $ \chi(k^{- \delta}  x) \cos ( a
x) $ is a Schwartz function, we may replace the Hankel function by
its Fourier integral formula  (\ref{HANKEL}) to obtain
\begin{equation}\begin{array}{lll}
  \int_{0}^{\infty}  \tilde{\chi}(k^{- \delta}  x) \cos ( a x) H_0^{(1)}(b x) dx
 & = & \int_1^{\infty} \{ \int_{0}^{\infty}  \chi(k^{- \delta}  x) \cos ( a x)
 e^{i b x t} dx\}  (t^2 - 1)^{-1/2} dt \\& & \\
&=&  k^{\delta} \sum_{\pm} \int_0^{\infty} \int_{0}^{\infty}
1_{[1, \infty]}(t) ( t^2 - 1)^{- 1/2}   \chi(  x) e^{i k^{
\delta}( b  t \pm a) x} dx. \end{array} \end{equation}

The phase has a (non-degenerate) critical point at $x = 0, t = \pm
a/b$ (since $b \not= 0$). However, with our assumptions on $(a,
b)$,  the critical point lies (slightly) outside of the support
$[1, \infty]$ of the  $dt$ integral, hence the operator
$$ L = k^{- \delta} \frac{1}{b t \pm a} D_x,$$
is well-defined on the support of the integral.  We cannot
integrate by parts in the $dt$ integral due to singularities in
the amplitude, but we can (and will)
 integrate by parts repeatedly in the $dx$ integral  with $L$,
which reproduces the phase.

The first time we integrate by parts with $L$, we
obtain
\begin{equation}\begin{array}{l}  \sum_{\pm} \int_0^{\infty} \int_{0}^{\infty}
1_{[1, \infty]} (t) (( t^2 - 1)^{- 1/2}   \chi(  x)
  \frac{1}{bt \pm a} D_x  e^{i k^{ \delta}( b  t \pm a) x} dxdt \\ \\
=    \int_{0}^{\infty} (( b  t \pm a)^{-1} 1_{[1, \infty]} (t) ( t^2 - 1)^{- 1/2}   dt \\ \\
+  \sum_{\pm} \int_0^{\infty} \int_{0}^{\infty}  1_{[1, \infty]}
(t) (( t^2 - 1)^{- 1/2}   (D_x \chi (  x))
  \frac{1}{bt \pm a}  e^{i k^{ \delta}( b  t \pm a) x} dxdt. \end{array}\end{equation}
We then integrate the second term  by parts repeatedly with $L$
and observe that no further boundary terms are picked up since
$(D_x \chi(  x)) \equiv 0 $ near $x = 0.$ We then break up the
$dt$ integral into $\int_1^2 \{ \cdots\} dt + \int_2^{\infty}
\{\cdots \} dt.$ On the $dt$- integral over $[2, \infty]$, the
factors of $\frac{1}{bt \pm a} $  obtained from the $R$ partial
integrations render  the integral  absolutely convergent and the
factors of $k^{- \delta}$ show that it is of order $O(|k|^{-R
\delta})$  with coefficients independent of $(a,b).$ In the
integral over $[1,2]$ we only have  $|bt \pm a| \geq |b - a|. $
After $R$ partial integrations we therefore obtain an estimate of
$k^{- R \delta} |b - a|^{-R}.$ for the integral. The same kind of
argument applies if we first differentiate the integral $n$ times
in $a$, without changing the remainder estimate.

To complete the proof, we note (with \cite{AG})  that
\begin{equation}   \int_{1}^{\infty} \sum_{\pm} (b  t \pm a)^{-1} ( t^2 - 1)^{- 1/2}     dt
= \frac{1}{\sqrt{b^2 - a^2}}. \end{equation}

\medskip

\end{proof}

\begin{cor}\label{INTFORMCOR} With the same assumptions as above, assume further that  $|a - b| \geq C \; k^{-1 + \delta}.$
 Then
$$ \begin{array}{l}  \int_{0}^{\infty}  \chi(k^{- \delta}  x) \cos ( a x) H_0^{(1)}(b x) dx  =
\frac{1}{\sqrt{b^2 - a^2}} + O(k^{R(1 - 2 \delta)}),
\end{array} $$
where the $O$-symbol is uniform. Thus, for any $\delta$ with $1
> \delta
> 1/2$, the remainder is rapidly decaying. The equation  may be
differentiated any number of times in $a$ with the same remainder
estimate.
  \end{cor}

We will also need a slight extension of this result:

\begin{prop}\label{INTFORMR} With the same notation and
assumptions as above, we have
$$ \begin{array}{l}  \int_{0}^{\infty}  \chi(k^{- \delta}  x) \cos ( a x)  H_0^{(1)} (b \sqrt{r^2 + u^2}) \cos
(a  u) du dx  = - i \frac{e^{- i r \sqrt{b^2 - a^2}}}{\sqrt{b^2 -
a^2}} + O( k^{R (1 - 2 \delta}),
\end{array} $$
where the $O$-symbol is uniform for $a$ in the interval above.
  \end{prop}

  \begin{proof}

The proof proceeds much as before except that we now use the
identity: (\cite{O}, (14.16)-(14.18))
$$\int_0^{\infty} H_0^{(1)} (b \sqrt{r^2 + u^2}) \cos
(a  u) du = - i \frac{e^{- i r \sqrt{b^2 - a^2}}}{\sqrt{b^2 -
a^2}},\;\;\; \mbox{valid if}\;\; b^2 - a^2 > 0.  $$ Details are
left to the reader.

\end{proof}

\subsection{The case of the unit disc}

The only computable example is the unit disc, and the reader might
wish to use it to check the calculations and estimaes.  Due to the
$S^1$ symmetry, $N(k)$ is a convolution operator and there are
simple commutation relations with the layer potentials. We have:

\begin{equation} \label{DISC} \left\{ \begin{array}{ll} N(k, \theta, \phi)
& = k \sum_{n \in \Z}  H^{(1)}_n (k)
J_n'(k  )  e^{i n (\theta - \phi)} \\
&
\\{\mathcal S}\ell (k, (r, \theta), \phi) & = \sum_{n \in \Z}
H^{(1)}_n (k) J_n(k r ) \cos n (\phi - \theta)\\ & \\ {\mathcal
D}\ell (k, (r, \theta), \phi) & = k \sum_{n \in \Z} H^{(1)}_n (k)
J_n'(k r ) \cos n (\phi - \theta)
\end{array} \right. \end{equation}
These layer potential identities follow from the  so-called Graf
addition theorem:
\begin{equation} \label{GRAF} G_0(k, x, y) = \sum_{n \in \Z} H^{(1)}_n (k \max \{|x|, |y|\}) J_n(k \min \{|x|, |y|\} )
\cos n \angle(x, y), \end{equation} which is valid for general
domains.

\section{\label{MICRO} Microlocalizing the trace}

As discussed in the introduction, we are going to use (\ref{POT})
to calculate the wave invariants at $\gamma^r$. In Lemma
\ref{EASYMICRO} we discussed the standard microlocalization of the
trace to $\gamma$. We now prove that we can use a microlocal
cutoff operator  on the boundary as well as  on the domain.
 This obvious sounding statement is not actually obvious
in the present approach since the layer potentials are not
standard semiclassical Fourier integral operators.

Let us first review Lemma \ref{EASYMICRO} in the language of layer
potentials. It asserts that
\begin{equation}\label{CUTOFF2}  \rho * Tr [{\mathcal D}(k + i \tau )
(I + N(k))^{-1 tr} {\mathcal S}(k + i \tau) ^{tr} \circ (1 -
 \chi)\sim 0, \end{equation}
 where $\chi(k)$ is a semiclassical cutoff in
$\Omega$ to a microlocal neighborhood of $\gamma$ which has the
form $\chi (r, y, k^{-1} D_y) $ near the boundary.  Here, $(r, y)$
are Fermi normal coordinates in $\Omega^c$ near the endponts of
$\gamma$ , with $r$ the distance to the boundary. Thus,
differentiations are only in tangential directions.

We now verify that the analogue of (\ref{CUTOFF2}) remains correct
if we use a microlocal cutoff on the boundary. The billiard map
$\beta$ is a cross section of the billiard flow, and in this cross
section a bouncing ball $\gamma$ corresponds to a  periodic orbit
which we denote  $\partial \gamma \in B^* \partial \Omega$ of
period $2$. We choose the boundary parametrization so that $\phi =
0$ is one of the endpoints of the segment (which we also denote
$\partial \gamma$) in $\Omega$. To microlocalize the boundary
operators to the periodic point,  we introduce a semiclassical
pseudodifferential cutoff operator $\cg (\phi, k^{-1} D_{\phi}) $
with complete symbol $\cg(\phi, \eta)$ supported in
$V_{\epsilon}:= \{(\phi, \eta): |\phi|, |\eta| \leq \epsilon\}$.

 \begin{prop} \label{MICRO} Suppose as above that supp $\hat{\rho}  \cap
 Lsp(\Omega) = \{ r L_{\gamma} \}$, and let $\cg (\phi,
 k^{-1} D_{\phi})$ be a semiclassical microlocal cutoff to
the billiard-map orbit corresponding to $\gamma$. Then:
  $$\begin{array}{l} \rho * Tr [{\mathcal D}(k + i \tau )
(I + N(k))^{-1 tr} {\mathcal S}\ell(k + i \tau) ^{tr}    \\ \\
\sim \rho * Tr [{\mathcal D}\ell(k + i \tau ) (I + N(k))^{-1 tr}
\circ \cg(k) \circ {\mathcal S}(k + i \tau) ^{tr}    .\end{array}
$$

 \end{prop}

 \begin{proof}

By (\ref{RTAU}) and by (\ref{INTPI})-(\ref{PIF}), we have:
$$\left\{ \begin{array}{ll} (i) & {\mathcal S}\ell(k + i \tau )^{tr} = \gamma  G_0(k + i \tau )
= \int_0^{\infty} e^{i t (k + i \tau ) }  r_{\Omega 2} E_0(t) dt,\\
\\(ii) &
 1_{\Omega} ({\mathcal D}\ell(k + i \tau) \circ (I + N(k + i \tau ))^{-1 }(x, q)
= \int_0^{\infty} e^{i t (k + i \tau )}  r_{\Omega \nu 2}
E_{D}(t)dt ,\end{array} \right.$$ where  $r_{\Omega \nu} u =
\partial_{\nu} u|_{\partial \Omega}. $   We have subscripted the restriction operators to clarify which
variables they operate on. The composition $ \DL \circ (I +
N(k))^{-1  }  \circ (1 - \cg(k)) \circ \SL^{tr}$ may be written
(with the relevant value of $k$)  as
\begin{equation} \int_0^{\infty} e^{i t (k + i \tau )} \{\int_0^{t} r_{\Omega \nu 2}
E_{D} (t - s) \circ (I - \cg)(y, |D_t|^{-1} D_y)
 r_{\Omega 1}  \circ  E_0(s) ds\} dt. \end{equation}

 We therefore have:
\begin{equation}\label{WFUB2}\begin{array}{l}  \rho *   \DL \circ (I +
N(k))^{-1  }  \circ (1 - \cg(k)) \circ \SL^{tr}
\\ \\
= \int_0^{\infty} \hat{\rho}(t) [\int_0^t r_{\Omega \nu 2} E_{D}
(t - s) \circ (I - \cg)(y, |D_t|^{-1} D_y)
 \circ r_{\Omega 1} E_0(s) ds] e^{i (k + i \tau) t}
 dt .
\end{array}
\end{equation}
The statement of the Proposition is equivalent to:
\begin{equation}
WF [Tr \int_0^t r_{\Omega \nu 2} E_{D} (t - s) \circ (I - \cg)(y,
|D_t|^{-1} D_y) \circ
 r_{\Omega 1}  E_0(s) ds] \cap [rL_{\gamma} - \epsilon, r
 L_{\gamma} + \epsilon] = \emptyset. \end{equation}

We now claim that the integrand
 $$V(s, t) =   Tr \; r_{\Omega \nu 2} E_{D}(t - s) \circ (I -
\cg)(y, |D_t|^{-1} D_y)\circ
 r_{\Omega 1}  E_0(s) $$
 is a smooth function for $t \in (r L_{\gamma} - \epsilon,
 rL_{\gamma} + \epsilon)$ and for $s \in (0,  rL_{\gamma} + \epsilon)$.
 Indeed, the singular support consists of $(s, t)$  such that there
 exists a closed billiard orbit of length $t$ outside a phase space neighborhood
 of $\gamma^r$ which consists of a
 straight line segment of length $s$ from a point $x \in \Omega$ to
 a boundary point $q$, followed by a generalized billiard orbit
of length $t - s$  from $q$ back to $x$. By our assumption on the
length spectrum, the only possible orbit with length $t\in (r
L_{\gamma} - \epsilon, rL_{\gamma} + \epsilon)$ is $\gamma^r$ of
length $r$; but the  cutoff has removed this orbit. Since the
integrand is smooth, the integral over $s \in [0, t]$ determines a
smooth function of $t \in L_{\gamma} - \epsilon, rL_{\gamma} +
\epsilon)$.

\end{proof}

\begin{rem} The Proposition is obvious in  the case of the unit disc, although in this
case it is only natural to  cutoff in the frequency variable since
all radial geodesics are bouncing ball orbits.  Lemma
\ref{EASYMICRO} and Proposition \ref{MICRO} are then  equivalent,
since the cutoff to the right of $\SL^{tr}$ has the form
$\chi(k^{-1} D_{\phi})$ and commutes with $\SL^{tr}$.

\end{rem}

\section{\label{REF} Regularizing the boundary integrals }

The purpose of this section is to analyze  the compositions in
(\ref{BINOMIAL}) with the semiclassical cutoff $\cg$ on the
boundary as semiclassical Fourier integral operators.
 Since the role of the imaginary part of the spectral
parameter  is not important here, we write it simply as $\tau$ and
only substitute $\tau \log k$ when it is needed in \S \ref{M}.

For ease of notation, we write the terms of  (\ref{BINOMIAL}) as
 \begin{equation} N_{\sigma} : =  N_{\sigma(1)} \circ
N_{\sigma(2)} \circ \cdots \circ N_{\sigma(M)}
\end{equation}  and we put
\begin{equation}    |\sigma| = \; \#
 \sigma^{-1} (0) = \; \mbox{the number of $N_0$ factors occurring in}\; N_{\sigma}.\end{equation}

\begin{prop}  \label{NPOWERFINAL}

 \noindent{\bf (A)} Suppose that $N_{\sigma}$  is not of the form $N_0^M$.
 Then for any integer  $R >0$, $N_{\sigma} \circ \cg (k + i \tau)$ may be expressed as the sum
 $$N_{\sigma}
= F_{\sigma}(k, \phi_1, \phi_2) + K_R,  $$  where $F_{\sigma}$ is
a semiclassical Fourier integral kernel of order $- |\sigma|$
associated to $\beta^{M - |\sigma|}$  and a remainder $K_R$, which
is a bounded kernel which is uniformly of order $k^{-R}$.
\medskip

\noindent{\bf (B)} $N_0^M \circ \cg \sim N_{0M} \circ \cg,$ where
$N_{0M}$ is a semiclassical pseudodifferential operator of order
$-M$.

\end{prop}

The proof will be broken up into a sequence of Lemmas. First we
will consider the compositions $N_0 \circ N_1,  N_1 \circ N_0$
without the cutoff $\cg$. Then we consider iterated compositions
with the cutoff. Finally we discuss the special term $N_0^M \circ
\cg. $

\subsection{The compositions $N_0 \circ N_1$}

Here we characterize the composition $N_0 \circ N_1$. Essentially
the same result holds for $N_1 \circ N_0$.

The composed kernel equals
\begin{equation}\label{IKJ}  \begin{array}{l} N_0 \circ N_1 (k + i \tau, \phi_{1},  \phi_{2}) :
 =  (k + i \tau) \int_{{\bf T}}
 \chi(k^{-1 + \delta} (\phi_{1} - \phi_{3} )) (1 - \chi(k^{-1 + \delta} (\phi_{2} - \phi_{3} ))) \\ \\
  H^{(1)}_1
 ((k \mu + i \tau) |q(\phi_{3 }) - q(\phi_{1})|)
\cos \angle (q(\phi_{3 }) - q(\phi_{1}), \nu_{q(\phi_{3 })} )
   N_1 ((k\mu  + i \tau, q(\phi_{2}),   q(\phi_{3})) . \end{array}\end{equation}

The somewhat technical nature of the following lemma is due in
part to the lack of a  cutoff to $\gamma$.

\begin{lem} \label{Onebadpairsym} For any $R \in {\bf N}$, there exists  an amplitude $A_{01} (k + i \tau, \phi_1, \phi_2)$
such that \\

\begin{itemize}

\item (i) $N_0 \circ N_1 (k + i \tau, \phi_{1},  \phi_{2}) ) =
(1 - \chi(k^{1 - \delta} (\phi_1 - \phi_2))
 k^{ (\frac{5}{2} - 3 \delta) }  e^{i k  |q(\phi_{1}) - q(\phi_{2 })|}
 A_{01} (k + i \tau, \phi_1, \phi_2) + K_R(k + i \tau, \phi_1, \phi_2)$;\\

\item (ii)  $A_{01} (k + i \tau, \phi_1, \phi_2) \in S^{0 }_{\delta}({\bf T}^2 \times \R).$\\

\item (iii) $K_R$ is a bounded kernel which is uniformly of order
$k^{-R}$.

\end{itemize}

\end{lem}

\begin{proof}

Using Lemma \ref{HANKSYM} and Proposition \ref{NSYM}, we  rewrite
$N_1$   in terms of phases and amplitudes. The proof is then based
on a change variables and on use of the explicit cosine transform
of the  Hankel function given in Proposition (\ref{INTFORM}) to
evaluate integrals involving the `difficult factor' $N_0(k + i
\tau)$, i.e.   $H^{(1)}_1 ((k \mu + i \tau) |q(\phi_{1}) -
q(\phi_{3} )|)$.

It is convenient to first change variables
$$\phi_3 \to \vartheta = \phi_1 - \phi_3,$$
 and then to change
variables
 $\vartheta \to u$, with:
\begin{equation}  u: = \left\{ \begin{array}{ll} |q(\phi_{3})
- q(\phi_{1})| , & \phi_{1} \geq \phi_{3} \\ & \\
- |q(\phi_{3}) - q(\phi_{1})|, & \phi_{1} \leq \phi_{3}
\end{array} \right.  = \left\{ \begin{array}{ll} |q(\phi_{1}
- \vartheta)
- q(\phi_{1})| , & \vartheta \geq 0 \\ & \\
- |q(\phi_{1} - \vartheta) - q(\phi_{1})|, & \vartheta \leq 0
\end{array} \right.\end{equation} In other words, we change from
the intrinsic distance along $\partial \Omega$ to chordal
distance. Due to the cutoff, the variable $u$ may be taken to
range over  $(- k^{-1 + \delta}, k^{-1 + \delta})$, so the change
of variables is well-defined and  smooth on the support of the
integrand.  The purpose of this change of variables is to simplify
the difficult factor in (\ref{IKJ}):

$$H^{(1)}_1 ((k\mu  + i \tau) |q(\phi_{3} ) - q(\phi_{1})|)
\to H^{(1)}_1 ((k\mu  + i \tau) |u|).$$

Now we consider the other factors. After the change of variables,
we have:

$$\left\{ \begin{array}{l}
 (i)   \cos (\angle q(\phi_{2}) - q( \phi_{2} + \vartheta),
\nu_{q(\phi_{2})})
 \to  |u| K(\phi_{2}, u), \;\; \mbox{with} \; K \; \mbox{smooth in } \; u;\\ \\
 (iii)   N_1 ((k\mu  + i \tau, q(\phi_{2}),   q(\phi_{1} - \vartheta))
 \to      e^{i k |q(\phi_{1 }) - q(\phi_{2} )|}
  e^{i k  u a}
A(k + i \tau,  \phi_{ 1}, \phi_{2}, u), \\ \\
\mbox{where}\;\; A_k \; \mbox{is a symbol in } \;k\; \mbox{ of
order } \;\; 1/2 \;\;\mbox{ and smooth in }\; u. \end{array}
\right.
$$
Here,  \begin{equation} \label{a} a = \sin \vartheta_{1, 2},\;\;\;
\mbox{with} \;\; \vartheta_{1, 2} = \angle (q(\phi_{2}) -
q(\phi_{1}), \nu_{q(\phi_{2})}).
\end{equation}

These statements follow  in a routine way from Proposition
\ref{NSYM} and from the following identities (cf. \cite{EP}):
\begin{equation} \label{Aone}\left\{\begin{array}{l} (i)\; q(\phi) - q(\phi') = (\phi - \phi') T(\phi') - \frac{1}{2}
\kappa (q(\phi))  (\phi - \phi')^2 \nu_{q(\phi')} + O((\phi -
\phi')^3) \\ \\
(ii) \; |q(\phi) - q(\phi') |^2 = (\phi - \phi')^2 + O((\phi - \phi')^4 ) \\ \\
(iii) \; |q(\phi) - q(\phi') | = |\phi - \phi'| + O(\phi - \phi')^3\\ \\
(iv) \; \langle (q(\phi) - q(\phi')), \nu_{q(\phi')} \rangle  =
\frac{-1}{2} (\phi - \phi')^2 \kappa(\phi') + O((\phi-
\phi')^3).\end{array}\right.
\end{equation} Here,   $T(\phi)$ denotes the unit tangent vector
at $q(\phi).$ It follows that $K(u) =  - 1/2 \kappa
(q(\phi_{j_0})) +  {\mathcal O}(|u|^2). $ For further details on
the claims above, we refer to the Appendix \cite{Z6}.

 Taking into account the factor $(k + i \tau)$ in front of the integral, it follows  that the composed kernel (\ref{IKJ})
 can be expressed in the form of Proposition (\ref{Onebadpairsym})
 (i), with
\begin{equation} \label{AK} \begin{array}{l}
A(k + i \tau, \phi_{1}, \phi_{2})  = \int_{- \infty}^{\infty}
\tilde{\chi}(k^{1 - \delta} u) (1 - \chi(k^{1 - \delta}(\phi_2 -
\phi_1 - u)) \\ \\
\times  |u| e^{i k a u} H_1^{(1)}((k + i \tau) |u|) G((k + i
\tau), u, \phi_{1}, \phi_{2}) du,\end{array}
\end{equation} where $G((k + i \tau), u , \phi_{1}, \phi_{2})$ is a symbol
 in $k$ of order $3/2$  and smooth in $u$,   and
where $\tilde{\chi}(k^{1 - \delta} u) = \chi(k^{1 - \delta}
(\phi_{1} - \phi_{3)}).$ The cutoff has been changed slightly
under the change of variables, but for notational simplicity we
retain the old notation for it.

We now  change variables again, $u' = k u$ (and then drop the
prime), to get
\begin{equation}\label{AKU}  A(k + i \tau, \phi_{1},
\phi_{2})  = k^{-2}  \int_{- \infty}^{\infty} \tilde{\chi}(k^{-
\delta} u) (1 - \chi(k^{1 - \delta}(\phi_2 - \phi_1 - k^{-1} u)))
|u| e^{i a u} H_1^{(1)}( b |u|) G(k + i \tau, \frac{u}{k},
\phi_{1}, \phi_{2})    du,
\end{equation}
with $b = 1 + i(\tau/k)$.

\begin{rem}  We note that the entire amplitude now has order $-1/2$ at
least formally. Since we have not yet microlocalized to $\gamma$,
we will temporarily obtain a worse estimate, but in the final step
we will show that this is the correct order. \end{rem}

   Since $|u| \leq k^{\delta}$ on the
support of the cutoff, we have $|\frac{u}{k} | \leq k^{-1 +
\delta}$, and then the Taylor expansion of $G(k + i \tau, u,
\phi_{1}, \phi_{2})$ at $u = 0$ produces an asymptotic series
$$k^{-3/2} \; G(k + i \tau, \frac{u}{k} , \phi_{1}, \phi_{2})  =
 \sum_{n = 0}^{p} k^{-n} u^n G_n (k + i \tau, \phi_{1}, \phi_{2}) + k^{-p} u^p
 R_p(k,  \frac{u}{k} , \phi_{1}, \phi_{2}),$$
where $G_n (k + i \tau,  \phi_{1}, \phi_{2})$ is a symbol of order
$0$ , and where $R_p$ is the remainder,
$$ R_p(k,  \frac{u}{k} , \phi_{1}, \phi_{2}) = \frac{1}{p!}
 \int_0^{1} (1 - s)^{p - 1} \; G^{(p)}(k + i \tau, s
\frac{u}{k} , \phi_{1}, \phi_{2}) ds. $$ Since $G(k + i \tau)$ is
a symbol of order $0$, $ G^{(p)}(k + i \tau, s , \phi_{1},
\phi_{2})$ is uniformly bounded  for $|s| \leq 1$ and hence
\begin{equation} \label{BDD} |R_p(k, \frac{u}{k} , \phi_{1},
\phi_{2})| \leq C_p \;\;\;\; \mbox{for}\;\; |u| \leq k^{\delta}.
\end{equation}

Second, we  Taylor expand the cutoff around $u = 0$ to one order:
$$(1 - \chi(k^{1 - \delta}(\phi_2 - \phi_1 - k^{-1} u)) = (1 - \chi(k^{1 - \delta}(\phi_2 - \phi_1
)) + k^{-  \delta} u S_1(k, \frac{u}{k}, \phi_1, \phi_2),$$ with
\begin{equation} \label{BDD2} S_1(k, \frac{u}{k}, \phi_1, \phi_2)
=  \int_0^{\frac{u}{k}}  \; \chi^{'}(k^{1- \delta} ( \phi_{1} -
\phi_{2} -  s\frac{u}{k})) ds.
\end{equation}

We then write:
$$G (1 - \chi) = (G_p + R_p)( (1 - \chi) + S_p) = G_p (1 - \chi) +
G_p S_p + [R_p ((1 - \chi) + S_p)], $$ where $G_p$ is the $p$th
Taylor polynomial of $G$.
 We claim that the first  terms is a  Fourier integral kernel of the
 type described in (i)-(iii); that the second kernel has this form but multiplied by
 cutoffs like $ \chi^{1} (k^{1 - \delta}(\phi_1 - \phi_2))$ which vanish except in a small
 band $C k^{-1 + \delta} \leq |\phi_1 - \phi_2| \leq 10 C k^{-1 + \delta}$;  and that the remaining two terms define an
 error of the type (iv).

 Consider the first term:
$$ \sum_{n = 1}^R G_n(k + i \tau, \phi_1, \phi_2) \; (1 - \chi(k^{1 - \delta}
 (\phi_1 - \phi_2)) )\; k^{-n - 2 } \int_{- \infty}^{\infty} \tilde{\chi}(k^{- \delta} u)
|u| u^n  e^{i a u} H_1^{(1)}( b |u|)    du.$$   The integral may
be expressed in the form
\begin{equation} \label{nthterm} \begin{array}{l}
\frac{\partial}{\partial b} \frac{\partial^n}{\partial a^n}
\int_{- \infty}^{\infty} \tilde{\chi}(k^{- \delta} u) e^{i  a u}
H_0^{(1)}(b |u|)    du |_{a = \sin \vartheta_{1, 2},  b = (1 +
i\tau/k)}. \end{array} \end{equation} Since the terms in the
finite part of the Taylor expansion of $(1 - \chi(k^{1 -
\delta}(\phi_2 - \phi_1 - k^{-1} u))$ vanish if $|\phi_2 - \phi_1|
\leq k^{-1 + \delta},$  we have $a = \sin \vartheta_{1, 2} \in (-1
+ k^{-1 + \delta}, 1 - k^{-1 + \delta})$. So we may apply
Proposition (\ref{INTFORM}) - Corollary (\ref{INTFORMCOR}) to
evaluate (\ref{nthterm}) asymptotically as
$$\begin{array}{l}  i^{-n}  \frac{\partial}{\partial b}  \frac{\partial^n}{\partial a^n}
 \int_{0}^{\infty} \tilde{\chi}(k^{- \delta} x)
\cos (a x)  H_0^{(1)}(b x)    dx |_{a = \sin
\vartheta_{1, 2}, b =  (1 + i\tau/k)} \\ \\
=   i^{-n} \frac{\partial^n}{\partial a^n} (1 - a^2)^{-3/2} +
O(k^{- n \delta} |(a - 1) + i
\frac{\tau}{k}|^{ - (3 + 2 n)}) \\ \\
=  P_n(a)(1 - a^2)^{-(3/2 + n)} + O(k^{- n \delta} |(a - 1) + i \frac{\tau}{k}|^{ - (3 + 2 n)})\\ \\
 =
(\cos \vartheta_{1, 2})^{-(3 +2 n)} + O(k^{- n \delta} k^{ (3 + 2
n)(1 -  \delta)}).
 \end{array} $$
Here,  $P_n$ as an $n$th degree polynomial which we will not need
to evaluate.

Let  us analyze the order in $k$ of this part of the  amplitude.
We now put in the factor of $k^{-n -2}$ in front of
(\ref{nthterm}).  We observe that $\cos \vartheta_{1, 2}$ can be
as small as $k^{-1 + \delta}$ on the support of the cutoffs in the
integral (which do not include a cutoff to $\gamma$). In view of
the factor
 $(\cos \vartheta_{1, 2})^{-(3 + 2 n)}$, we estimate the
 contribution of the $n$th term to the amplitude as
\begin{equation} \label{COSORD} k^{-n - 2} (\cos \vartheta_{1, 2})^{-(3 + 2n)}
\leq k^{-n - 2 } k^{(3 + 2 n) (1 - \delta)} = k^{n ( 1 - 2 \delta)
+  ( 1 - 3 \delta)}. \end{equation} Each further derivative in
$k^{-1} D$ gives a further factor of $k^{-\delta}.$ Since $\delta
> 1/2$ the terms decrease in order with $n$, and since there are
only a finite number of terms, we conclude that $k^{-3/2} A(k + i
\tau, \phi_1, \phi_2)$ lies in
 $S^{( 1- 3 \delta)}_{\delta}.$ This proves (i) - (ii) modulo
 the remainder estimate.

 \begin{rem} In Lemma \ref{NPOWER}, we will make use of the
 cutoff $\cg$ at the end to eliminate these factors of $\cos
 \vartheta_{1,2}.$ \end{rem}

 We next turn to the term $G_p S_p$, which is the most difficult
 of the remainder terms. As in (\ref{nthterm}) of the previous step, the key point is
 to analyze the integrals
\begin{equation} \label{nthterm2} \begin{array}{l} i^{-n} k^{ - \frac{1}{2}}
\frac{\partial}{\partial b} \frac{\partial^n}{\partial a^n}
\int_0^1 \int_{- \infty}^{\infty} \tilde{\chi}(k^{- \delta} u)
\chi'(k^{1 - \delta} (\phi_1 - \phi_2 - s \frac{u}{k})  e^{i a u}
H_0^{(1)}(b |u|) du ds |_{a = \sin \vartheta_{1, 2}, b = (1 +
i\tau/k)}.
\end{array}
\end{equation}
To deal with this integral, we need to adapt the method of  proof
of Proposition (\ref{INTFORM}) - Corollary (\ref{INTFORMCOR})  to
take into account  the additional cutoff.

As before, we have:
\begin{equation}\begin{array}{l}
  \int_0^1 \int_{0}^{\infty}  \tilde{\chi}(k^{- \delta}  x) \chi'(k^{1 - \delta} (\phi_1 - \phi_2 - s \frac{u}{k})
    \cos ( a x) H_0^{(1)}(b x) dx ds
\\ \\
=  k^{\delta} \sum_{\pm} \int_0^{\infty} \int_0^1
\int_{0}^{\infty} 1_{[1, \infty]}(t) ( t^2 - 1)^{- 1/2}   \chi( x)
\chi'(k^{1 - \delta} (\phi_1 - \phi_2) - s x) e^{i k^{ \delta}( b
t \pm a) x} dx ds dt. \end{array} \end{equation}   The phase is
the same as in Proposition (\ref{INTFORM}), and we again integrate
by parts in $dx$ with $ L = k^{- \delta} \frac{1}{b t \pm a} D_x.$

The first time we integrate by parts with $L$, we obtain
\begin{equation}\begin{array}{l}
  \chi'(k^{1 - \delta}
(\phi_1 - \phi_2)) \int_{0}^{\infty} (( b  t \pm a)^{-1} 1_{[1, \infty]} (t) ( t^2 - 1)^{- 1/2}   dt \\ \\
+  \sum_{\pm} \int_0^1 \int_0^{\infty} \int_{0}^{\infty} 1_{[1,
\infty]} (t) (( t^2 - 1)^{- 1/2}   (D_x [\chi (  x) \chi'(k^{1 -
\delta} (\phi_1 - \phi_2) - s x)])
  \frac{1}{bt \pm a}  e^{i k^{ \delta}( b  t \pm a) x} dxdt ds. \end{array}\end{equation}
We then continue to integrate the second term  by parts repeatedly
with $L$. But unlike the case of Proposition (\ref{INTFORM}), we
do pick up boundary terms each time when  $D_x$ falls on the
second cutoff factor. For instance, the next iterate produces the
boundary term:
\begin{equation}\begin{array}{l}
  k^{-\delta} \chi''(k^{1 - \delta}
(\phi_1 - \phi_2)) \int_{0}^{\infty} (( b  t \pm a)^{-2} 1_{[1, \infty]} (t) ( t^2 - 1)^{- 1/2}   dt \\ \\
+  k^{-\delta} \sum_{\pm} \int_0^1  \int_0^{\infty}
\int_{0}^{\infty} 1_{[1, \infty]} (t) (( t^2 - 1)^{- 1/2}   (D_x^2
[\chi (  x) \chi'(k^{1 - \delta} (\phi_1 - \phi_2) - s x)])\\ \\
  (\frac{1}{bt \pm a})^2  e^{i k^{ \delta}( b  t \pm a) x} dxdt ds. \end{array}\end{equation}
  The boundary terms are similar to the kind in the $G_p \chi$
  terms, except multiplied by $ \chi'(k^{1 - \delta}$ and higher
  derivatives. This establishes the claimed form for the first two
  terms.

  As for the integral, we break  up the $dt$ integral
into $\int_1^2 \{ \cdots\} dt + \int_2^{\infty} \{\cdots \} dt$ as
we did in Proposition (\ref{INTFORM}).  The same observations show
that after $R$ partial integrations, the remainder term is bounded
by $k^{- R \delta} |b - a|^{-R}.$ for the integral. The boundary
terms have the form of $\frac{1}{\sqrt{b^2 - a^2}}$ times the
series $\sum_{j = 1}^R k^{-\delta j} \chi^{(j)} (k^{1 -
\delta}(\phi_1 - \phi_2)). $

 We now sketch the  estimate of the last two remainder terms. There are two
 kinds of terms, those from $G$ which descend in  steps of
 $k^{-1}$ and those of $\chi$ which descend in steps of $k^{-
 \delta}$. If we retain the principal term of $\chi$ and the first $n$ terms of $G$, then
 we obtain a remainder of
\begin{equation}\label{AKU}  k^{ - n - 2} \int_{- \infty}^{\infty} \tilde{\chi}(k^{- \delta} u)
  |u| u^n R_n (u/k, \phi_{1}, \phi_{2})  e^{i  a u} H_1^{(1)}( |u|)    du \end{equation}
  of the $n$th order remainder in the Taylor expansion of the
  amplitude.
  By (\ref{BDD}), we can  bound the $R_n$ and exponential factors by a
  uniform constant $C$.  We can also bound
  the   Hankel factor by $|u|^{-1/2}$ for $|u| \geq 1$; the singularity at $u = 0$ is cancelled by
  the $u^n$.  Thus,
  the integral is bounded by $k^{ - n  - 2} \int_{0}^{k^{\delta}}
  x^{n + 1/2} dx \sim k^{ - n - 2} k^{\delta(n + 3/2)} = k^{n (-1
  + \delta) -1/2}.$ So for $n$ sufficiently large we have an
  arbitrarily small remainder. We can choose $n$ to obtain the
  order stated in (iv).

 \end{proof}

\subsection{Proof of (A)}

By  iterating  the Lemma and  composing with $\cg$ we can improve
Lemma \ref{NPOWER} to the statement in (A) of the Proposition.
First, some more notation. Any  term  $N_{\sigma}$  of
(\ref{BINOMIAL}) can be expressed as a product
\begin{equation}   N_0^{s_r} \circ N_1^{t_r} \cdots  N_0^{s_1} N_1^{t_1},
\end{equation}
of blocks of $N_0$ and $N_1$, with $\sum_{j = 1}^r s_j = |\sigma|$
and $\sum_{j = 1}^r (s_j + t_j) = M$. The number $r = r(\sigma)$
counts the number of blocks of $N_0$. We now compose with $\cg$
and successively eliminate the blocks $N_0^{s_j}$ from right to
left using one factor of $N_1$ to the right.  There is a slight
notational problem since it could happen that $t_1 = 0$, in which
case one should instead use the factor of $N_1$ to the left in
eliminating a block of $N_0$ . Since the process is analogous in
that case, we will assume for simplicity of notation  that $t_1
\not= 0.$

\begin{lem}  \label{NPOWER} For any term $N_{\sigma}$ of (\ref{BINOMIAL})
  except for $N_0^M$, and for any $R \in {\bf N}$, there exist
  $(s_1, \dots, s_r)$ as above and
    amplitudes $A_{j} (k + i \tau, \phi_1,
  \phi_2),$ (with $j = 1, \dots, r$)
  of order $- s_j$
such that:

\begin{itemize}

\item (i) $N_{\sigma} \circ \cg  =  N_r^{t_r - 1} \circ M_{r} \circ N_1^{t_{r - 1} - 1} \circ
\cdots N_1^{t_2 - 1} \circ M_1 \circ N_1^{t_1- 1} \circ \cg +
K_R^j$, where each  $M_j$
 is  a semiclassical Fourier integral kernel of the
form $$M_j(k + i \tau, \phi_1, \phi_2) = k^{1/2}  (1 - \chi(k^{1 -
\delta} (\phi_2 - \phi_2))  e^{i k |q(\phi_{1}) - q(\phi_{2 })|}
 A_j (k + i \tau, \phi_1, \phi_2), $$
 with   $A_j (k + i \tau, \phi_1, \phi_2) \in S^{- s_j }_{\delta}({\bf T}^2 \times \R).$

\item (iv) $K_R$ is a bounded kernel which is uniformly of order
$k^{-R}$.

\end{itemize}

\end{lem}

\begin{proof}  We work from right to left using
the argument of  Lemma \ref{Onebadpairsym} repeatedly to remove
all of the $N_0$ factors in each block.   This can be done because
we only used knowledge of the phase and of the order of the
amplitude to obtain (i) - (iii) of the Lemma . In each such $M_j$
we may choose $R$ large enough so that the remainder for this
block, when composed with the remaining factors of $N_0, N_1$
satisfies (iii). We then  define \begin{equation} F_{\sigma} =
N_r^{t_r - 1} \circ M_{r} \circ N_1^{t_{r - 1} - 1} \circ \cdots
N_1^{t_2 - 1} \circ M_1 \circ N_1^{t_1- 1} \circ \cg .
\end{equation} It is a semiclassical Fourier integral kernel with phase
\begin{equation} {\mathcal L}_{\sigma}(\phi_1, \dots, \phi_{M -
|\sigma|}) = |q(\phi_1) - q(\phi_2)| + \cdots + |q(\phi_{M -
|\sigma|} ) - q(\phi_{M - |\sigma| - 1)}|. \end{equation}

Now let us use the composition with  $\cg$. Since $N_r^{t_r - 1}
\circ M_{r} \circ N_1^{t_{r - 1} - 1} \circ \cdots N_1^{t_2 - 1}
\circ M_1 \circ N_1^{t_1- 1}$ is a semiclassical Fourier integral
operator, its composition with $\cg$ microlocalizes the kernel to
the periodic orbit of $\beta$ corresponding to $\gamma$. That is,
on the support of the cutoff $\cg$, critical points correspond to
Snell paths in which each link points roughly in the direction of
$\gamma$. In calculating the order of the amplitude in Lemma
\ref{Onebadpairsym}, we had to take into account the factors of
$\cos \vartheta_{1, 2}$, which we observed could  be as small as
$k^{-1 + \delta}$. However, upon microlocalizing to the periodic
orbit, there is a uniform lower bound for $\cos \vartheta_{1, 2}$,
and hence these factors are bounded above. As remarked in the
proof of Lemma \ref{Onebadpairsym}, each removal of $N_0$
decreases the order by $1$ on the set where $\cos \vartheta_{1,
2}$ is bounded uniformly below. Thus, the order of $F_{\sigma}$ is
$-|\sigma|$. This improves the order estimate to the statement in
(A) and completes the proof of this part of the Proposition.

\end{proof}

\subsection{The term $N_0^M \circ \cg$}

Our first step in proving (B) of the Proposition is:

\begin{lem} \label{ONEN0} $N_0 \circ \cg $ is a semiclassical pseudodifferential operator of order
$-1$.   \end{lem}

\begin{proof}

The composed kernel equals
\begin{equation}\label{0CHI}  \begin{array}{l} N_0 \circ \cg (k + i \tau, \phi_{1},  \phi_{2}) :
 =  (k + i \tau) \int_{{\bf T}}
 \chi(k^{-1 + \delta} (\phi_{1} - \phi_{3} )) (1 - \chi(k^{-1 + \delta} (\phi_{2} - \phi_{3} ))) \\ \\
  H^{(1)}_1
 ((k \mu + i \tau) |q(\phi_{3 }) - q(\phi_{1})|)
\cos \angle (q(\phi_{3 }) - q(\phi_{1}), \nu_{q(\phi_{3 })} ) \;
  \cg(k, q(\phi_2), q(\phi_3)) d \phi_3, \end{array}\end{equation}
  where
  $$  \cg(k, q(\phi_2), q(\phi_3)) = k \;  \int_{\R} e^{i k (\phi_2 -
  \phi_3) \cdot \eta } \cg (\phi_2, \eta) d \eta $$
  where $\cg(\phi, \eta)$ now denotes (with a slight abuse of notation) the symbol of $\cg$,
   a smooth phase space cutoff to $(0,0), \beta(0,0)$, the
  coordinates of the periodic point of period $2$.

  We make the same change of variables as in Lemma \ref{Onebadpairsym}, which again takes

$$\left\{ \begin{array}{l}

(i)
  H^{(1)}_1 ((k\mu  + i \tau) |q(\phi_{3} ) - q(\phi_{1})|)
\to H^{(1)}_1 ((k\mu  + i \tau) |u|); \\ \\
(ii)  \cos (\angle q(\phi_{2}) - q( \phi_{2} + \vartheta),
\nu_{q(\phi_{2})})
 \to  |u| K(\phi_{2}, u), \;\; \mbox{with} \; K \; \mbox{smooth in } \; u;\\ \\
 (iii) e^{i k (\phi_2 -
  \phi_3) \cdot \eta } \cg (\phi_2, \eta) \to e^{i k (\phi_2 -
  \phi_1)} e^{i k \vartheta   \cdot \eta } = \cg (\phi_2, \eta)  e^{i k (\phi_2 -
  \phi_1)} e^{i k u  \cdot \eta } A(k, u, \eta) \cg  (\phi_2, \eta), \end{array} \right.
$$
where $A(k, u, \eta)$ is an amplitude of order $0$.

Then (\ref{0CHI} ) becomes
\begin{equation}\label{0CHII}  \begin{array}{l}   (k + i \tau) (1 - \chi(k^{-1 + \delta} (\phi_{2} - \phi_{1} )) \int_{\R}
e^{i k (\phi_2 - \phi_1)} a(k + i \tau, \phi_1, \phi_2, \eta) \chi_{\gamma} (\phi_2, \eta)d \eta,\\ \\
\mbox{with} \;\; a(k + i \tau, \phi_1, \phi_2, \eta) = \int_{{\bf
T}}
 \chi(k^{-1 + \delta} u )
 H^{(1)}_1 ((k\mu  + i \tau) |u|)
 |u|  \; e^{i k u  \cdot \eta } A_1 (k, u, \eta)
   d u d \eta, \end{array}\end{equation}
   for another amplitude of order $0$. This is the same kind of
   integral we analyzed in (\ref{AKU}) and  (\ref{nthterm}) of  Lemma
   \ref{Onebadpairsym}, and we obtain the same description of
   $a(k, \phi_1, \phi_2)$ as for $A(k + i \tau, \phi_1, \phi_2)$
   except that the parameter $a$ in (\ref{nthterm}) is now $\eta$.
   Thus, we obtain powers of $(1 - \eta^2)^{-1}$, which blow up on
   the unit sphere bundle of $\partial \Omega$. These points of
   course correspond to tangential (or grazing) rays, and on the support of the cutoff $\cg$
$(1 - \eta^2)^{-1}$ is uniformly bounded. Thus, we obtain the
statement of the Lemma precisely as in Lemma
   \ref{Onebadpairsym} and its improvement Lemma \ref{NPOWER} with the
   cutoff. The change of variables eliminated two powers of $k$,
   leaving an amplitude of order $0$;
   an amplitude of order $0$ defines  a semiclassical  pseudodifferential
   operator of order $-1$.

\end{proof}

To complete the proof of $B$ it suffices to iterate Lemma
\ref{ONEN0}. On each application, we have a new cutoff operator,
but there is no essential change in the argument. This completes
the proof of Proposition \ref{Onebadpairsym}.

\section{\label{RINT} Regularizing the interior integral }

We further need regularize the integrals involving the  outer
factors of ${\mathcal D} \ell (k + i \tau) $ and of $ {\mathcal S}
\ell (k + i \tau)^{tr} $. Since we are taking a  trace, we can
(and will) cycle the factor of ${\mathcal D} \ell (k + i \tau)$ to
the right of ${\mathcal S} \ell^{tr} (k + i \tau) \chi(k) $ to
obtain an operator
\begin{equation}\label{COMPO}  \SL \chi(k) \DL: L^2(\partial \Omega) \to L^2(\partial
\Omega),\;\; \SL^{tr} \chi(k) \DL(q, q') = \int_{\Omega}
G_0(\lambda, q, x) \chi(k^{-1} D_x, x)
\partial_{\nu} G_0(\lambda, x, q') dx. \end{equation}
on  $\partial \Omega$ (or $S^1$, after parametrizing it). Here,
$\chi(k)$ is a semiclassical cutoff to a neighborhood of a
periodic orbit $\gamma$.  For simplicity we will assume that
$\gamma$ is a bouncing ball orbit, although the same method would
apply to a general periodic reflecting ray.

We note that use of the inside/outside duality as in  \cite{Z5}
would in effect make this section unnecessary, at the expense of
forcing a stronger hypothesis on the simplicity of the length
spectrum. Indeed, the integrals over the inside/outside would add
up to the integral in (\ref{COMPO}) but over $\R^2$ instead of
$\Omega$. This integral is easily evaluated to be $N'(k + i
\tau).$  We will now show that the integral over $\Omega$ (with
the cutoff in place) produces a similar kind of semi-classical
Fourier integral operator.

Before stating the precise results, we give  some heuristics on
the composition (\ref{COMPO}). In Proposition \ref{LSYM}, we
described the layer kernels away from their diagonal
singularities. We are now including the latter singularities.  In
addition to $\Gamma$ of (\ref{GAMMA}), the wave front description
now  also includes the relation
\begin{equation} \label{DELTA} \Delta_s = \{(q, \xi, q, \eta):
\xi|_{\partial \Omega} = \eta  \} \subset T^* \Omega \times T^*
\partial \Omega, \end{equation}  which carries the singularities of the kernel.
Intuitively,   $\SL, \DL$ are singular Fourier integral kernels
associated to the union (which we write as a sum) of the two
canonical relations $\Delta_s + \Gamma$, and hence the composition
$\SL^t \circ \chi(k) \circ \DL$ should be associated to the
composition
$$(\Delta_s + \Gamma)^t \circ (\Delta_s + \Gamma) = \Delta_s^t \circ
\Delta_s + \Delta_s^t \circ \Gamma + \Gamma^t \circ \Delta_s +
\Gamma^t \circ \Gamma. $$It is clear that
$$\Gamma^t \circ \Gamma = \Delta_s^t \circ \Gamma =  \Gamma^t \circ \Delta_s = \Gamma_{\beta},$$
the graph of the billiard map. Also, $$\Delta_s^t \circ \Delta_s =
\Delta_{\partial}, $$ the diagonal of $T^* \partial \Omega.$ We
therefore expect the composition to contain these two components,
precisely as $N'(k + i \tau)$ does in the combined inside/outside
case.  The following proposition confirms this. It also produces a
cutoff function, which confirms Proposition \ref{MICRO}.

\begin{prop} \label{COMP} $\SL \chi(k) \DL \sim \cg (k D_{\theta}, \theta) [D_0 + D_1], $
where $D_1$ is a semiclassical Fourier integral operator of order
$-1 $ associated to the billiard map,  where $D_0$ is a
semiclassical pseudodifferential operator of order $-3$, and where
$\cg (k D_{\theta}, \theta)$ is a microlocal cutoff to the
$\beta$-orbit of $(0,0)$.
\end{prop}

\begin{proof} We break up each term as in (\ref{RADCUT}) and then
(\ref{TANCUT}) and analyze each one separately.

 \subsubsection{The most regular terms}

 The regular  terms are those of the form
\begin{equation} \label{I}   \int_{\Omega}  (1 - \chi_{\partial
\Omega}^{k^{-1 + \delta}}(x)) )^2 G_0(\lambda, q, x) \chi(k^{-1}
D_x, x)
\partial_{\nu} G_0(\lambda, x, q') dx. \end{equation}
or where the cutoff has the form $ (1 - \chi_{\partial
\Omega}^{k^{-1 + \delta}}(x)) ) \chi_{\partial \Omega}^{k^{-1 +
\delta}}(x))$.

\begin{lem} \label{REG}  The regular  terms of the form  (\ref{I}) define  Fourier integral
kernels of the form $$k^{1/2}  e^{i (k + i \tau) |q (\phi) -
q(\phi')|} \cg(q, q', \frac{q - q'}{|q - q'|})  A(k, \tau,\phi,
\phi')$$ where $A$ is a semiclassical amplitude  of order $- 1$,
and where $\cg$ is a cutoff to $\gamma$. Thus, (\ref{I}) defines a
semiclassical Fourier integral operator of order $-1$ associated
to $\beta$.
\end{lem}

\begin{proof}

Each factor of the product is described by  Proposition
\ref{LSYM}. Further, we  may take to be the  product of a
frequency cutoff
 $\chi(k^{-1} D_{\theta})$ and a spatial cutoff $\psi_{\gamma}(x)$
 to a strip around $\gamma$.
 Thus, the integral has the form:
 \begin{equation} \label{Ia}  (k + i\tau) \int_{\Omega}  (1 - \chi_{\partial
\Omega}^{k^{-1 + \delta}}(x)) )^2 e^{i ((k + i \tau) [ |x - q| +
|q' - x|]} \psi_{\gamma}(x) \chi(\frac{x - q'}{|x - q'|}) A(k, x,
q, q') dx
\end{equation} where $A$ is a semiclassical smooth amplitude of
order $-1$.  The stationary phase set is defined by
$$\{x : d_x |x - q| = - d_x |x - q'| \iff \frac{x - q}{|x - q|} = -  \frac{x - q'}{|x - q'|}\}. $$
First, we see that $x \in \overline{q q'}$ (the line segment
between $q, q'$).    The ray $\overline{q q'}$ is constrained by
the cutoff to point in the nearly vertical direction in the strip
containing $\gamma$, hence $|q - q'| \geq C
> 0$ on the support of the cutoff.  (Note that the phase is the
sum, not the difference, of the distances since we are taking the
transpose, not the adjoint, of ${\mathcal S}\ell(k + i \tau)$. For
the adjoint, the critical point equation would force $q = q'$).

The ray $\overline{q q'}$ is thus a critical manifold of the
phase, and the phase equals $|q - q'|$ along it. To prove Lemma
\ref{REG}, we show that this critical manifold is non-degenerate
and determine the amplitude by stationary phase. We choose
rectangular coordinates $(s, t)$ oriented so that the $t$ axis is
the ray $\overline{q q'}$ and so that the $s$-axis is orthogonal
to it. Then $q = (0, b), q' = (0, b')$ for some $b, b' \in \R$ and
the phase may  be written $\Psi =(s^2 + (t - b)^2)^{1/ 2} + (s^2 +
(t - b')^2)^{1/ 2}$. Hence, $\Psi'_s = s [(s^2 + (t - b)^2)^{- 1/
2} + (s^2 + (t - b')^2)^{- 1/ 2}]$ and on the stationary phase set
$s = 0$ the second derivative is simply $\Psi''_{ss}(0, t)  = [((t
- b)^2)^{- 1/ 2} + ((t - b')^2)^{- 1/ 2}]$. It is clear that
$\Psi''_{ss}(0, t) $ is uniformly bounded below and since $|b -
b'| \geq C$ on the support of the cutoff it is also uniformly
bounded above. Stationary phase introduces a factor of $k^{-1/2}$
so the resulting amplitude is a product of $k^{1/2}$ with an
amplitude of order  $-1$ and hence the composition defines a
semiclassical Fourier integral operator of order $- 1$.

\end{proof}

\subsubsection{Singular terms}

Therefore we only need to consider the integral \begin{equation}
\label{III} \int_{\Omega} [ \chi_{\partial \Omega}^{k^{-1 +
\delta}}(x)) ]^2 G_0(\lambda, q, x) \chi(k^{-1} D_x, x)
\partial_{\nu} G_0(\lambda, x, q') dx.
\end{equation}  With no essential loss of generality, we redefined the cutoff
to remove the square. We write $x = (r, \theta)$ and then break up
the integral into the sum of four terms corresponding to the
cutoffs:

 \begin{itemize}

\item (i) $   \chi_{\partial \Omega}^{k^{-1 +
\delta}}(r)) \chi(k^{1 -
 \delta}(\theta - \phi)  \chi(k^{1 -
 \delta}(\theta - \phi'); $

\item (ii)  $
  \chi_{\partial \Omega}^{k^{-1 + \delta}}(r))
  (1 -  \chi(k^{1 -
 \delta}(\theta - \phi))  \chi(k^{1 -
 \delta}(\theta - \phi');$

\item (iii) $    \chi_{\partial \Omega}^{k^{-1 +
\delta}}(r)) \chi(k^{1 -
 \delta}(\theta - \phi)  (1 - \chi(k^{1 -
 \delta}(\theta - \phi')); $

\item (iv) $
  \chi_{\partial \Omega}^{k^{-1 + \delta}}(r))
  (1 -  \chi(k^{1 -
 \delta}(\theta - \phi)) (1 - \chi(k^{1 -
 \delta}(\theta - \phi');$

 \end{itemize}

 \subsubsection{Codimension zero case} This refers to case (iv).
 There are no singularities in the integrand due to the cutoff. We
 can use the WKB approximation in each Green's function and apply the cutoff operator in the
 smooth variables $\phi, \phi'$ to  obtain
 integrals of the form:
 \begin{equation} \label{III (iv)} \begin{array}{l} \int_0^{\epsilon} \int_0^{2 \pi}
  \chi_{\partial \Omega}^{k^{-1 + \delta}}(r))
  (1 -  \chi(k^{1 -
 \delta}(\theta - \phi)) (1 - \chi(k^{1 -
 \delta}(\theta - \phi') ) \\ \\
 \chi(\nabla_x |x(r,\theta) - q')   e^{i ((k + i \tau) [ |x(r, \theta) - q| +
|q' - x(r, \theta) |]} A(k, r, \theta, q, q') r dr
d\theta.\end{array} \end{equation} There is no essential
difference to the regular terms in (\ref{I}). We therefore have
the same result.

\begin{lem} Integral (\ref{III (iv)}) defines a Fourier integral kernel
of order $-1$ of the same form as Lemma (\ref{REG}).
\end{lem}

\subsubsection{Codimension one case}

This applies to cases (ii) - (iii), which are quite similar
although not identical. We do case (ii); case (iii) is similar.

 We have:
\begin{equation} \label{III(ii, iii)} (ii) =  \begin{array}{l} \int_0^{\epsilon} \int_0^{2\pi}
  \chi_{\partial \Omega}^{k^{-1 + \delta}}(r))
  (1- \chi(k^{1 -
 \delta}(\theta - \phi'))  \chi(k^{1 -
 \delta}(\theta - \phi)   \\ \\ G_0(\lambda, q(\phi), r, \theta) \chi(k^{-1} D_x, x)
\partial_{\nu}  G_0(\lambda,
(r, \theta), q(\phi')) r dr d\theta \end{array} \end{equation}

\begin{lem} Integrals (\ref{III(ii, iii)}) define  Fourier integral kernels of the
type $$k^{1/2} \cg(q, q', \frac{q - q'}{|q - q'|})  e^{i (k + i
\tau) |q(\phi) - q(\phi')|} A(k, \tau, q(\phi), q(\phi')),$$ with
$A$ an amplitude of order $-2$, i.e. they are semiclassical
Fourier integral operators of order $-2$ associated to $\beta$.
\end{lem}

\begin{proof} The proof is reminiscent of that of Lemma \ref{Onebadpairsym} but is somewhat more
complicated. For the sake of brevity, we concentrate on the new
details and do not discuss the error estimate, which is similar to
that in the proof of Lemma \ref{Onebadpairsym}).

We substitute  the WKB approximation for $
\partial_{\nu} G_0(\lambda, (r, \theta),
q(\phi'))$ (see Proposition \ref{LSYM}) but not for $G_0(\lambda,
q(\phi), r, \theta) $ and apply the  cutoff operator to the WKB
expression. The integral of concern is thus:
\begin{equation}\label{CODIMONE0}  \begin{array}{ll} (\ref{CODIMONE0}) :
=&
  (k + i \tau) \;\int_0^{\epsilon} \int_{{\bf T}}  \chi(k^{1 - \delta}  r)
  \chi(k^{1 - \delta} (q(\phi_0) + r \nu_{q(\phi_0)} - q(\phi_1)))\\ &
  \\&
  H^{(1)}_0 ((k  + i \tau)  |q(\phi) + r \nu_{q(\phi)} -
  q(\theta)|) \; \chi(\frac{q(\theta) + r \nu_{q(\theta)} -
q(\phi')}{|q(\theta) + r \nu_{q(\theta)} - q(\phi')|} \cdot
T_{\theta} ))\\ & \\ &
  A (k + i \tau,  \theta, r, \phi')  e^{i (k + i \tau) | q(\theta) + r
\nu_{q(\theta)} - q(\phi') )|} (1 - \chi(k^{1 - \delta} |q(\theta)
- q(\phi') |) r   d \theta dr,\end{array}
\end{equation}
where $A \in S^{-1/2}_{\delta}$.

We cannot use the WKB expression for the  factor  $H^{(1)}_0 ((k +
i \tau) \mu |q(\theta) + r \nu_{q(\theta)} - q(\phi)|)$, so we
deal with it by
  changing variables and explicitly integrating. The  change of variables is
  given by $(r,  \phi)  \to (r, u)$, with :
\begin{equation} \label{uone}  u(\theta): = \left\{ \begin{array}{ll} (|q(\theta) + r \nu_{q(\theta)}
 - q(\phi)|^2 - r^2)^{1/2} , & \phi \geq \theta \\ & \\
- (|q(\theta) + r \nu_{q(\theta)} - q(\phi)|^2 - r^2)^{1/2} , &
\theta \geq \phi  \end{array} \right. \end{equation} Here, $\phi -
\theta$ and $u$ range only over  $(- k^{-1 + \delta}, k^{-1 +
\delta})$. We claim that  $u(\phi)$ is smooth and invertible with
uniform bounds on derivatives as $r$ varies on $[0, \epsilon_0].$

For the sake of brevity, we only verify this in the basic  case of
a circle of radius $a$ and refer to   \cite{Z6} for the routine
extension to the general case,
 which only requires showing  that the quadratic
approximation is sufficient to determine the smoothness of the
change of variables.  We recall here that we have microlocalized
the integral to $\gamma$ and that $\partial \Omega$ has
non-vanishing curvature at the reflection points of $\gamma$.  We
then have:
$$(|q(\theta) + r \nu_{q(\theta)} - q(\phi)|^2 - r^2) = 2 (1 - r) (1 - \cos (\phi - \theta)).$$
For $|\theta - \phi| \leq - k^{-1 + \delta}, (1 - \cos (\theta -
\phi))$ has a smooth square root, given by the stated formula.

We now consider the effect of this change of variables on the
remaining factors. For the exponential factor, we have
$$\left\{ \begin{array}{l}
 e^{i (k + i \tau) | q(\theta) + r
\nu_{q(\theta)} - q(\phi') )|}  \to e^{i k  |q(\phi) - q(\phi'
)|}
  e^{i k   (u a + r \sqrt{1 - a^2}) } A_1(k,  \phi, \phi', u, r),
 \end{array}\right.$$
where $a = \sin \angle(q(\phi) - q(\phi'), \nu_{q(\phi)}) $ and
where $A_1(k, \mu, \phi, \phi', u, r)$ is a polyhomogeneous symbol
of order $0$  in $k$. We note that $F(0,0) = 1, G(0,0) = - 1/2
\kappa(q(\phi_1)).$ Also, the cutoff transforms as:
$$\cg(\frac{q(\theta) + r \nu_{q(\theta)} -
q(\phi')}{|q(\theta) + r \nu_{q(\theta)} - q(\phi')|} \cdot
T_{\theta} )) \to \cg (\frac{q(\phi) - q(\phi')}{|q(\phi) -
q(\phi')|} \cdot T_{\theta} )) K(k, r, q, \phi, \phi',\theta),
$$
where $K(k, r, q, \phi, \phi',\theta)$ is a smooth semiclassical
amplitude whose leading order term equals $1$.

After multiplying together  these amplitudes, and rescaling  the
variables $u \to k u, r \to k r$. In the regime $ \phi - \theta =
O(k^{-1 + \delta})$, we obtain we obtain an amplitude  $A(k, \phi,
\phi')$ such that:
 $ (\ref{CODIMONE0})  =  \cg(q, q', \frac{q - q'}{|q - q'|})  e^{i k |q(\phi) - q(\phi')|} A(k, \phi,\phi');$
with
\begin{equation} \label{COV} \begin{array}{ll} A (k, \phi, \phi')
=& k^{-1} \int_0^{\infty} \int_{- \infty}^{\infty} \chi (k^{-
\delta} u) \chi (k^{ - \delta} r) e^{i  \mu (u a + r \sqrt{1 -
a^2})}\;  \\ & \\
& H_0^{(1)}(b \sqrt{u^2 + r^2})
 A_1 ((k + i \tau), u, \phi, \phi',r, \mu)) )   du dr,\end{array}
\end{equation}
where $b = \mu (1 + i \tau/k)$ and where $A_1$ is an amplitude of
order $-1/2$. As a check on the order,  we note that the normal
derivative contributed a factor of $(k + i \tau)$ and the change
of variables put in $k^{-2}$, leading to the stated power of $k$.

To complete the proof, we need to show that  $
 A \in  S^{-1/2}_{\delta}({\bf T}\times \R  ).$
We analyse this integral (\ref{COV}), working by induction on the
Taylor expansions of $ A_1$ in the integral
$$\int_0^{\infty} \int_{- \infty}^{\infty} \chi (k^{-
\delta} u) \chi (k^{ - \delta} r) e^{i  (u a_1 + r a_2)}\;
 H_0^{(1)}(b \sqrt{u^2 + r^2})    du dr.$$

 We recall from  Proposition (\ref{INTFORMR}) that:
$$J(r; a, b) := \int_0^{\infty} H_0^{(1)} (b \sqrt{r^2 + u^2}) \cos
(a  u) du = - i \frac{e^{- i r \sqrt{b^2 - a^2}}}{\sqrt{b^2 -
a^2}},\;\;\; \mbox{valid if}\;\; b^2 - a^2 > 0.  $$ Taylor
expanding the $\cos \angle q(\theta) + r \nu_{q(\theta)} -
q(\phi), \nu_{q(\phi)} )$ leads to  the following integrals:
$$\label{INTCOV} \left\{ \begin{array}{l}  \int_{- \infty}^{\infty}
   e^{- i   a u }  H_0^{(1)}(b  \sqrt{u^2 + r^2})
 \frac{ r   }{\sqrt{u^2 + r^2}} d u = \frac{\partial}{b \partial r} J_k(r; a, b)  \\ \\
 \int_{- \infty}^{\infty}  e^{- i   a u }  H_1^{(1)}(b  \sqrt{u^2 + r^2})
 \frac{ u^2   }{\sqrt{u^2 + r^2}} d u =[ \frac{\partial}{\partial
 b} - r \frac{\partial}{b  \partial r}]   J_k(r; a, b)
\end{array}\right. $$
As in Proposition (\ref{INTFORM}), the cutoff factor gives lower
order terms in $k$.

 We then integrate in $dr$. The basic integrals are :
 \begin{equation} \label{INTCOV} \left\{ \begin{array}{lll}  \int_{0}^{\infty}
   e^{ i   a_2 r } \chi (k^{ - \delta} r)  \frac{\partial}{b \partial r}\frac{e^{- i r \sqrt{b^2 - a^2}}}{\sqrt{b^2 -
a^2}} dr & = &    \frac{-i }{b}  \int_{0}^{\infty}
   e^{ i   a_2 r } \chi (k^{ - \delta} r) e^{- i r \sqrt{b^2 - a^2}} dr \\& &  \\
 \int_{0}^{\infty}  e^{ i   a_2 r }\chi (k^{ - \delta} r)  [ \frac{\partial}{\partial
 b} - r \frac{\partial}{b  \partial r}]  \frac{e^{- i r \sqrt{b^2 - a^2}}}{\sqrt{b^2 -
a^2}} dr&  = &   \frac{\partial}{\partial
 b} \int_{0}^{\infty}  e^{ i   a_2 r }\chi (k^{ - \delta} r)   \frac{e^{- i r \sqrt{b^2 - a^2}}}{\sqrt{b^2 -
a^2}} dr \\ & \\ & + & \frac{1}{b} \int_{0}^{\infty}  e^{ i a_2 r
}\chi (k^{ - \delta} r)  r e^{- i r \sqrt{b^2 - a^2}}  dr
\end{array}\right. \end{equation}

To leading order in $k$ the exponentials $ e^{ i   a_2 r }$ and
$\frac{e^{- i r \sqrt{b^2 - a^2}}}{\sqrt{b^2 - a^2}}$ cancel, and
we find the integrals grow at the rates $k^{-2 + \delta} \cos
\vartheta_{M 1}^{-1}.$ By differentiating in $(a_1, a_2)$, can
obtain any term in the Taylor expansion of $F$ with remainder
estimate. Putting in the higher order terms in the Taylor
expansions of $F$ just adds lower order terms in $k^{-1}$,
producing a symbol expansion as in the previous cases. Due to the
cutoff factor $\cg (\frac{q(\phi) - q(\phi')}{|q(\phi) -
q(\phi')|} \cdot T_{\theta} ))$, the factors of  $\cos
\vartheta_{M 1}^{-1}$ are bounded above.    The details are now
similar to those in the proof of Lemma \ref{Onebadpairsym}.

\end{proof}

\subsubsection{Codimension two}

The final integral we must consider is

\begin{equation} \label{IV} \begin{array}{l}  \int_0^{\epsilon} \int_{S^1}   \chi_{\partial
\Omega}^{k^{-1 + \delta}}(x)) \chi(k^{-1 + \delta}(\theta -
\phi))) \chi(k^{-1 + \delta}(\theta - \phi')))  G_0(\lambda, q, x)
 \chi(k^{-1} D_x, x)
\partial_{\nu} G_0(\lambda, x, q') r dr d\theta.\end{array}
\end{equation}

In this case, we cannot use the WKB expansion for either Green's
function and must use integral formulas for products of Hankel
functions. This is the most complicated case, and it is the one
producing the pseudodifferential term in Proposition \ref{COMP}.

\begin{lem} Integral (\ref{IV}) defines a semiclassical pseudodifferential  operator
of order $- 3$  with kernel of  the form $a(\phi, D_{\phi}) \circ
\cg$, i.e. with kernel of the form
$$ k \; \int_{\R} \cg(\phi, \eta)) e^{i k
  \eta  (\phi - \phi')} A_k(\phi, \phi', \eta) d \eta,$$
  where $A $ is a semiclassical amplitude of order $- 3$.
\end{lem}

\begin{proof}
The integral we are considering is:
\begin{equation}\label{CODIMTWO}  \begin{array}{ll}  (\ref{IV})
=& (k + i \tau)
  \int_0^{\epsilon} \int_{{\bf T}}  \chi(k^{1 - \delta}  r)
  \chi(k^{1 - \delta} (q(\phi) + r \nu_{q(\theta)} - q(\phi)))\\ &
  \\&  \chi(k^{1 - \delta} (q(\theta) + r \nu_{q(\theta)} -
  q(\phi')))\;
  H^{(1)}_0 ((k  + i \tau)  |q(\theta) + r \nu_{q(\theta)} -
  q(\phi)|)\\ & \\ &  \chi(k^{-1} D_{\theta}, \theta, r)
H^{(1)}_1 ((k + i \tau)  |q(\theta) + r \nu_{q(\theta)} -
q(\phi')|)  \angle q(\theta) + r \nu_{q(\theta)} - q(\phi'),
\nu_{q(\theta)} )  r d \theta dr .\end{array} \end{equation}

We substitute the Fourier integral formula for $\chi(k^{-1}
D_{\theta}, \theta, r)$ to obtain:
\begin{equation}\label{CODIMTWO}  \begin{array}{ll} (\ref{IV})  :
=& (k + i \tau)
  \int_0^{\epsilon} \int_{{\bf T}} \int_{{\bf T}}  \int_{\R}  \chi(k^{1 - \delta}  r)
  \chi(k^{1 - \delta} (q(\phi) + r \nu_{q(\theta)} - q(\phi)))\\ &
  \\&  \chi(k^{1 - \delta} (q(\theta) + r \nu_{q(\theta)} -
  q(\phi')))\;
  H^{(1)}_0 ((k  + i \tau)  |q(\theta) + r \nu_{q(\theta)} -
  q(\phi)|)\\ & \\ &  \chi(p_{\theta}, \theta, r) e^{i k
  p_{\theta} (\theta - \theta')}
H^{(1)}_1 ((k + i \tau)  |q(\theta') + r \nu_{q(\theta')} -
q(\phi')|) \\ & \\ &  \angle q(\theta') + r \nu_{q(\theta')} -
q(\phi'), \nu_{q(\theta')} )  r d \theta dr d \theta' d p_{\theta}
.\end{array}
\end{equation}

We then   make  the  change of variables
 $(\theta, \theta', r) \to (u, u', r)$ defined  in (\ref{uone})
with respect to the pairs $(\theta; \phi, r) \to (u, r)$ and
$(\theta'; \phi', r) \to (u', r)$.  Thus,
\begin{equation} \label{uone}  u: = \left\{ \begin{array}{ll}
 (|q(\theta) + r \nu_{q(\theta)} - q(\phi)|^2 - r^2)^{1/2} , & \phi \geq \theta \\ & \\
- (|q(\theta) + r \nu_{q(\theta)} - q(\phi)|^2 - r^2)^{1/2} , &
\theta \geq \phi  \end{array} \right. \end{equation} while $u'$ is
defined similarly with $(\theta', \phi')$ replacing $(\theta,
\phi)$.   We note that $u \sim (\theta - \phi)$ since
$$u^2 = |q(\theta) - q(\phi)|^2 - 2 r \kappa(\phi) |q(\theta) -
q(\phi)| (\theta - \phi) \cdot  \sin \angle \frac{q(\theta) -
q(\phi)}{|q(\theta) - q(\phi)|}, \nu_{q(\theta) })+ \cdots$$

This gives
\begin{equation}\label{CODIMONE0CH}  \begin{array}{ll}
(\ref{CODIMTWO}) =&
  \int_0^{\epsilon} \int_{-\infty}^{\infty} \int_{-\infty}^{\infty}  \int_{-\infty}^{\infty}  \chi(k^{1 - \delta}  r)
  \chi(k^{ 1 - \delta} u)
   \chi(k^{1 - \delta} u')  \\ &
  \\& \chi(p_{\theta}, \theta(u, r, \phi), r) e^{i k
  p_{\theta} (\theta(u, \phi, r) - \theta'(u', \phi', r))}\\ & \\
  &
  H^{(1)}_1 ( (k + i \tau)  \sqrt{u^2 + r^2})
H^{(1)}_0 ((k + i \tau) \sqrt{(u')^2 + r^2})  \\ & \\ & (\frac{r
F(u',r) + (u')^2 G(u', r)}{\sqrt{(u')^2 + r^2}} ) B_k(u, u', r,
\phi, \phi') d u du' dr dp_{\theta} .\end{array}
\end{equation}
Here, $B_k$ is a smooth, polyhomogeneous amplitude. We recall that
$\chi(p_{\theta}, u(\theta, r, \phi), r) $ is compactly supported
in $p_{\theta}$ so that there is no problem of convergence of the
$d p_{\theta}$ integral.

We then change to scaled  variables  $ k^{-1} u,  k^{-1} u',
k^{-1} r$ and expand
\begin{equation} \theta (k^{-1} u, k^{-1} r, \phi) = k^{-1} u +
\phi + \cdots
\end{equation}
to obtain
\begin{equation} \begin{array}{ll} (\ref{CODIMTWO}) =& \chi(p_{\theta}, 0, \phi)) e^{i k
  p_{\theta} (\phi - \phi')} A_k(\phi, \phi', p_{\theta}),
  \end{array}
  \end{equation}
  where
$$ \begin{array}{l}  A_k(\phi, \phi', p_{\theta})=  k^{-3} \int_0^{\epsilon} \int_{-\infty}^{\infty} \int_{-\infty}^{\infty}   \chi(k^{- \delta}  \rho)
  \chi(k^{ - \delta} u)
   \chi(k^{ - \delta} u' )\\  \\
   \chi'(p_{\theta}, \theta (k^{-1} u^*, k^{-1} r^*, \phi), k^{-1} r) e^{i
  p_{\theta} (u - u')}\\  \\
  H^{(1)}_1 ( b_1  \sqrt{(u')^2 + r^2}) (\frac{r   }{\sqrt{(u')^2 + r^2}} ) H^{(1)}_0 (b_0 \sqrt{u^2 + r^2})
  (\frac{r
F(u',r) + (u')^2 G(u', r)}{\sqrt{(u')^2 + r^2}} ) \\ \\ B_k(k^{-1}
u,
  k^{-1} r, \phi, \phi', k^{-1} u') d
u du' dr. \end{array} $$

We would like to prove that $A_k$ is a symbol.  As above, we work
by induction on the Taylor expansion of $F, G, B_k $. Polynomials
in $(u, u', r)$ may be expressed as sums of derivatives with
respect to parameters $(a_0, a_1, b_0, b_1)$ at $a_0 = a_1 = $ of
the integrals
$$ \begin{array}{l} k^{-3} \int_0^{\epsilon} \int_{-\infty}^{\infty} \int_{-\infty}^{\infty}   \chi(k^{- \delta}  \rho)
  \chi(k^{ - \delta} u)
   \chi(k^{ - \delta} u' )\\  \\
   \chi(p_{\theta}, \theta (k^{-1} u, k^{-1} r, \phi), k^{-1} r)  e^{i
  p_{\theta} (u - u')}\\  \\
  H^{(1)}_1 ( b_1  \sqrt{(u')^2 + r^2}) (\frac{r   }{\sqrt{(u')^2 + r^2}} ) H^{(1)}_0 (b_0 \sqrt{u^2 + r^2}) d
u du' dr, \end{array} $$ or the same with $(\frac{(u')^2
}{\sqrt{(u')^2 + r^2}} )$ replacing $(\frac{r   }{\sqrt{(u')^2 +
r^2}} )$.

We now do the $du_1$ and $du_0$ integrals first as in the $\phi_0
\sim \phi_1$ case.  The zeroth order terms in the usual Taylor
expansions produce  $dr$ integrals of a form similar
 to (\ref{INTCOV}):
 \begin{equation} \label{INTCOVTWO} \left\{ \begin{array}{l}   \frac{1}{\sqrt{b_0^2 -
a_0^2}}  \frac{-i }{b_1}  \int_{0}^{\infty}
   e^{ i   a_2 r } \chi (k^{ - \delta} r) e^{- i r [\sqrt{b_1^2 - a_1^2} } dr \\  \\
   \frac{1}{\sqrt{b_0^2 -
a_0^2}}  \frac{\partial}{\partial
 b_1} \int_{0}^{\infty}  e^{ i   a_2 r }\chi (k^{ - \delta} r)   \frac{e^{- i r \sqrt{b_1^2 - a_1^2}}}{\sqrt{b_1^2 -
a_1^2}} dr \\  \\   \frac{1}{\sqrt{b_0^2 - a_0^2}}
\frac{k^{-3}}{b} \int_{0}^{\infty}  e^{ i a_2 r }\chi (k^{ -
\delta} r)  r e^{- i r \sqrt{b_1^2 - a_1^2}} dr
\end{array}\right.\end{equation}
The rest proceeds as in the proof of Lemma \ref{Onebadpairsym}. We
conclude that

\end{proof}

\section{\label{M} Tail estimate}

In this section we provide the remainder estimate. As discussed in
the Introduction, the remainder comes from the Neumann series:
\begin{equation} (I +  N(k + i \tau))^{-1} = \sum_{M = 0}^{M_0} (-1)^M N(k + i \tau)^M \;
+\; {\mathcal R}_{M_0},\;\;\; {\mathcal R}_{M_0} = N(k + i
\tau)^{M_0 + 1} \; (I +  N(k + i \tau))^{-1}. \end{equation}

It is now important to take the  imaginary part to be of the form
$\tau \log k.$ We further specify $\rho$ to have the following
properties: $\hat{\rho}(t) = \hat{\rho}_0( t - L)$ where
$\hat{\rho}_0 \in C_0^{\infty}(\R)$ is non-negative and supported
in an $\epsilon$-interval around $0$.

\begin{lem} \label{MASPECT}  For any $R$, there  exists $M_0 = M_0(R)$ and
$$ Tr \rho *  {\mathcal D} \ell (k + i \tau \log k) {\mathcal
R}_{M_0}(k + i \tau \log k)  {\mathcal D} \ell (k + i \tau \log
k)^{tr} \chi(k)  = O(k^{-R}). $$
\end{lem}

\begin{proof}

By Proposition \ref{MICRO}, we can insert the cutoff $\cg$ into
the trace to obtain:
\begin{equation} \begin{array}{l} Tr 1_{\Omega} \rho \; * \; {\mathcal D} \ell (k + i \tau \log k)
(I + N(k + i \tau \log k)^{-1} N^{M_0} \cg(k)  \\ \\ \circ (k + i
\tau \log k) {\mathcal S} \ell ^t(k + i \tau \log k) \chi(k)
\end{array} \end{equation} where ${\mathcal S}^t(k, q, x) =
r_{\Omega} G_0(k, x, y)$ is the transpose of the single layer
potential. The cutoffs may be inserted at various points, and for
simplicity of notation we often write only one in (the other can
then also be inserted).

 We regard the trace as a Hilbert-Schmidt inner product
for Hilbert-Schmidt operators from $L^2(\partial \Omega)$ to
$L^2(\Omega)$.  For ease of notation, we  the remainder at step
$M_0 + 1$ and the trace is then
\begin{equation}  \rho \; * \; \langle {\mathcal D} \ell (k + i \tau \log k)
(I + N(k + i \tau \log k)^{-1} N(k + i \tau \log k) ,  [ N^{M_0}
(k + i \tau \log k) \cg {\mathcal S} \ell ^t(k + i \tau \log k)
]^* \rangle_{HS}.
\end{equation} Since $|\rho| = \rho_0$ is a probability measure,
we can estimate the convolution integral by Cauchy-Schwarz
inequality as
\begin{equation}  \rho_0 \; * \; |\langle {\mathcal D} \ell (k + i \tau \log k)
(I + N(k + i \tau \log k)^{-1} N(k + i \tau \log k),  [ N^{M_0} (k
+ i \tau \log k) \cg {\mathcal S} \ell ^t(k + i \tau \log k) ]^*
\rangle_{HS}|^2. \end{equation} We further apply Schwartz'
inequality to the inner product to obtain the upper bound
\begin{equation}  \begin{array}{l} \rho_0 \; * \; || {\mathcal D} \ell (k + i \tau \log k)
(I + N(k + i \tau \log k)^{-1} N(k + i \tau \log k) ||^2_{HS}\;\;
\\ \\ \cdot  ||  N^{M_0} (k + i \tau \log k) \cg {\mathcal S} \ell ^t(k + i \tau
\log k) \chi(k) ||_{HS}^2. \end{array} \end{equation}

We now separately estimate each factor. The first estimate is
crude but sufficient for our purposes.

\begin{lem} \label{HSUSED}  $|| {\mathcal D} \ell (k + i \tau \log k)
(I + N(k + i \tau \log k)^{-1} N ||_{HS} = O(\frac{k^{1 +
\epsilon}}{\tau}). $ \end{lem}

\begin{proof}  We first use the inequality  $|| A B ||_{HS} \leq ||A|| ||B||_{HS}$ where
$||A||$ is the operator norm to bound
$$|| {\mathcal D} \ell (k + i \tau \log k)
(I + N(k + i \tau \log k)^{-1} N ||_{HS} \leq ||  {\mathcal D}
\ell (k + i \tau \log k) (I + N(k + i \tau \log k)^{-1} ||\; || N
||_{HS}.$$ Here, $|| \cdot ||$ denotes the $L^2(\partial \Omega)
\to L^2(\Omega)$ operator norm.

By Proposition (\ref{N}) (ii), we have  we have $$|| N^2(k + i
\tau \log k)||_{HS} = O(k^{1/2}).
$$ By  (\ref{INTPI})-(\ref{PIF}) the  norm $||  {\mathcal D} \ell (k + i \tau \log k)
(I + N(k + i \tau \log k)^{-1} ||$ is the norm of the Poisson
operator. We claim that
\begin{equation} \label{POISEST} ||PI(k + i \tau)||_{L^2(\partial \Omega) \to
L^2(\Omega)} \leq \tau^{-1} k^{1/2 + \epsilon}. \end{equation}
This estimate is probably crude but it is sufficient for our
purposes. We are not aware of any prior estimates on mapping norms
of the Poisson kernel in the $k$ aspect.

 We first rewrite the Poisson
integral of  (\ref{INTPI}) as
\begin{equation} PI(k, x, q ) = r_{\Omega} X_{\nu}  G_{\Omega}(k, x, y)  \end{equation}
where $X_{\nu}$ is any smooth vector filed on $\Omega$ which
agrees with $\partial_{\nu}$ on $\partial \Omega.$ The operator
$r_{\Omega} X_{\nu}$ is (roughly) of order $3/2$, so we rewrite
the composition as
\begin{equation} \label{PI} PI(k + i \tau, x, q) = r_{\Omega} X_{\nu} \Delta_{\Omega}^{-3/4 - \epsilon}
\Delta_{\Omega}^{3/4 + \epsilon} G_{\Omega}(k + i \tau, x, y).
\end{equation} Here, $\Delta_{\Omega}$ is the Dirichlet Laplacian
and the fractional powers are defined by the method of Seeley
\cite{S1}. Both $X_{\nu}$ and $r_{\Omega}$ operate in the $y$
variable.

To prove (\ref{POISEST}), we use the (well known) fact that
$$\Delta_{\Omega}^{-3/4 - \epsilon}: L^2(\Omega) \to H^{3/2 + \epsilon}(\Omega)$$
is a bounded operator. Furthermore,
$$r_{\Omega} X_{\nu} : H^{3/2 + \epsilon}(\Omega) \to L^2 (\partial \Omega)$$
is bounded for any $\epsilon > 0$. (See \cite{LM}, Theorem 9.4 for
proof of the continuity for $\epsilon > 0$ and Theorem  9.5 for
proof of lack of continuity if $\epsilon = 0$).  So the first
factor is bounded independently of $k + i \tau$.

For the second, we use that $\Delta_{\Omega}^{3/4 + \epsilon}
G_{\Omega}(k + i \tau, x, y)$ is the kernel of
$\Delta_{\Omega}^{3/4 + \epsilon}  (\Delta_{\Omega} + (k + i
\tau)^2)^{-1}$. This is clearly a bounded
 normal operator on $L^2(\Omega)$, so its $L^2(\Omega) \to L^2(\Omega)$ operator  norm is given by
\begin{equation} ||\Delta_{\Omega}^{3/4 + \epsilon}
 (\Delta_{\Omega} + (k + i \tau)^2)^{-1}||_{L^2(\Omega) \to L^2(\Omega)}\;
\; = \max_{j} \frac{\lambda_j^{3/2 + \epsilon}}{|\lambda_j^2 + (k + i \tau)^2|}. \end{equation}
It is elementary to maximise this function, and one finds that the order of magntitude
of the  maximum occurs
when $\lambda_j \sim k$, and it then has the form
\begin{equation}  \frac{k^{3/2 + \epsilon}}{|\tau k|} \sim \frac{k^{1/2 + \epsilon}}{\tau}. \end{equation}
This completes the proof of the Claim and hence of the Lemma.

\end{proof}

We now give the crucial estimate. It explains why we did not need
sharp estimates in the previous step.

\begin{prop}\label{CRUCIAL}  For any $R$, there exists $M_0$ such that:
$$||\rho \; *\;   N^{M_0} (k + i \tau \log k) \cg {\mathcal S} \ell ^t(k + i \tau \log k) \chi(k) ||_{HS}^2
= O(k^{- R}). $$ \end{prop}

\begin{proof}

\end{proof}
We  estimate $||{\mathcal S} \ell (k + i \tau \log k) \chi(k^{-1}
D_x)||$ by a power of $k$ as above, and thus reduce the
proposition to estimating
\begin{equation}| || \rho * N(k + i \tau \log k)^{M_0} \cg ||_{HS}^2 \end{equation}
where the $HS$ norm is now on Hilbert-Schmidt operators on $L^2(\partial \Omega)$.
We write out the Hilbert-Schmidt norm square as the trace:
\begin{equation}\label{TRACE}  Tr \rho * N(k + i \tau \log k)^{M_0} \cg \cg^*  N(k + i \tau \log k)^{* M_0}  \end{equation}
where the trace is on $L^2(\partial \Omega).$

By  Proposition \ref{NPOWERFINAL},
\begin{equation} \label{TRACEFORM} (\ref{TRACE}) = \sum_{\sigma_1, \sigma_2:
Z_{M_0} \to \{0, 1\}} Tr \rho * F_{\sigma_1} \cg \cg^*
F_{\sigma_2}^* \end{equation} plus errors which may be assumed to
be $O(k^{-R})$. We recall that $F_{\sigma}$ is a Fourier integral
operator of order $- |\sigma|$ associated to $\beta^{M_0 -
|\sigma|}$. In addition, the phase of $F_{\sigma}$  has the form
\begin{equation} i k ({\mathcal L}_{\sigma_1} - {\mathcal L}_{\sigma_2})
- \tau \log k  ({\mathcal L}_{\sigma_1} + {\mathcal
L}_{\sigma_2}),
\end{equation}
where   $${\mathcal L}_{\sigma}(q_1, \dots, q_{M_0 - |\sigma|}) =
|q_1 - q_2| + \cdots + |q_{M_0 - |\sigma| - 1} - q_{M_0 -
|\sigma|}|
$$ is the length of an $M_0 - |\sigma|$-link.
The sign difference in the first term reflects the composition
of $N^{M_0}$ with its adjoint, and since the second term comes
from the real part of the phase there is no sign change.

 We now estimate the traces by applying
apply stationary phase. We are only interested in the order of the
trace and not in the coefficients, so we argue qualitatively.  The
terms of the stationary phase expansion of $Tr \rho * F_{\sigma_1}
\cg \cg^* F_{\sigma_2}^* $ correspond  to the fixed points of the
canonical transformation  $\beta^{M_0 - |\sigma_1|} \circ \beta^{-
(M_0 - |\sigma_2|}) = \beta^{- |\sigma_1| +|\sigma_2|} $ which lie
in the support of the cutoffs $\cg$. By assumption, the only
periodic point of $\beta$ in the support of the cutoff is the
period $2$ orbit corresponding to $\gamma$. Thus, other critical
points (which we will call {\it general} critical points) can only
occur when $|\sigma_1| = |\sigma_2|$, which we henceforth write as
$|\sigma$. Equivalently, the general critical points correspond to
a closed path obtained by first following  any $M_0 -
|\sigma|$-link Snell path from a variable point $q_1 \in
\partial \Omega$ to an  endpoint $q_0$, and then reversing along  the same path back to $q_1$.
The general critical points thus form a  non-degenerate critical
manifold parametrized by $(q_1, q_0) \in
\partial \Omega \times \partial \Omega.$

Now consider critical points in the support of $\cg$.  It vanishes
unless the first link points roughly in the direction of the first
link of $\gamma$. Since critical paths are Snell, this forces all
links to point roughly in the directions of links of $\gamma$. It
follows that all links of critical paths (including $\gamma^r$ and
the general ones) have lengths $\sim C L_{\gamma}$ for some
absolute constant $C
> 0$. Note that  $C = 1/2$ for $\gamma$, so this is approximately correct for all links.
Thus,
the imaginary part of the phase introduces the damping factors
\begin{equation}\label{DAMP}  e^{- 2 C \tau \log k (M_0 - |\sigma|)  L_{\gamma}} \end{equation}
into the stationary phase expansion. It follows that
\begin{equation} \label{TRORD} Tr \rho * F_{\sigma_1}
\cg \cg^* F_{\sigma_2}^* \sim   k^{- 2 C \tau  (M_0 - |\sigma|)
L_{\gamma}} k^{- 2 |\sigma|}  k^{- (M_0 - |\sigma|) + 1}
\int_{\partial \Omega \times
\partial \Omega} \alpha \end{equation}
for some smooth density $\alpha$.  We arrived at this order due to
:
\begin{itemize}

\item  the order $-|\sigma|$ of $F_{\sigma}$  (Proposition
\ref{NPOWERFINAL});

\item The damping factor (\ref{DAMP});

\item Application of stationary phase to an $2M_0 - (|\sigma_1| +
|\sigma_2)$-fold integral with a $2$-dimensional non-degenerate
stationary phase manifold.

\end{itemize}

 Given $R$,  we need to choose $(M_0, \tau)$ so that
$$ \begin{array}{l}   - 2 C \tau  (M_0 - |\sigma|)
L_{\gamma} - 2 |\sigma|  - (M_0 - |\sigma|) + 1 \leq -   R \\ \\
\iff (M_0 - |\sigma|) (   2 C \tau L_{\gamma} + 1) +  2 |\sigma|
\geq R + 1.
\end{array} $$ for every $\sigma$. For the case $|\sigma| = M_0$,  it
suffices to have $M_0 \geq (R + 1)/2$. Otherwise, $M_0 - |\sigma|
\geq 1$ and it suffices to pick $\tau \geq \frac{R + 1}{C
L_{\gamma}}.
 $ With these choices of $(M_0, \tau)$ the inequality is true for
 all $\sigma$.

\end{proof}

\end{proof}

\section{Completion of the proof of Theorem \ref{SUM}}

To complete the proof it suffices to determine the trace
asymptotics of the regularized finite sums
 \begin{equation} \label{POT3} \begin{array}{l}
    \sum_{M =
 0}^{M_0} \sum_{\sigma: \{ 1, \dots, M\} \to \{0, 1\}} Tr \rho*
N_{\sigma(1)} \circ N_{\sigma(2)}  \circ \cdots \circ
N_{\sigma(M)} \circ \cg \\ \\ \circ
 {\mathcal S} \ell^{tr}(k + i \tau \log k ) \circ \chi(k)  \circ  {\mathcal D} \ell(k + i \tau \log k).\end{array} \end{equation}
We recall (see \ref{RRHO})  that $\rho * A(k)$ is short for
\begin{equation} \int_{\R} \rho (k - \mu) (\mu + i \tau)  A(\mu ) d \mu.  \end{equation}
\begin{prop} \label{SPM} Let $\gamma$ be a periodic $m$-link reflecting ray, and let $\gamma^r$ be its $r$th iterate.
Let $\hat{\rho} \in C_0^{\infty}(r L_{\gamma} - \epsilon, r
L_{\gamma} + \epsilon)$ be a smooth cutoff, equal to one near $r
L_{\gamma}$ and containing no other lengths in its support.  Then
there exist coefficients $a_{\gamma^{\pm r}, \sigma, j}$ such that
$$\begin{array}{l} Tr\;  \rho*\;
N_{\sigma(1)} \circ N_{\sigma(2)} \circ \cdots \circ N_{\sigma(M)}
\cg
 {\mathcal S} \ell^{tr}(k + i \tau) \circ \chi(k)  \circ  {\mathcal D} \ell(k + i
 \tau) \\ \\
 \sim \left\{ \begin{array}{ll}   e^{ (ik - \tau \log k) r L_{\gamma} }\; k^{ - |\sigma| } \{\sum_{j = 1}^{R}
( a_{\gamma^r, \sigma,  j} + a_{\gamma^{-r} , \sigma,  j}) k^{-j}
+ O(k^{-R} )\}, &  ( M \geq mr ), \\ &  \\  O_R( k^{-R}), & M < m
r\end{array} \right.
\end{array}. $$
Here, $|\sigma| = \{j: \sigma(j) = 0\}$.  In the special case
where $M = m r$ and $\sigma(j) = 1$ for all $j$, we write
$a_{\gamma^{\pm r}, \sigma, j)} = b_{\gamma^{\pm r}; j}$.
 \end{prop}

 The proof consists of a sequence of Lemmas.
We begin by collecting the results of Propositions
\ref{NPOWERFINAL} and \ref{COMP}.

\begin{lem} \label{SOCINT}
 $Tr\;  \rho*\;
N_{\sigma(1)} \circ N_{\sigma(2)} \circ  \circ \cdots \circ
N_{\sigma(M)} \cg
 {\mathcal S} \ell^{tr}(k + i \tau \log k) \circ \chi(k)  \circ(k)  {\mathcal D} \ell(k + i
 \tau)$ is a finite sum of  oscillatory integrals of the form
$$ \begin{array}{l} k^{ (M - |\sigma| + 3)/2} \int_{\R} \int_{\R} \int_{{\bf T}^{M + 1 -
|\sigma|}} e^{i k [(1 - \mu) t + \mu {\mathcal L}_{\sigma}
(q(\phi_1), \dots, q(\phi_{M - |\sigma|} )]} e^{- \tau \log k
{\mathcal
L} (q(\phi_1), \dots, q(\phi_{M - |\sigma|}))} \\ \\
\chi(\overline{q(\phi_1) - q(\phi_2)}, \phi_1) A(k \mu, \phi_1,
\dots, \phi_{M - |\sigma|}) \hat{\rho}(t) dt d\mu d\phi_1 \cdots d
\phi_{M - |\sigma|},
\end{array} $$
where $$\left\{\begin{array}{l}{\mathcal L}_{\sigma} (q(\phi_1),
\dots, q(\phi_{M - |\sigma|} ) = |q(\phi_1) - q(\phi_2)| + \cdots
+ |q(\phi_{M - |\sigma|} - q(\phi_1)|\\ \\
\chi(\overline{q(\phi_1) - q(\phi_2)}, \phi_1) =  \end{array}
\right., $$ and where $A(k, \phi_1, \dots, \phi_{M - |\sigma|})
\in S^{-|\sigma|  }_{\delta}.$
\end{lem}

\begin{proof}

  There are two somewhat
different cases, namely the case where $|\sigma| < M$ and the case
where $|\sigma| = M$.

\subsubsection{Case (i): \; $|\sigma| < M$}

By  Proposition  \ref{NPOWERFINAL},
\begin{equation} \label{NSTRING1}  N_{\sigma(1)} \circ
N_{\sigma(2)} \circ \cdots \circ N_{\sigma(M)}(k, \phi_1, \phi_2)
= F_{\sigma}(k, \phi_1, \phi_2),  \end{equation} where
$F_{\sigma}$ is a semiclassical Fourier integral kernel of the
form \begin{equation} \label{F}  F_{\sigma}(k, \phi_1, \phi_2) =
e^{i (k + i \tau) |q(\phi_{1}) - q(\phi_{2 })|} \chi(k^{1 -
\delta} (\phi_1 - \phi_2)) A_{\sigma}(k, \phi_1, \phi_2),
\end{equation}
with  $A_{\sigma} \in I^{ - |\sigma|}({\bf T}^2 )$. By
Proposition \ref{COMP}, the full composition  \begin{equation}
\label{FULL} N_{\sigma(1)} \circ N_{\sigma(2)} \circ \cdots \circ
N_{\sigma(M)}
 {\mathcal S} \ell^{tr}(k + i \tau) \cg \circ \chi(k)  \circ {\mathcal D} \ell(k + i
 \tau) \end{equation}
 has the form
 \begin{equation} \label{FULLFORM} F_{\sigma}(k) \circ  \cg \circ [D_0 + D_1],\end{equation}
 where the operators $D_0, D_1$ are from Proposition \ref{COMP}.
 $D_1$ puts in an extra factor of $T$ and an amplitude of order
 $-1$.
 All of these operators are semiclassical Fourier integral
 operators with symbols in $S^*_{\delta}$ and therefore they can
 be composed in the standard way.  Therefore we have:
 \begin{equation} \label{FCUT} F_{\sigma}(k) \circ \cg  {\mathcal S} \ell^{tr}(k + i \tau \log k )
  \circ \chi(k)  \circ  {\mathcal D} \ell(k + i \tau \log k)\in I^{- |\sigma| - 1}_{\delta} ({\bf T}^2). \end{equation}
  We then unravel $\rho *$ and recall that the factor $(k - \mu)$
  raises the order by one. Finally we change variables $\mu \to k
  \mu$ which again raises the order by one.

\subsubsection{{\bf Case (ii)}: $|\sigma| = M$}
By  Proposition \ref{NPOWERFINAL}, $N_0^M \circ \cg$ is a $-M$th
order semiclassical pseudodifferential operator, so it suffices to
consider  composition of the form  $A_{-M } \circ \cg [D_0 +
D_1].$ The statement is clear in this case.

\end{proof}

The next step is to show that the stationary phase method applies
to
 oscillatory integrals in $I_{\delta}^{-r}(Y, \Phi)$. This is
 almost obvious, but for the sake of completeness we include the
 proof.

\begin{lem}  \label{MSP}   Let $I_k(a,
\Phi) \in I^{-r}_{\delta}(Y, \Phi)$, and  let $C_{\Phi}$ denote
the set of critical points of $\Phi$, and assume that $C_{\Phi}$
is a non-degenerate critical manifold of codimension $q$. Suppose
that $1/2 < \delta < 1 .$
 Then:
$$I_k(a, \Phi) \sim \left \{ \begin{array}{ll}  k^{-r - q/2}  \sum_{j = 0}^{R} a_{ j} k^{-j} + Rem_{rm, R} (a, \Phi, R)k^{-R - 1},&
C_{\Phi} \not= \emptyset \\ & \\
Rem_{M, R} (a, \Phi, R)k^{-R }, & C_{\Phi} = \emptyset \end{array}
\right.$$ where $a_{ j}$ is a polynomial in the jet of the
amplitude $a$ and phase $\Phi$ at $C_{\Phi},$ and where the
remainder $Rem_{R}(a, \Phi, k) \leq ||a||_{C^{3R}} +
||\Phi||_{C^{3R}}.$
\end{lem}

\begin{proof}

We consider an oscillatory integral
$$I_k(a, \Phi) =  \int_{Y} e^{i k \Phi}  A(k , y)  d y,$$
with $\Phi \in S^0_{\delta}$, with $A(k, y) \in S^0_{\delta}(Y)$.
By assumption, the critical set $C_{\Phi} = \{y:  \nabla_{y} \Phi
= 0 \}$ of the phase  is a non-degenerate critical manifold. We
choose a cutoff $\psi$ supported near $C_{\Phi}$ and write
$$I(a, \Phi) =  I(a \psi, \Phi) + I(a (1 - \psi), \Phi).$$
We now show that if  $a \in S^r_{\delta}$, for some $r$, then $
I(a (1 - \psi), \Phi) = O (k^{-R})$ for any $R > 0$. The implicit
constant is of linear growth in $r$.

In the usual way, we  integrate by parts repeatedly with the
operator:
$$\begin{array}{l} {\mathcal L}_{y} = \frac{1}{k |\nabla_y \Phi|^2} \nabla_y \Phi \cdot \nabla_y, \end{array}$$
that is,  we  apply the transpose
$${\mathcal L}_{y}^t = {\mathcal L}_{y} + \frac{1}{k} \nabla \cdot (\frac{\nabla_y \Phi}{|\nabla \Phi|^2}) $$
to the amplitude. The second term is a scalar multiplication.

We first observe  that   the coefficients of ${\mathcal L}$ belong
to $ S^{0}_{\delta}(Y \times \R \times \R).$ Indeed, by assumption
$L \in S_{\delta}^0(Y)$. In the expression
$$\nabla \cdot (\frac{\nabla \Phi}{|\nabla \Phi|^2}) = \frac{\Delta \Phi }{|\nabla \Phi|^2} +
\frac{\nabla^2 \Phi(\nabla \Phi, \nabla \Phi)}{|\nabla \Phi|^4} =
\frac{\Delta_y L }{|\nabla \Phi|^2} + \frac{\nabla^2_y \L (\nabla
L , \nabla L) }{|\nabla \Phi|^4},$$ it is then  obvious that the
numerator belongs to $S^{0}_{\delta}(Y ).$ Since the denominator
$|\nabla \Phi|^2 =   |\nabla_{y} L|^2 $ is bounded below on
supp$(1 - \psi),$ it follows that the coefficients belong to
$S^{0}_{\delta}(Y).$

We now verify that each partial integration lowers the symbol
order by one unit of $k^{- \delta},$ i.e. that  $(k^{-1} {\mathcal
L}^t )^R  A(k,  y) \in S^{- R \delta}_{\delta}(Y )$. We prove this
by induction on $R$.  As $R \to R + 1$, we apply one of two terms
of $L^t$. We know that each differentiation
 improves the symbol order by one unit of $k^{-\delta}.$ But each coefficient multiplies
by an element of  $S^{0}_{\delta}(Y ),$ hence preserves symbol
order.

The  remainder has the form
\begin{equation} |Rem_{M, R} (a, \Phi, R)| \leq \sup_{(t, \mu, x, \phi)} \max_{\alpha: \alpha \leq 2 R} |D^{\alpha}  \rho(t) (1 - \psi_T)(y) \chi(1 - \mu) A_k((t, \mu, x, \phi)|.\end{equation}
Since $ \rho(t) (1 - \psi_T)(y) \chi(1 - \mu) A_k((t, \mu, x,
\phi) \in S^0(\R \times \R \times Y)$, the right side is $O(k^{(1
- \delta) 2 R}).$ If we choose $\delta$ satisfying $1 > \delta >
1/2 $, then $k^{-R} k^{(1 - \delta) 2 R} = k^{- ( 2 \delta - 1)
R}$ gives a negative power of $k$.  This is sufficient for the
proof of the remainder estimate.

This completes the proof.

\end{proof}

Combining Lemmas \ref{SOCINT}-\ref{MSP}, we obtain

 \begin{equation} \label{EXPANSION}  \begin{array}{l}   k^{ (M - |\sigma| + 3)/2} \int_{\R} \int_{\R} \int_{{\bf T}^{M -
|\sigma| + 1}} e^{i k [(1 - \mu) t + \mu L(q(\phi_1), \dots,
q(\phi_{M - |\sigma|} )]} e^{- \tau \log k {\mathcal L}
(q(\phi_1), \dots, q(\phi_{M - |\sigma| + 1}))} \\ \\
\chi(\overline{q(\phi_1) - q(\phi_2)}, \phi_1) A(k \mu, \phi_1,
\dots, \phi_{M - |\sigma|}) \hat{\rho}(t) dt d\mu
d\phi_1 \cdots d \phi_{M - |\sigma| + 1}\\ \\
 \sim \left\{ \begin{array}{ll}   e^{ (ik - \tau \log k) r L_{\gamma} }\; k^{ - |\sigma| } \{\sum_{j = 1}^{R}
( a_{\gamma^r, \sigma,  j} + a_{\gamma^{-r} , \sigma,  j}) k^{-j}
+ O(k^{-R} )\}, &  ( M \geq mr ), \\ &  \\  O_R( k^{-R}), & M < m
r\end{array} \right.
\end{array}. \end{equation}

Indeed, the critical point set of the phase $\Phi = t (1 - \mu) +
\mu {\mathcal L} $  is given by:
\begin{equation} C_{\Phi} = \{(t, \mu, \phi_1, \dots, \phi_{M -
|\sigma|}): \mu = 1, t = {\mathcal L}(\phi_1, \dots, \phi_{M -
\sigma|}), \overline{q(\phi_1), \dots, q(\phi_{M - |\sigma|}) } =
\gamma^r\}. \end{equation} Clearly, the stationary phase set is
empty if $M < m r$. When $M = rm$ the phase is non-degenerate and
we obtain the result stated in Proposition \ref{SPM} by stationary
phase. This  completes the proof of Theorem \ref{SUM}.

\end{document}